\documentclass{article}

\usepackage{microtype}
\usepackage{graphicx}
\usepackage{subfigure}
\usepackage{booktabs} % for professional tables

\usepackage{hyperref}

\usepackage[accepted]{icml2024}

\usepackage[utf8]{inputenc} % allow utf-8 input
\usepackage[T1]{fontenc}    % use 8-bit T1 fonts
\usepackage{hyperref}       % hyperlinks
\usepackage{url}            % simple URL typesetting
\usepackage{booktabs}       % professional-quality tables
\usepackage{amsfonts}       % blackboard math symbols
\usepackage{nicefrac}       % compact symbols for 1/2, etc.
\usepackage{microtype}      % microtypography
\usepackage{xcolor}         % colors

\usepackage{amsmath}
\usepackage{amssymb}
\usepackage{mathtools}
\usepackage{amsthm}

\usepackage[capitalize,noabbrev]{cleveref}

\theoremstyle{plain}
\newtheorem{theorem}{Theorem}[section]
\newtheorem{proposition}[theorem]{Proposition}
\newtheorem{lemma}[theorem]{Lemma}
\newtheorem{corollary}[theorem]{Corollary}
\theoremstyle{definition}
\newtheorem{definition}[theorem]{Definition}
\newtheorem{assumption}[theorem]{Assumption}
\theoremstyle{remark}

\usepackage{verbatim}
\usepackage[utf8]{inputenc} % allow utf-8 input
\usepackage[T1]{fontenc}    % use 8-bit T1 fonts
\usepackage{hyperref}       % hyperlinks
\usepackage{url}            % simple URL typesetting
\usepackage{booktabs}       % professional-quality tables
\usepackage{nicefrac}       % compact symbols for 1/2, etc.
\usepackage{microtype}      % microtypography
\usepackage{xcolor}         % colors
\usepackage{amsmath}
\usepackage{amssymb}
\usepackage[colorinlistoftodos,bordercolor=orange,backgroundcolor=orange!20,linecolor=orange,textsize=scriptsize]{todonotes}
\usepackage{enumitem}
\usepackage{pifont}
\usepackage{algorithm,algorithmicx,algpseudocode}
\usepackage{wrapfig}

\usepackage{makecell}
\usepackage{caption}
\usepackage{multirow}

\usepackage{colortbl}
\definecolor{bgcolor}{rgb}{0.8,1,1}
\definecolor{bgcolor2}{rgb}{0.8,1,0.8}
\definecolor{niceblue}{rgb}{0.0,0.19,0.56}

\usepackage{hyperref}
\hypersetup{colorlinks,linkcolor={blue},citecolor={niceblue},urlcolor={blue}}

\usepackage{threeparttable}

\usepackage{tcolorbox}
\usepackage{pifont}
\definecolor{mydarkgreen}{RGB}{39,130,67}
\definecolor{mydarkred}{RGB}{192,47,25}
\newcommand{\green}{\color{mydarkgreen}}
\newcommand{\red}{\color{mydarkred}}
\newcommand{\cmark}{{\green\ding{51}}}
\newcommand{\xmark}{{\red\ding{55}}}

\allowdisplaybreaks

\def\R{\mathbb{R}}
\newcommand{\E}{{\mathbb E}}

\def\C{\mathcal C}

\def\X{\mathcal X}

\def\R{\mathbb R}
\def\E{\mathbb E}

\DeclareMathOperator{\gap}{gap}

\newcommand{\cO}{\mathcal{O}}

\newcommand{\EndProof}{\begin{flushright}$\square$\end{flushright}}

\def\<#1,#2>{\langle #1,#2\rangle}

\newcommand{\EE}[1]{\mathbb{E} \left[ #1 \right]}

\icmltitlerunning{Sarah Frank-Wolfe: Methods for Constrained Optimization with Best Rates and Practical Features}

\begin{document}

\twocolumn[
\icmltitle{Sarah Frank-Wolfe: Methods for Constrained Optimization \\ with Best Rates and Practical Features}



\icmlsetsymbol{equal}{*}

\begin{icmlauthorlist}
\icmlauthor{Aleksandr Beznosikov}{1,2,3,4}
\icmlauthor{David Dobre}{5}
\icmlauthor{Gauthier Gidel}{5,6}
\end{icmlauthorlist}

\icmlaffiliation{1}{Innopolis University, Russia}
\icmlaffiliation{2}{Skolkovo Institute of Science and Technology, Russia}
\icmlaffiliation{3}{Moscow Institute of Physics and Technology, Russia}
\icmlaffiliation{4}{Yandex, Russia}
\icmlaffiliation{5}{Universite de' Montreal and Mila, Canada}
\icmlaffiliation{6}{Canada CIFAR AI Chair}

\icmlcorrespondingauthor{Aleksandr Beznosikov}{anbeznosikov@gmail.com}

\icmlkeywords{Machine Learning, ICML}

\vskip 0.3in
]



\printAffiliationsAndNotice{}  % leave blank if no need to mention equal contribution

\begin{abstract}
The Frank-Wolfe (FW) method is a popular approach for solving  optimization problems with structured constraints that arise in machine learning applications. In recent years, stochastic versions of FW have gained popularity, motivated by large datasets for which the computation of the full gradient is prohibitively expensive. In this paper, we present two new variants of the FW algorithms for stochastic finite-sum minimization. Our algorithms have the best convergence guarantees of existing stochastic FW approaches for both convex and non-convex objective functions.
Our methods do not have the issue of permanently collecting large batches, which is common to many stochastic projection-free approaches. Moreover, our second approach does not require either large batches or full deterministic gradients, which is a typical weakness of many techniques for finite-sum problems. The faster theoretical rates of our approaches are confirmed experimentally.
\end{abstract}

\section{Introduction}

Empirical risk minimization is a cornerstone for training supervised machine learning models such as various regressions, support vector machine, and neural networks \citep{shalev2014understanding}. 
We consider a constrained problem of this type:
\begin{equation}
\label{eq:main_problem}
\min_{x \in \X \subset \R^d} f(x) = \frac{1}{n} \sum\limits_{i=1}^n f_i (x)\, .
\end{equation}
% where the constraint set $\X \subset \R^d$ is assumed to be convex. convex, that projecting onto this set is expensive, and also admits a fast linear minimization oracle (LMO). 
The objective function of \eqref{eq:main_problem} has the form of a finite sum. 
Typically, this setting corresponds to the sum of the losses of the model with parameters $x$ applied to a large number of data points, indexed by $i=1,\ldots,n$. 
Because $n$ is large, calculating the full gradient of $f$ is expensive. Therefore, stochastic methods which very rarely resort to calling $\nabla f$ (or avoid it altogether) are of particular importance. 
In this problem setting, we assume the set $\X \subset \R^d$ to be convex, that projecting onto this set is expensive, and that it also admits a fast linear minimization oracle (LMO). 
% Regarding the set $\X \subset \R^d$, we assume that it is convex, that projecting onto this set is expensive, and also admits a fast linear minimization oracle (LMO). 

The study of methods for \eqref{eq:main_problem} that do not require projections has a history of more than half a century. 
Arguably the most popular projection-free method is the Frank-Wolfe (a.k.a, Conditional Gradient algorithm~\citep{frank1956algorithm}). 
This method maintains sparse iterates and only requires a linear minimization oracle that takes into account the specificity of the constraint set $\X$.
In particular, the classical version of the method considers a linear approximation of the function at the current point $x^k$, and minimizes this approximation on the set $\X$:
\begin{align}
\label{eq:fw}
\begin{split}
s^k &= \arg\min_{s \in \mathcal{X}} \< \nabla f(x^k), s - x^k>, \\
x^{k+1} &= x^k + \eta_k (s^k - x^k),
\end{split}
\end{align}
where $\eta_k$ is parameter-free and equal to $\tfrac{2}{k+2}$.

In the last decade, the Frank-Wolfe-type approaches have attracted increasing interest in the machine-learning community because of their good performance on sparse problems, or on problems where the constraints are complex but structured (e.g., various $\ell_p$ balls, trace norms), having applications to submodular optimization~\citep{bachSub2011}, vision~\citep{miech2017learning,bojanowski2014weakly}, and variational inference~\citep{krishnan2015barrier}.
% multitask learning, multi-class classification, matrix learning, recommendation systems, signal
% processing (for details and many more examples, see \cite{hazan2012projection, pmlr-v22-dudik12, pmlr-v28-jaggi13, lacoste2015global, hazan2016variance, negiar2020stochastic}). 

Due to the significant increase in dataset size and complexity within the machine learning community, stochastic algorithms are of great interest and are the focus of this paper. 
% As noted earlier, from a machine learning perspective, stochastic algorithms are of great interest. % Indeed, it is this type of methods that we focus on in this paper. 
In particular, we seek to answer the two following questions:
% \vspace{-0.45cm}
\begin{quote} 
    \textit{
    1. Can we improve upon the convergence rates of existing approaches for \eqref{eq:main_problem}?
    \vspace{0.15cm}
    \\
    2. Can we avoid the computation of both full gradients and large batches of stochastic gradients?}
\end{quote}
% \vspace{-0.25cm}

\clearpage

\renewcommand{\arraystretch}{2}
\begin{table*}[!t]
% \vskip-3pt
    \centering
    % \small
%    \scriptsize
\captionof{table}{Summary of the results on projection free methods for \textbf{stochastic constrained minimization problems}.}
% \vspace{-0.3cm}
    \label{tab:comparison0}   
    % \scriptsize
    \small
\resizebox{\linewidth}{!}{
  \begin{threeparttable}
    \begin{tabular}{ccccccc}
    \hline
    \multirow{2}{*}{\textbf{Reference}} & \multicolumn{2}{c}{\textbf{Convex case complexity}} & \multicolumn{2}{c}{\textbf{Non-convex case complexity}} & \multirow{2}{*}{\textbf{No full gradients?}} & \multirow{2}{*}{\textbf{No big batches?}}
    \\
     & \textbf{SFO} & \textbf{LMO} & \textbf{SFO} & \textbf{LMO}  &  &  \\
    \hline
    \begin{tabular}{@{}c@{}}
    % Frank \& Wolfe 
    \cite{frank1956algorithm}\tnote{{\color{blue}(1)}}  \\[-3mm] 
    % Lacoste-Julien 
    \citep{lacoste2016convergence}\tnote{{\color{blue}(1)}} \end{tabular} & $\cO \left( \frac{n}{\varepsilon}\right)$ & $\cO \left( \frac{1}{\varepsilon}\right)$  &  $\cO \left( \frac{n}{\varepsilon^2}\right)$ & $\cO \left( \frac{1}{\varepsilon^2}\right)$ & \xmark & \xmark
    \\\hline
    % Hazan \& Kale 
    \cite{hazan2012projection} & $\mathcal{O} \left( \frac{1}{\varepsilon^4}\right)$ & $\mathcal{O} \left( \frac{1}{\varepsilon^2}\right)$ &  \multicolumn{2}{c}{\xmark} & \cmark & \xmark
    \\\hline
    % Lan \& Zhou 
    \cite{doi:10.1137/140992382} & $\mathcal{O} \left( \frac{1}{\varepsilon^2}\right)$ & $\mathcal{O} \left( \frac{1}{\varepsilon}\right)$ &  \multicolumn{2}{c}{\xmark}  & \cmark & \xmark
    \\\hline
    % Lan \& Zhou 
    \cite{doi:10.1137/140992382} \tnote{{\color{blue}(1)}} & $\mathcal{O} \left( \frac{n}{\sqrt{\varepsilon}}\right)$ & $\mathcal{O} \left( \frac{1}{\varepsilon}\right)$ &  \multicolumn{2}{c}{\xmark}  & \xmark & \xmark
    \\\hline
    % Reddi et.al. 
    \cite{reddi2016stochastic} Alg. 2 & \multicolumn{2}{c}{\xmark}  &  $\mathcal{ O} \left( \frac{1}{\varepsilon^4}\right)$ & $\mathcal{ O} \left( \frac{1}{\varepsilon^2}\right)$  & \cmark & \xmark
    \\\hline
    % Reddi et.al. 
    \cite{reddi2016stochastic} Alg. 3 & \multicolumn{2}{c}{\xmark}  &  $\mathcal{ O} \left( n + \frac{n^{2/3}}{\varepsilon^2}\right)$ & $\mathcal{ O} \left( \frac{1}{\varepsilon^2}\right)$   & \xmark & \xmark
    \\\hline
    % Reddi et.al. 
    \cite{reddi2016stochastic} Alg. 4 & \multicolumn{2}{c}{\xmark} &  $\mathcal{ O} \left( \frac{n}{\varepsilon^2}\right)$ \tnote{{\color{blue}(2)}} & $\mathcal{ O} \left( \frac{1}{\varepsilon^2}\right)$  & \xmark &  \xmark
    \\\hline
    % Hazan \& Luo 
    \cite{hazan2016variance} & $\mathcal{\tilde O} \left( n + \frac{1}{\varepsilon^2}\right)$ & $\mathcal{O} \left( \frac{1}{\varepsilon}\right)$  &  \multicolumn{2}{c}{\xmark} &  \xmark & \xmark
    \\\hline
    % Qu et.al. 
    \cite{pmlr-v80-qu18a} Alg. 3 & \multicolumn{2}{c}{\xmark}  & $\mathcal{ O} \left( \frac{1}{\varepsilon^4}\right)$ \tnote{{\color{blue}(3)}} & $\mathcal{ O} \left( \frac{1}{\varepsilon^4}\right)$ \tnote{{\color{blue}(3)}} & \cmark & \xmark
    \\\hline
    % Qu et.al. 
    \cite{pmlr-v80-qu18a} Alg. 4 & \multicolumn{2}{c}{\xmark}  & $\mathcal{ O} \left( n + \frac{n^{2/3}}{\varepsilon^2}\right)$ \tnote{{\color{blue}(3)}} & $\mathcal{ O} \left( \frac{1}{\varepsilon^4}\right)$ \tnote{{\color{blue}(3)}} & \xmark & \xmark
    \\\hline
    % Yurtsever et.al. 
    \cite{pmlr-v97-yurtsever19b} & $\mathcal{O} \left( n + \frac{1}{\varepsilon^2}\right)$ & $\mathcal{O} \left( \frac{1}{\varepsilon}\right)$  &  $\mathcal{O}\left( n+ \frac{\sqrt{n}}{\varepsilon^2}\right)$ & $\mathcal{O}\left( \frac{1}{\varepsilon^2}\right)$ & \xmark & \xmark
    \\\hline
    % Gao \& Huang 
    \cite{pmlr-v119-gao20b} Alg. 1 & \multicolumn{2}{c}{\xmark} &  $\mathcal{O}\left( n+ \frac{\sqrt{n}}{\varepsilon^2}\right)$ & $\mathcal{O}\left( \frac{1}{\varepsilon^2}\right)$ &  \xmark & \xmark
    \\\hline
    % Gao \& Huang 
    \cite{pmlr-v119-gao20b} Alg. 2 & \multicolumn{2}{c}{\xmark} &  $\mathcal{O}\left( n + \frac{\sqrt{n}}{\varepsilon^2}\right)$ \tnote{{\color{blue}(3)}} & $\mathcal{O}\left( \frac{1}{\varepsilon^4}\right)$ \tnote{{\color{blue}(3)}} &  \xmark & \xmark
    \\\hline
    % Mokhtari et.al. 
    \cite{mokhtari2020stochastic} & $\cO \left( \frac{1}{\varepsilon^3}\right)$ & $\cO \left( \frac{1}{\varepsilon^3}\right)$ &  \multicolumn{2}{c}{\xmark} & \cmark & \cmark
    \\\hline
    % Zhang et.al. 
    \cite{zhang2020one} & $\cO \left( \frac{1}{\varepsilon^2}\right)$ & $\cO \left( \frac{1}{\varepsilon^2}\right)$ &  \multicolumn{2}{c}{\xmark} & \cmark & \cmark
    \\\hline
    % N{\'e}giar et.al. 
    \cite{negiar2020stochastic}\tnote{{\color{blue}(4)}} & $\cO \left( \frac{n}{\varepsilon}\right)$ \tnote{{\color{blue}(5)}} & $\cO \left( \frac{n}{\varepsilon}\right)$ \tnote{{\color{blue}(5)}} &  \multicolumn{2}{c}{convergence without rate} & \cmark & \cmark
    \\\hline
    % Lu \& Freund 
    \cite{lu2021generalized}\tnote{{\color{blue}(4)}} & $\cO \left( \frac{n}{\varepsilon}\right)$ & $\cO \left( \frac{n}{\varepsilon}\right)$  &  \multicolumn{2}{c}{\xmark} & \cmark & \cmark
    \\\hline
    % Akhtar \& Rajawat 
    \cite{9483167} & $\cO \left( \frac{1}{\varepsilon^2}\right)$  & $\cO \left( \frac{1}{\varepsilon^2}\right)$ &  \multicolumn{2}{c}{\xmark} & \cmark & \cmark
    \\\hline
    % Weber \& Sra 
    \cite{weber2022projection} Alg.2  & $\mathcal{ O} \left( \frac{1}{\varepsilon^4}\right)$ & $\mathcal{ O} \left( \frac{1}{\varepsilon^2}\right)$ & $\mathcal{ O} \left( \frac{1}{\varepsilon^4}\right)$ &  $\mathcal{ O} \left( \frac{1}{\varepsilon^2}\right)$  & \cmark   & \xmark
    \\\hline
    % Weber \& Sra 
    \cite{weber2022projection} Alg.3  & $\mathcal{ O} \left( n + \frac{n^{2/3}}{\varepsilon^2}\right)$ & $\mathcal{ O} \left( \frac{1}{\varepsilon^2}\right)$ & $\mathcal{ O} \left( n + \frac{n^{2/3}}{\varepsilon^2}\right)$ & $\mathcal{ O} \left( \frac{1}{\varepsilon^2}\right)$   & \xmark   & \xmark
    \\\hline
    % Weber \& Sra 
    \cite{weber2022projection} Alg.4  & $\mathcal{ O} \left( \frac{1}{\varepsilon^3}\right)$ & $\mathcal{ O} \left( \frac{1}{\varepsilon^2}\right)$  & $\mathcal{ O} \left( \frac{1}{\varepsilon^3}\right)$    & $\mathcal{ O} \left( \frac{1}{\varepsilon^2}\right)$ & \cmark   & \xmark
    \\\hline
    % Hou et.al. 
    \cite{hou2022distributed} & $\cO \left( \frac{1}{\varepsilon^2}\right)$ & $\cO \left( \frac{1}{\varepsilon^2}\right)$ &  $\cO \left( \exp\left(\frac{1}{\varepsilon}\right)\right)$ & $\cO \left( \exp\left(\frac{1}{\varepsilon}\right)\right)$ & \cmark & \cmark
    \\\hline
    \cellcolor{bgcolor2}{(This paper) Alg. \ref{alg:fw_sarah}} & \cellcolor{bgcolor2}{$\mathcal{\tilde O}\left( n + \frac{\sqrt{n}}{\varepsilon}\right)$} & \cellcolor{bgcolor2}{$\mathcal{\tilde O}\left( \sqrt{n} + \frac{1}{\varepsilon}\right)$}  &  \cellcolor{bgcolor2}{$\mathcal{O}\left( n + \frac{\sqrt{n}}{\varepsilon^2}\right)$} & \cellcolor{bgcolor2}{$\mathcal{O}\left( \frac{1}{\varepsilon^2}\right)$} & \cellcolor{bgcolor2}{\xmark} & \cellcolor{bgcolor2}{\xmark}
    \\\hline
    \cellcolor{bgcolor2}{(This paper) Alg. \ref{alg:fw_sarah}} & \cellcolor{bgcolor2}{$\mathcal{\tilde O}\left( n + \frac{\sqrt{n}}{\varepsilon}\right)$} & \cellcolor{bgcolor2}{$\mathcal{\tilde O}\left( n + \frac{\sqrt{n}}{\varepsilon}\right)$}  &  \cellcolor{bgcolor2}{$\mathcal{O}\left( \frac{n}{\varepsilon^2}\right)$} & \cellcolor{bgcolor2}{$\mathcal{O}\left( \frac{n}{\varepsilon^2}\right)$} & \cellcolor{bgcolor2}{\xmark} & \cellcolor{bgcolor2}{\cmark}
    \\\hline
    \cellcolor{bgcolor2}{(This paper) Alg. \ref{alg:fw_zerosarah}} & \cellcolor{bgcolor2}{$\mathcal{\tilde O}\left( n + \frac{\sqrt{n}}{\varepsilon}\right)$}  & \cellcolor{bgcolor2}{$\mathcal{\tilde O}\left( \sqrt{n} + \frac{1}{\varepsilon}\right)$} &  \cellcolor{bgcolor2}{$\mathcal{O}\left( n + \frac{\sqrt{n}}{\varepsilon^2}\right)$} & \cellcolor{bgcolor2}{$\mathcal{O}\left( \frac{1}{\varepsilon^2}\right)$} & \cellcolor{bgcolor2}{\cmark} & \cellcolor{bgcolor2}{\xmark}
    \\\hline
    \cellcolor{bgcolor2}{(This paper) Alg. \ref{alg:fw_zerosarah} } & \cellcolor{bgcolor2}{$\mathcal{\tilde O}\left( n + \frac{\sqrt{n}}{\varepsilon}\right)$} & \cellcolor{bgcolor2}{$\mathcal{\tilde O}\left( n + \frac{\sqrt{n}}{\varepsilon}\right)$}  &  \cellcolor{bgcolor2}{$\mathcal{O}\left( \frac{n}{\varepsilon^2}\right)$} & \cellcolor{bgcolor2}{$\mathcal{O}\left( \frac{n}{\varepsilon^2}\right)$} & \cellcolor{bgcolor2}{\cmark} & \cellcolor{bgcolor2}{\cmark}
    \\\hline
    %%%%%%%%%%%%%%
    \end{tabular}   
    \begin{tablenotes}
    {\small   
    \item [] \tnote{{\color{blue}(1)}} fully deterministic; \tnote{{\color{blue}(2)}} the authors give a different complexity, but it seems to us that their proof contains an error, we try to correct it (see Appendix \ref{sec:error}); \tnote{{\color{blue}(3)}} in the original papers, the authors give better results, e.g. $\mathcal{O} (\nicefrac{\sqrt{n}}{\varepsilon})$ instead of $\mathcal{O} (\nicefrac{\sqrt{n}}{\varepsilon^2})$, but this results violate the lower bounds (see Table 1 from \cite{pmlr-v139-li21a}]), this is due to the difference in the convergence criterion: in \cite{pmlr-v139-li21a}, the authors use $\|\nabla f\|^2 \sim \varepsilon^2$ and in \cite{pmlr-v80-qu18a, pmlr-v119-gao20b} -- $\|\nabla f\|^2 \sim \varepsilon$; \tnote{{\color{blue}(4)}} only for linear models; \tnote{{\color{blue}(5)}} the authors give a rate in the form $\kappa/\varepsilon$, where $\kappa$ is a special constant, which is equal to $n$ in the worst case. 
    
    {\em Notation:} $\varepsilon$ = accuracy of the solution, $n$ =  size of the dataset, SFO = stochastic first-order oracle, LMO = linear minimization oracle.     
 % \item [] 
    }
\end{tablenotes}    
    \end{threeparttable}
    % }
}
% \vspace{-0.3cm}
\end{table*}

\clearpage

\section{Related Works and Our Contributions} \label{sec:related}

After~\citet{frank1956algorithm} proposed the FW algorithm, many works improved its theory and extended it to special cases~\citep{levitin1966constrained,demianov1970approximate,dunn1978conditional,patriksson1993partial}. 
About ten years ago, \citet{pmlr-v28-jaggi13, lacoste2015global} developed more robust and practical versions of the original FW method, motivated by ML applications with sparsity and structured constraints (see~\citep{braun2022conditional} for a detailed historical survey).

Motivated by applications with large datasets, the theory of stochastic methods for unconstrained (or projection-friendly) optimization problems has built upon the highly successful SGD method \citep{robbins1951stochastic,nemirovski2009robust} to obtain faster methods for finite sum-problems. 
Particularly, many so-called variance-reduced variants of SGD have been proposed, including SAG/SAGA \citep{defazio2014saga,schmidt2017minimizing,qian2019saga}, SVRG \citep{NIPS2013_ac1dd209,allen2016improved,yang2021accelerating},
MISO \citep{mairal2015incremental},
SARAH \citep{nguyen2017sarah,nguyen2021inexact,nguyen2017stochastic,hu2019efficient,li2021zerosarah}, 
SPIDER \citep{fang2018spider}, STORM \citep{cutkosky2019momentum}, PAGE \citep{pmlr-v139-li21a}, and many others.

% \looseness=-1
Extensive research in the theory of deterministic Frank-Wolfe-type methods and stochastic methods for unconstrained problems has led to the development of stochastic versions of projection-free algorithms. 
\citet{hazan2012projection}~proposed an algorithm for online stochastic optimization.
\citet{doi:10.1137/140992382}~developed a projection-free version using sliding.
\citet{hazan2016variance,reddi2016stochastic, pmlr-v80-qu18a, pmlr-v97-yurtsever19b,pmlr-v119-gao20b, pmlr-v89-shen19b} proposed modifications of the Frank-Wolfe method using variance reduction techniques, namely SVRG, SAGA and SPIDER.
\citet{mokhtari2020stochastic, 9483167, hou2022distributed} used the idea of momentum to deal with stochasticity. 
\citet{negiar2020stochastic} and \citet{lu2021generalized} explored stochastic methods for linear predictors. \cite{weber2022projection} extended the results of \citet{reddi2016stochastic} from convex sets to manifolds. 
We summarize and compare the convergence rate of each method in Table~\ref{tab:comparison0}.
Note that for the SAGA-related methods, we report a slightly different result from the one reported by~\citet{reddi2016stochastic} because we believe that their proof contains a slight inaccuracy (see App.~\ref{sec:error} for more details). We also do not include the approach from \citep{pmlr-v89-shen19b} in Table~\ref{tab:comparison0}, since this method uses the hessian of the target function.

% \looseness=-1

% \textbf{Our contributions}

Next, we detail our contributions which can be divided into four parts.
% See Tabel~\ref{tab:comparison0} for a comparison.

$\bullet$ \textbf{The best rates in the convex case.} Our convergence guarantees are better than the classical deterministic method \citep{frank1956algorithm, doi:10.1137/140992382} as well as the stochastic methods from \citep{negiar2020stochastic, lu2021generalized, weber2022projection} in terms of dataset size $n$. Moreover, the theoretical rates of our methods also surpasses the rest existing results from \citep{hazan2012projection, doi:10.1137/140992382, hazan2016variance, pmlr-v97-yurtsever19b,pmlr-v119-gao20b, mokhtari2020stochastic, 9483167, weber2022projection} in terms of the accuracy $\varepsilon$. 

$\bullet$ \textbf{No need for full gradients.} Many stochastic methods, especially for finite-sum problems, require the calculation of some full gradients. 
This makes these techniques less practical because even the infrequent computation of the deterministic gradient can slow down the convergence. 
Some methods for constrained problems also have this disadvantage \citep{reddi2016stochastic, hazan2016variance, pmlr-v80-qu18a, pmlr-v97-yurtsever19b,pmlr-v119-gao20b, weber2022projection}. 
While Algorithm~\ref{alg:fw_sarah} also requires the computation of the full gradient, Algorithm~\ref{alg:fw_zerosarah} removes this issue and uses only stochastic gradients. 
Note that this modification does not affect convergence: Algorithm~\ref{alg:fw_sarah} and Algorithm~\ref{alg:fw_zerosarah} have the same theoretical guarantees. 

$\bullet$ \textbf{Small batches.} 
Many methods that avoid the computation of the full gradient still
use large fixed batch sizes~\citep{hazan2012projection, doi:10.1137/140992382, reddi2016stochastic, pmlr-v80-qu18a, pmlr-v119-gao20b, weber2022projection} or batch sizes that geometrically increase with iteration number \citep{hazan2016variance, pmlr-v97-yurtsever19b}, which, like with the collection of full gradients, is a rather strong limitation on the practical applicability of the method. Conversely, our algorithms are guaranteed to converge with all sizes of batches. Large batches are required to get a better dependence in $n$ for the non-convex case. In particular, it is possible to improve the LMO estimate by a factor of $n$, and the estimate on SFO by a factor of $\sqrt{n}$.  
Methods dealing with fixed small batches are either only analyzed for linear predictors \citep{negiar2020stochastic, lu2021generalized}, or have slower convergence rates than our approach~\citep{mokhtari2020stochastic, 9483167, hou2022distributed}.

% \looseness=-1
$\bullet$ \textbf{Non-convex analysis.} 
We give convergence results not only for the convex problem, but also in the case where the target function $f$ in \eqref{eq:main_problem} is non-convex. In this setting, our oracle complexity results are the first to be non-exponential (in $\epsilon$) with small mini-batches, and are state-of-the-art with large mini-batches.

\section{Notation and Assumptions}

We use $\< x,y > = \sum_{i=1}^d x_i y_i$ to denote the standard inner product of vectors $x,y\in\R^d$, where $x_i$ corresponds to the $i$-th component of $x$ in the standard basis in $\R^n$. With this notation we can introduce the standard $\ell_2$-norm in $\R^d$ in the following way: $\|x\| = \sqrt{\< x, x >}$.  We write $[n] =\{1,2,\dots,n\}$. Calls of the stochastic oracle means computing the gradient $\nabla f_i$ for some $i \in [n]$.

In order to prove convergence results, we state the following standard assumptions on the problem \eqref{eq:main_problem}. The first two assumptions relate to the target function $f$, and the third relates to the constraint set $\X$.

We start with the assumption that the gradients of both the function $f$ and all terms $\{f_i\}^n_{i=1}$ are smooth. This assumption is standard in the optimization literature and widely used in the analysis of Frank-Wolfe-type methods.
\begin{assumption} \label{as:lip}
The function $f: \X \to \R$, is $L$-smooth on $\X$, i.e., there exists a constant $L > 0$ such that 
\begin{equation*}
\| \nabla f(x) - \nabla f(y)\| \leq L \|x-y \|
\,,\quad \forall \,x,y \in \X.
\end{equation*}
Each function $f_i: \X \to \R$, $i \in [n]$, is $L_i$-smooth on $\X$, i.e., there exists a constant $L_i > 0$ such that,
\begin{equation*}
\| \nabla f_i(x) - \nabla f_i(y)\| \leq L_i \|x-y \|
\,,\quad \forall \,x,y \in \X.
\end{equation*}
We also define the constant $\tilde L$ as $\tilde L^2 := \tfrac{1}{n} \sum_{i=1}^n L_i^2$. By convexity of $x\mapsto x^2$, it is easy to prove that $\tilde L \geq L$.
\end{assumption}

The second assumption is the convexity of the function $f$.
\begin{assumption} \label{as:conv}
The function $f: \X \to \R$, is convex, i.e., 
\begin{equation*}
f(x) \geq f(y) + \< \nabla f(y) , x - y>,
\, \quad \forall \,x,y \in \X.
\end{equation*}
\end{assumption}

Note that we consider both convex and non-convex cases of the function $f$. But even if $f$ is convex, we do not additionally assume that the terms $\{f_i\}_{i=1}^n$ are convex, hence in general they can be non-convex. Naturally, it is common to find settings with convex $f$ and convex $f_i$, but formulations with non-convex $f_i$ also arise, e.g., in PCA \citep{garber2015fast, shamir2015stochastic, allen2016improved}.
% This phenomenon is also believed occur in neural networks \citep{allen2016variance, NIPS2013_ac1dd209}.

The next assumption is also typical and found in all works on projection-free methods.

\begin{assumption} \label{as:set}
The set $\X$ is convex and compact with a diameter $D$, i.e., for any $x,y \in \X$,
\begin{equation*}
\| x - y \| \leq D.
\end{equation*}
\end{assumption}

%%%%%%%%%%%%%%%%%%%%%%%%%%%%%%%%%%%%%%%%%%
%%%%%%%%%%%%%%%%%%%%%%%%%%%%%%%%%%%%%%%%%% 

\section{Main Part}

In this section, we present two new algorithms and their convergence guarantees. 

\subsection{State-of-the-art complexity with Sarah Frank-Wolfe} \label{sec:sarah}

\begin{algorithm*}[!h]
	\caption{Sarah Frank-Wolfe}
	\label{alg:fw_sarah}
	\hspace*{\algorithmicindent} {\bf Parameters:} step sizes $\{\eta_k\}_{k\geq 0}$, probability $p$, batch size $b$;\\
	\hspace*{\algorithmicindent} {\bf Initialization:} choose  $x^0 \in \mathcal{X}$; $g^0 = \nabla f(x^0)$;
	\begin{algorithmic}[1]
		\For{$k=0,1,2,\ldots K-1$}
            \State Compute $s^k = \arg\min_{s \in \mathcal{X}} \< g^k, s>$; \label{alg1:line2}
            \State Update $x^{k+1} = x^k + \eta_k (s^k - x^k)$ with $\eta_k$ \label{alg1:line3}
            \State Generate batch $S_k$ with size $b$;
		    \State Update 
            $g^{k+1} = \begin{cases}
            \nabla f(x^{k+1}), & \text{with probability } p, \\
            g^k + \frac{1}{b} \sum\limits_{i \in S_k} [\nabla f_{i}(x^{k+1}) -  \nabla f_{i}(x^{k})], & \text{with probability } 1-p,
            \end{cases}$; \label{alg1:line4}
		\EndFor
	\end{algorithmic}
\end{algorithm*}
 
Previously, \citet{reddi2016stochastic, hazan2016variance, weber2022projection} proposed to modify the classical Frank-Wolfe algorithm \eqref{eq:fw} using the SVRG technique \citep{NIPS2013_ac1dd209}. 
% The authors of \cite{reddi2016stochastic, hazan2016variance, weber2022projection} propose to modify the classical Frank-Wolfe algorithm \eqref{eq:fw} using the SVRG technique \cite{NIPS2013_ac1dd209}. 
The essence of these modifications is to change the deterministic gradient in the Conditional Gradient method to some stochastic gradient $g^k$, e.g., calculated according to the SVRG approach:
\begin{equation*}
    g^k = \nabla f_{i_k} (x^k) - \nabla f_{i_k} (w^k) + \nabla f (w^k),
\end{equation*}
where $i_k$ is randomly generated from $[n]$, and $w^k$ is rarely taken equal to $x^k$, much more often equal to $w^{k-1}$. Therefore, when we update $w^k = x^k$, we sometimes need to consider the full deterministic gradient. The update rule for $w^k$ can be deterministic (as in the original version) or randomized, known as the loopless approach \citep{kovalev2020don}. 
Meanwhile, there are other variance-reduced methods, such as SARAH \citep{nguyen2017sarah}:
\begin{equation}
    \label{eq:sarah}
    g^k = \nabla f_{i_k} (x^k) - \nabla f_{i_k} (x^{k-1}) + g^{k-1}, 
\end{equation}
where $i_k$ is also randomly generated from $[n]$, and $g^k$ is rarely taken equal to $\nabla f(x^k)$ rather than \eqref{eq:sarah}. As noted in the original paper on SARAH, this method has better convergence guarantees and smoother convergence paths with less oscillations than SVRG, making SARAH preferred in both theory and practice. As a result, the SARAH update is also used in the Conditional Gradient method for non-convex problems \citep{pmlr-v97-yurtsever19b, pmlr-v119-gao20b, weber2022projection}. In these works, the authors call \eqref{eq:sarah} the SPIDER technique. We also use SARAH as a base for Algorithm \ref{alg:fw_sarah}, but unlike \citet{pmlr-v97-yurtsever19b, pmlr-v119-gao20b, weber2022projection}, its' loopless version \citep{pmlr-v139-li21a}.
First we give the convergence of Algorithm~\ref{alg:fw_sarah} in the convex case. 
\begin{theorem} \label{th:main1}
Let $\{x^k\}_{k\geq0}$ denote the iterates of Algorithm~\ref{alg:fw_sarah} for solving problem \eqref{eq:main_problem}, which satisfies Assumptions~\ref{as:lip}--\ref{as:set}. Let $x^*$ be the minimizer of $f$. Then for any $K$ one can choose $\{ \eta_k \}_{k \geq 0}$ as follows:
\begin{align*}
    &\text{if}~~  K \leq \frac{2}{p}, && \eta_k = \frac{p}{2}, \\
    &\text{if}~~  K > \frac{2}{p} ~~  \text{ and } ~~  k < \left\lceil \frac{K}{2} \right\rceil, &&  \eta_k = \frac{p}{2}, \\
    &\text{if}~~  K > \frac{2}{p} ~~  \text{ and } ~~  k \geq \left\lceil \frac{K}{2} \right\rceil, && \eta_k = \frac{2}{(\nicefrac{4}{p} + k - \lceil \nicefrac{K}{2}\rceil)}.
\end{align*}
For this choice of $\eta_k$, we have the following convergence:
\begin{align*}
    \EE{f(x^{K}) - f(x^*)} =
    \mathcal{O}\Bigg(& \frac{f(x^{0}) - f(x^*) }{p} \exp\left(-\frac{K p}{4} \right)
    \\
    &+ \left[1 
    +  \frac{\tilde L}{L}\sqrt{\frac{1-p}{p b}}\right] \frac{LD^2}{K}\Bigg).
\end{align*}
\end{theorem}

% \begin{proof}[Proof Sketch]
% \color{red}We could add a proof sketch if we have space.
% \end{proof}

See the full proof in Section \ref{sec:proof1}. 
% In the convergence assessment from Theorem \ref{th:main1}, the parameter $p$ is not picked, let us do it. 
We highlight an important detail that the convergence is proved not only in terms of $f(x^K) - f(x^*)$, but also for the Lyapunov function, which includes additionally $\| g^K - \nabla f(x^K) \|^2$. Therefore, Theorem \ref{th:main1} gives guarantees that $\| g^K - \nabla f(x^K) \|^2 \sim \frac{1}{K}$ and hence $g^K$ becomes a good approximation of the full gradient $\nabla f(x^K)$.
It is also worth pointing out that the results of Theorem \ref{th:main1} depends on $\| g^0 - \nabla f(x^0) \|^2$ in the general case (see Section \ref{sec:proof1}), but due to our initialization of Algorithm \ref{alg:fw_sarah}, it is equal to zero.

To choose $p$, one can note that for each iteration, we on average compute the stochastic gradient $(pn + (1-p)\cdot 2b)$ times: with probability $p$ we need the full gradient, with probability $(1-p)$ -- a batch of size $b$ in two points $x^{k+1}$ and $x^k$. 
If we take $p$ close to $1$, the guarantees in Theorem \ref{th:main1} gives faster convergence, but the oracle complexity per iteration increases. For example, if we take $p = 1$, we simply obtain a deterministic method, and the estimates for convergence and the number of gradient calculations reproduce the results for the classical Frank-Wolfe method. On the other hand, if $p$ tends to $0$, the number of stochastic gradient calls per iteration decreases, but the iterative convergence rate drops.
It is optimal to choose $p$ based on the condition: $pn = 2(1-p)b$, i.e. $p = \tfrac{2b}{n+2b}$. From Theorem \ref{th:main1}  we can also obtain an estimate on the required number of linear minimizations (LMO complexity). It is equal to the number of iterations of Algorithm \ref{alg:fw_sarah}. 
Then, the following corollary holds.
\begin{corollary} \label{cor:main1}
Under the conditions of Theorem \ref{th:main1}, Algorithm \ref{alg:fw_sarah} with $p = \tfrac{2b}{n+2b}$ achieves an $\varepsilon$ suboptimality in expectation with 
% $\bullet$ the stochastic oracle calls:
\begin{align*}
    &\mathcal{O}\left( \frac{n}{b} \log \frac{1}{\varepsilon}   + \left[1
    +  \frac{\tilde L \sqrt{n}}{b L} \right] \frac{LD^2}{\varepsilon}\right) \; 
    \text{LMO calls, and}\\
    &\mathcal{O}\left( n \log \frac{1}{\varepsilon}   + \left[b 
    +  \frac{\tilde L \sqrt{n}}{L} \right] \frac{LD^2}{\varepsilon}\right) \,
    \text{stoch. oracle calls.} 
\end{align*}
\end{corollary}
For any $b \leq \tfrac{\tilde L \sqrt{n}}{L}$, the estimate of the number of calls for the stochastic oracles does not change. Given that $\tilde L \geq L$, the smallest batch size $b = 1$ is appropriate for us. 
In this setting, the required number of the stochastic gradient computations is $\mathcal{\tilde O} \left(n + \tfrac{\sqrt{n} \tilde L D^2}{\varepsilon}\right)$. This result is the best in the literature around stochastic projection-free methods, especially since it does not require using large batches (see Table \ref{tab:comparison0}). Meanwhile, in the general case (not necessarily projection-free), these estimates can be explicitly improved up to $\mathcal{\tilde O} \left(n + \sqrt{\tfrac{n \tilde L D^2}{\varepsilon}}\right)$ \cite{allen2018katyusha}. Whether this is possible in the case of Frank-Wolfe type methods is an attractive question for future consideration.
For $p = \tfrac{2b}{n+2b}$ and $b =  \sqrt{n}$,  the LMO complexity is $\mathcal{\tilde O} \left(\sqrt{n} + \tfrac{\tilde L D^2}{\varepsilon}\right)$, this result is optimal up to an additional factor $\sqrt{n}$ (see Section 2.1.2 from \citep{braun2022conditional}). With $b = 1$, the LMO complexity equals $\mathcal{\tilde O} \left(n + \tfrac{\sqrt{n} \tilde L D^2}{\varepsilon}\right)$. Note that for many practical examples, the LMO complexity is not the computational bottleneck since the solution of linear minimization problems has a closed-form solution (see, e.g., Algorithm 2 from \citep{bellet2015distributed}).
It is also important to notice that, based on the above choices for $\eta_k$, $p$, and $b$, both our method  and the original Frank-Wolfe are independent of the objective function parameters (e.g., $L$ or $\tilde L$).

Next, we prove the convergence of Algorithm \ref{alg:fw_sarah} for the non-convex objective function $f$. We use the \textit{Frank-Wolfe gap} function \citep{pmlr-v28-jaggi13} as a criterion for convergence: 
\begin{equation}
    \label{eq:gap}
    \gap (y) = \max_{x \in \X} \, \langle \nabla f(y), y - x \rangle.
\end{equation}
Such a criterion is standard in the analysis of algorithms for the constrained problems with non-convex functions \citep{lacoste2016convergence, reddi2016stochastic}. It is easy to check that $\gap(y) \geq 0$ for any $y \in \X$. Moreover, a point $y \in X$ is stationary for \eqref{eq:main_problem} if and only if $\gap(y) = 0$. \cite{lacoste2016convergence} notes that the Frank-Wolfe gap is a meaningful measure of non-stationarity, and also an affine invariant generalization of the more standard convergence criterion $\|\nabla f (y) \|$ that is used for unconstrained non-convex problems. Then the following theorem is valid.

\begin{algorithm*}[!t]
	\caption{Saga Sarah Frank-Wolfe}
	\label{alg:fw_zerosarah}
	\hspace*{\algorithmicindent} {\bf Parameters:} step sizes $\{\eta_k\}_{k \geq 0}$; momentum $\lambda$; batch size $b$;\\
	\hspace*{\algorithmicindent} {\bf Initialization:} choose  $x^0 \in \mathcal{X}$; $g^0 = \nabla f(x^0)$ or $g^0 = \nabla f_{i_0} (x^0)$; $y^0_i = \nabla f_i (x^{0})$ or $y^0_i = 0$;
 % for $i \in [n]$;
	\begin{algorithmic}[1]
		\For{$k=0,1,2,\ldots K-1$}
            \State Compute $s^k = \arg\min_{s \in \mathcal{X}} \< g^k, s>$; \label{alg2:line2}
            \State Update $x^{k+1} = x^k + \eta_k (s^k - x^k)$ with $\eta_k$ \label{alg2:line3}
            \State Generate batch $S_k$ with size $b$
		\State Update \text{\small{$g^{k+1} = \frac{1}{b} \sum\limits_{i \in S_k}[\nabla f_{i}(x^{k+1}) - \nabla f_{i}(x^{k})] + (1 - \lambda) g^k +        \lambda\left( \frac{1}{b} \sum\limits_{i \in S_k}[\nabla f_{i} (x^k) - y^k_{i}] + \frac{1}{n}\sum\limits_{j=1}^n y_j^k\right)$}};              \label{alg2:line4}
            \State Update $y^{k+1}_i = 
            \begin{cases}
            \nabla f_i (x^{k+1}), & i \in S_k, \\
            y^k_i, & i \notin S_k,
            \end{cases}$; \label{alg2:line5}
		\EndFor
	\end{algorithmic}
\end{algorithm*}

\begin{theorem} \label{th:main2}
Let $\{x^k\}_{k\geq0}$ denote the iterates of Algorithm \ref{alg:fw_sarah} for solving problem \eqref{eq:main_problem}, which satisfies Assumptions \ref{as:lip},\ref{as:set}. Let $x^*$ be the global (may be not unique) minimizer of $f$. Then, if we choose $\eta_k = \frac{1}{\sqrt{K}}$, we have the following convergence:
\begin{align*}
    \E &\left[\min_{0\leq k \leq K-1}\gap(x^k)\right]
    \\&=
    \mathcal{O} \Bigg(\frac{f(x^0) - f(x^*)}{\sqrt{K}} 
    + \frac{LD^2}{\sqrt{K}} \left[1 
    +  \frac{\tilde L}{L} \sqrt{\frac{(1-p)}{pb}}\right]\Bigg).
\end{align*}
\end{theorem}
See the proof in Section \ref{sec:proof2}. 
In this case, the optimal choice of the parameter $p$ is the same as in Corollary \ref{cor:main1} of Theorem \ref{th:main1}. 
\begin{corollary} \label{cor:main2}
Under the conditions of Theorem \ref{th:main2},  Algorithm \ref{alg:fw_sarah} with $p = \tfrac{2b}{n+2b}$ achieves an $\varepsilon$ suboptimality in expectation with  
\begin{align*}
    &\mathcal{O}\left( \left[\frac{h^0}{\varepsilon}\right]^2 + \left[\frac{LD^2}{\varepsilon}\right]^2 \left[1 +  \frac{\tilde L^2 n}{L^2 b^2}\right]\right) \; 
    \text{LMO calls, and}\\
    &\mathcal{O}\left( b\left[\frac{h^0}{\varepsilon}\right]^2 + \left[\frac{LD^2}{\varepsilon}\right]^2 \left[b +  \frac{\tilde L^2 n}{L^2 b}\right]\right)\,  \text{stoch. oracle calls,} 
\end{align*}
where $h^0 := f(x^0) - f(x^*).$
\end{corollary}
First, we substitute $b = 1$ in the previous result.
This gives us an oracle complexity of $\mathcal{O}\left( \tfrac{n (\tilde L^2 + L^2) D^4 + (h^0)^2}{\varepsilon^2}\right)$, which corresponds to the Frank-Wolfe complexity. 
If we choose the batch size $b = \tfrac{\sqrt{n} \tilde L}{L}$ that minimizes the expression $\left[b +  \tfrac{\tilde L^2 n}{L^2 b}\right]$, then we need $\mathcal{O}\left( n+ \tfrac{\sqrt{n} \tilde L [L^2 D^4 + (h^0)^2]}{L \varepsilon^2}\right)$ stochastic gradient calls. 
If one wishes to avoid using $\tilde L$ and $L$ constants when selecting $b$, it is possible to take $b = \sqrt{n}$, and the complexity then becomes $\mathcal{O}\left(n + \tfrac{\sqrt{n}[(\tilde L^2 + L^2) D^4 + (h^0)^2]}{\varepsilon^2}\right)$. 
As noted earlier, both of these estimates are the best result of the projection-free methods for the non-convex setup (see Table \ref{tab:comparison0}). Moreover, they are optimal according to the lower estimates from Table 1 of \cite{pmlr-v139-li21a}. 
Additionally, with $p = \tfrac{2b}{n+2b}$ and $b =1$, the LMO complexity is equal to $\mathcal{O}\left( \tfrac{n (\tilde L^2 + L^2) D^4 + (h^0)^2}{\varepsilon^2}\right)$, and with $b = \sqrt{n}$, it is $\mathcal{O}\left( \tfrac{ (\tilde L^2 + L^2) D^4 + (h^0)^2}{\varepsilon^2}\right)$. In terms of lower bounds, only the already mentioned LMO estimates $\mathcal{O}\left( \tfrac{1}{\varepsilon} \right)$ for the convex setting are known. One important detail to note is that in both the convex (see the discussion after Corollary \ref{cor:main1}) and non-convex cases, LMO estimates is better when the batches are chosen quite large. In Section \ref{sec:batches}, we give reasoning about this.
%%%%%%%%%%%%%%%%%%%%%%%%%%%%%%%%%%%%%%%%%%%%%%%%%%%%
%%%%%%%%%%%%%%%%%%%%%%%%%%%%%%%%%%%%%%%%%%%%%%%%%%%%

\subsection{Avoiding full gradient computations with Saga Sarah Frank-Wolfe} \label{sec:saga_sarah}

The idea of Algorithm \ref{alg:fw_zerosarah} is to use a combination \citep{li2021zerosarah} of the SARAH and SAGA \citep{defazio2014saga} approaches. Both SAGA and SARAH are some of the main variance-reduced methods for finite-sum minimization problems. An important feature of SAGA is that it does not use full gradient calculations, but it has worse convergence guarantees than SARAH (see, e.g., Table 2 in \citep{nguyen2017sarah}). The synergy of SARAH and SAGA brings together the strengths of both methods. 

The essence of the SAGA method is similar to SVRG, but where SVRG collects the full gradients at some reference points, SAGA instead maintains a "sliding" version of the full gradient. The gradient $\nabla f(w)$ at the reference point $w$ may become obsolete after a small number of iterations, but in the course of the algorithm we compute newer stochastic gradients for some $i \in [n]$, and one can leverage them to calculate a more recent approximation of the full gradient. To do this, SAGA introduces additional vectors $\{y\}_{i=1}^n$; each such $y_i$ keeps the latest version of the gradient $\nabla f_i$ (implemented in line \ref{alg2:line5}).
The $\tfrac{1}{n}\sum_{j=1}^n y_j^k$ term is the aforementioned approximation of the full gradient. 
The calculation of $g^{k}_{\text{SAGA}}$ is $g^{k}_{\text{SAGA}} = \big( \tfrac{1}{b} \sum_{i \in S_k}[\nabla f_{i} (x^k) - y^k_{i}] + \tfrac{1}{n}\sum_{j=1}^n y_j^k\big)$ (similarly to SVRG and SARAH).
Line \ref{alg2:line4} provides a combination of $g^{k}_{\text{SARAH}}$ and $g^{k}_{\text{SAGA}}$.

Recall that the average number of the stochastic oracle calls per iteration of Algorithm \ref{alg:fw_sarah} is $(pn + (1-p)\cdot 2b)$, and that Algorithm \ref{alg:fw_zerosarah} requires $2b$ computations of the stochastic gradients each iteration.
In lines \ref{alg2:line4} and \ref{alg2:line5}, one need to calculate $\nabla f_i(x^{k+1})$, $\nabla f_i(x^{k})$ for $i \in S_k$. Therefore, if $b \leq \tfrac{n}{2}$, then for any $p \neq 0$, the complexity of one iteration of Algorithm \ref{alg:fw_zerosarah} is better than that of Algorithm~\ref{alg:fw_sarah}.

In summary, Algorithm \ref{alg:fw_zerosarah} does not collect full gradients and has a bit better iteration complexity, but is required to use $n$ extra vectors $\{ y_i \}_{i=1}^n$, requiring an additional $\mathcal{O}(nd)$ memory cost compared to Algorithm~\ref{alg:fw_sarah}. Once can note that the methods from \citep{negiar2020stochastic, lu2021generalized} also use an extra memory size of $\mathcal{O}(nd)$.

For the convex target function $f$, Algorithm \ref{alg:fw_zerosarah} satisfies the following convergence theorem.
\begin{theorem} \label{th:main3}
Let $\{x^k\}_{k\geq0}$ denote the iterates of Algorithm \ref{alg:fw_zerosarah} for solving problem \eqref{eq:main_problem}, which satisfies Assumptions \ref{as:lip}--\ref{as:set}. Let $x^*$ be the minimizer of $f$. Then for any $K$ one can choose $\{ \eta_k \}_{k \geq 0}$ as follows:
\begin{align*}
    &\text{if}~~  K \leq \frac{4n}{b}, && \eta_k = \frac{b}{4n}, \\
    &\text{if}~~  K > \frac{4n}{b} ~~  \text{ and } ~~  k < \left\lceil \frac{K}{2} \right\rceil, &&  \eta_k = \frac{b}{4n}, \\
    &\text{if}~~  K > \frac{4n}{b} ~~  \text{ and } ~~  k \geq \left\lceil \frac{K}{2} \right\rceil, && \eta_k = \frac{2}{(\nicefrac{8n}{b} + k - \lceil \nicefrac{K}{2}\rceil)},
\end{align*}
and $\lambda = \tfrac{b}{2n}$.
For this choice of $\eta_k$ and $\lambda$, we have the following convergence:
\begin{align*}
    \hspace{-0.2cm}\EE{f(x^{K}) - f(x^*)} 
    =
    \cO\Bigg(& \frac{n \left[f(x^{0}) - f(x^*) \right] }{b}  \exp\left(-\frac{bK}{8n} \right)
    \\
    &+ \left[1 
    +  \frac{\tilde L \sqrt{n}}{L b}\right] \frac{LD^2}{K}\Bigg).
\end{align*}
\end{theorem}
See the proof in Section \ref{sec:proof3}. As in the case of Theorem \ref{th:main1}, here the proof is also obtained in terms of the Lyapunov function containing $\| g^K - \nabla f(x^K)\|^2$, we can guarantee that $g^K$ tends to $\nabla f(x^K)$. The guarantees of Theorem \ref{th:main3} depend on $\| g^0 - \nabla f(x^0)\|^2$ and $\sum_{i=1}^n\| y^0_i - \nabla f_i (x^0)\|$, but because of initialization we put them equal to $0$ again.
Since we do not need to choose $p$ for Algorithm \ref{alg:fw_zerosarah} we proceed directly to the corollary on the oracle complexity. 
\begin{corollary}
Under the conditions of Theorem \ref{th:main3}, Algorithm \ref{alg:fw_zerosarah} achieves an $\varepsilon$ suboptimality in expectation with
\begin{align*}
    &\mathcal{O}\left( \frac{n}{b} \log \frac{1}{\varepsilon}   + \left[1
    +  \frac{\tilde L \sqrt{n}}{b L} \right] \frac{LD^2}{\varepsilon}\right) \; 
    \text{LMO calls, and}\\
    &\mathcal{O}\left( n \log \frac{1}{\varepsilon}   + \left[b 
    +  \frac{\tilde L \sqrt{n}}{L} \right] \frac{LD^2}{\varepsilon}\right) \,  \text{stoch. oracle calls.} 
\end{align*}
\end{corollary}
This result is exactly the same as Corollary \ref{cor:main1}, so we obtain the same conclusions for choosing the size of $b$ as that Corollary \ref{cor:main1}. In particular, in this case with  $b=1$, the method also have oracle complexity $\mathcal{\tilde O} \left(n + \tfrac{\sqrt{n} \tilde L D^2}{\varepsilon}\right)$ -- the best among the works on stochastic projection-free methods. The findings on the LMO complexity is also consistent with Theorem \ref{th:main1}. The reasoning around optimality is also consistent with that given after Corollary \ref{cor:main1}. Note also that Algorithm \ref{alg:fw_zerosarah}, just the same as Algorithm \ref{alg:fw_sarah} and the classical Frank-Wolfe method, is free of the target function's parameters.

In the following theorem for the non-convex case of $f$, as in Theorem \ref{th:main2}, we use \eqref{eq:gap} to estimate convergence.

\begin{theorem} \label{th:main4}
Let $\{x^k\}_{k\geq0}$ denote the iterates of Algorithm \ref{alg:fw_zerosarah} for solving problem \eqref{eq:main_problem}, which satisfies Assumptions \ref{as:lip},\ref{as:set}. Let $x^*$ be the global (may be not unique) minimizer of $f$ on $\X$. Then, if we choose $\eta_k = \frac{1}{\sqrt{K}}$ and $\lambda = \tfrac{b}{2n}$, we have the following convergence:
\begin{align*}
    \E \left[\min_{0\leq k \leq K-1}\gap(x^k)\right]
    =
    \mathcal{O}\Bigg(&
    \frac{f(x^0) - f(x^*)}{\sqrt{K}}
    \\
    &+ \frac{LD^2}{\sqrt{K}} \left[1 
    +  \frac{\tilde L \sqrt{n}}{Lb}\right]\Bigg).
\end{align*}
\end{theorem}
See the proof in Section \ref{sec:proof4}.
\begin{corollary}
Under the conditions of Theorem \ref{th:main4}, Algorithm \ref{alg:fw_zerosarah} achieves an $\varepsilon$ suboptimality in expectation with
\begin{align*}
    &\mathcal{O}\left( \left[\frac{h^0}{\varepsilon}\right]^2 + \left[\frac{LD^2}{\varepsilon}\right]^2 \left[1 +  \frac{\tilde L^2 n}{L^2 b^2}\right]\right)\; 
    \text{LMO calls, and}\\
    &\mathcal{O}\left( b\left[\frac{h^0}{\varepsilon}\right]^2 + \left[\frac{LD^2}{\varepsilon}\right]^2 \left[b +  \frac{\tilde L^2 n}{L^2 b}\right]\right)\,  \text{stoch. oracle calls, } 
\end{align*}
where $h^0 := f(x^0) - f(x^*).$
\end{corollary}
We once again have the same results as for Algorithm \ref{alg:fw_sarah} in Corollary \ref{cor:main2}, resulting in the same analysis for choosing the batch sizes and for the LMO complexity.

%%%%%%%%%%%%%%%%%%%%%%%%%%%%%%%%%%%%%%%%
%%%%%%%%%%%%%%%%%%%%%%%%%%%%%%%%%%%%%%%%
%%%%%%%%%%%%%%%%%%%%%%%%%%%%%%%%%%%%%%%%
%%%%%%%%%%%%%%%%%%%%%%%%%%%%%%%%%%%%%%%%

\section{Experiments}

\begin{figure*}[h!]
% \vspace{-0.3cm}
\begin{minipage}{0.24\textwidth}
  \centering
\includegraphics[width =  \textwidth ]{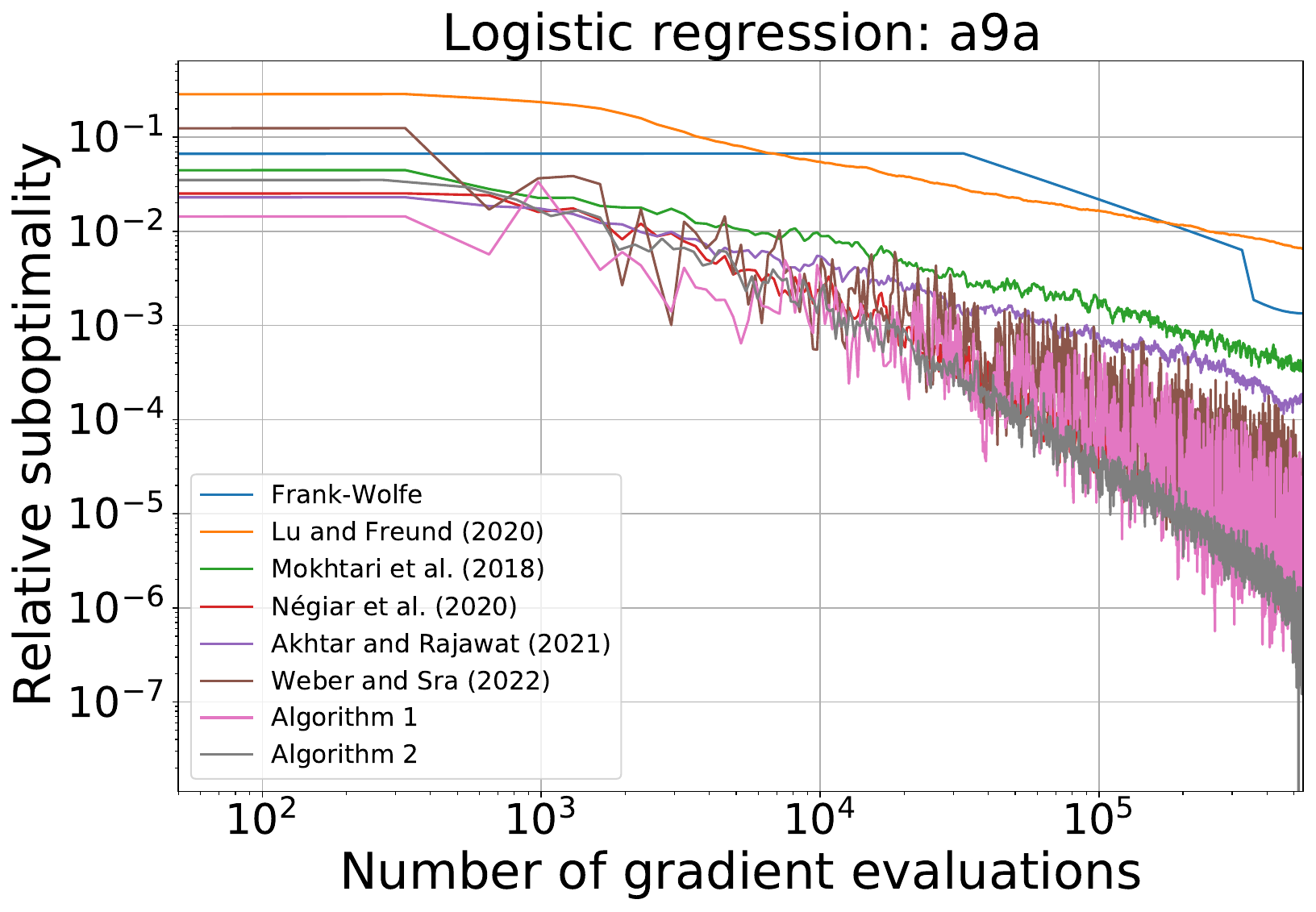}
\end{minipage}%
\begin{minipage}{0.24\textwidth}
  \centering
\includegraphics[width =  \textwidth ]{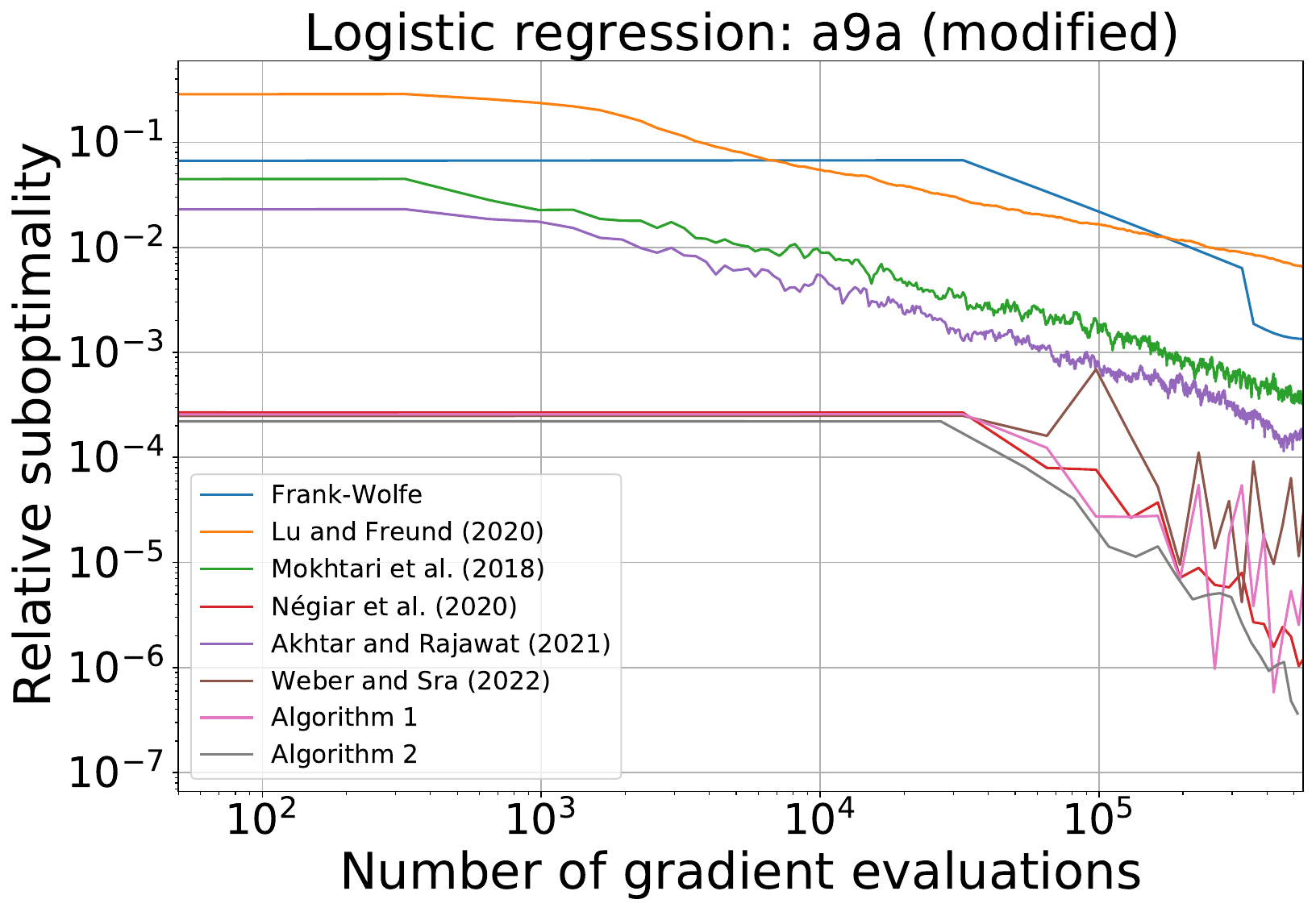}
\end{minipage}%
\begin{minipage}{0.24\textwidth}
  \centering
\includegraphics[width =  1\textwidth ]{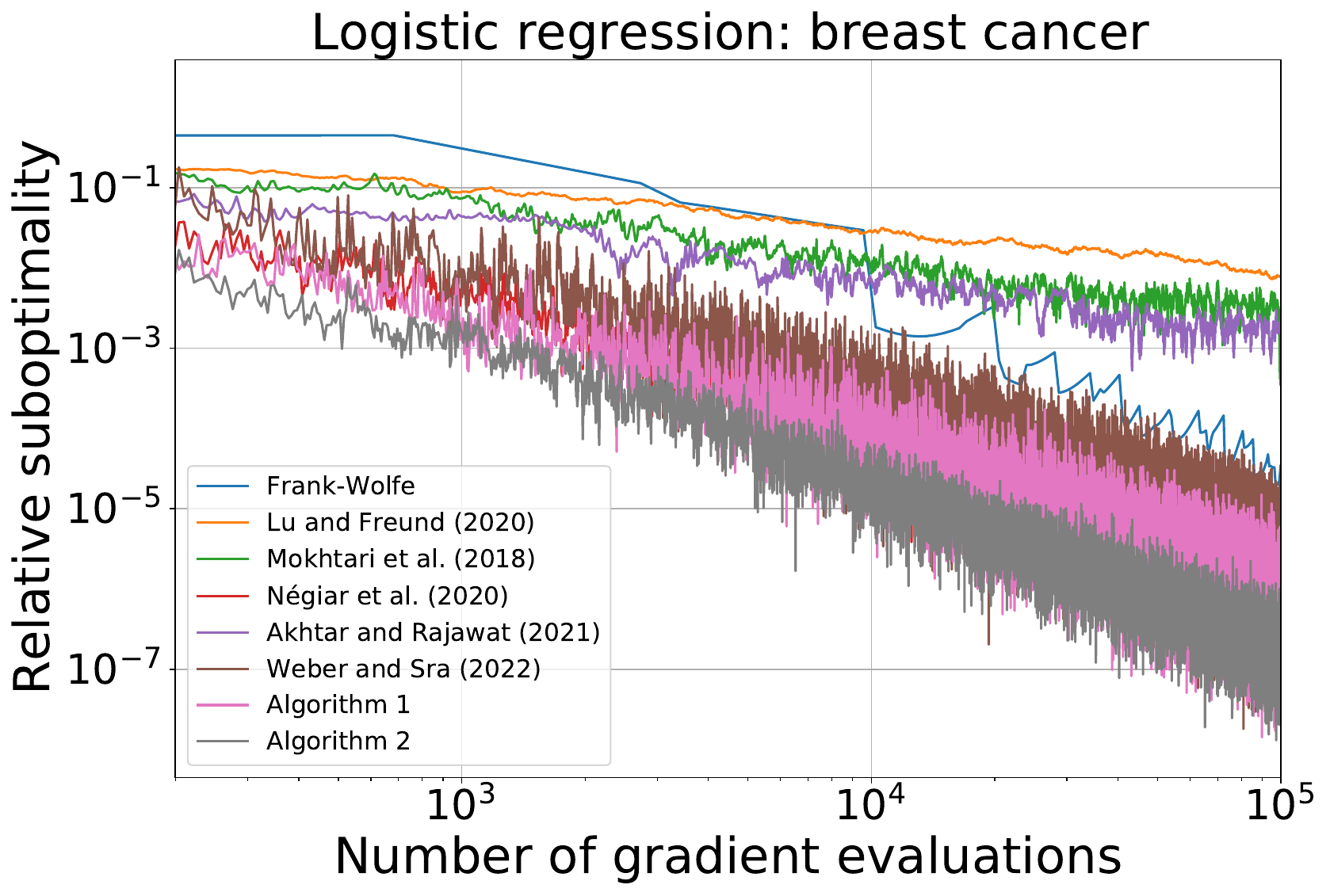}
\end{minipage}%
\begin{minipage}{0.24\textwidth}
  \centering
\includegraphics[width =  1\textwidth ]{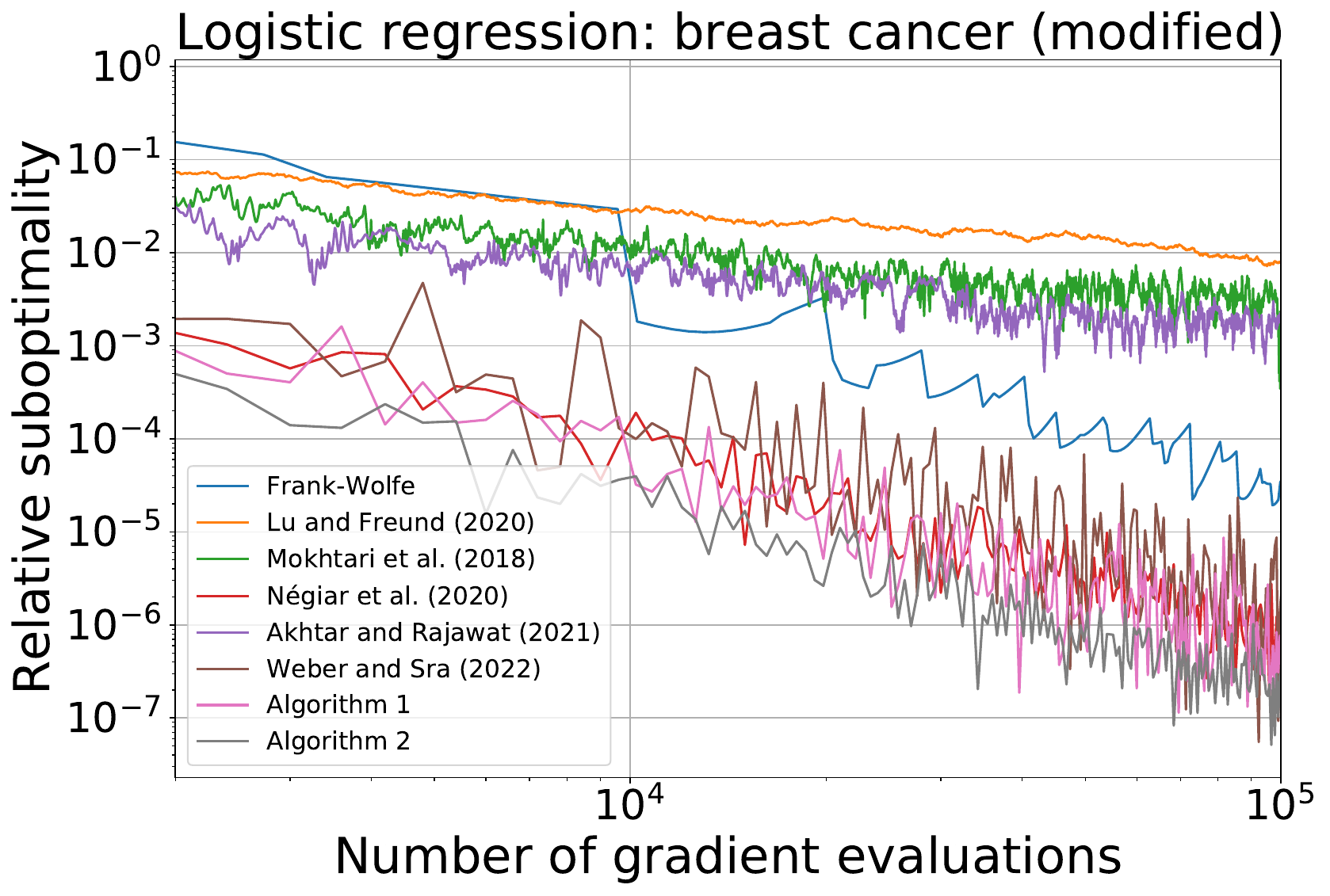}
\end{minipage}%
\\
\begin{minipage}{0.49\textwidth}
  \centering
(a) \texttt{a9a}
\end{minipage}%
\begin{minipage}{0.49\textwidth}
\centering
  (b) \texttt{breast cancer}
\end{minipage}%
\\
\begin{minipage}{0.24\textwidth}
  \centering
\includegraphics[width =  1\textwidth ]{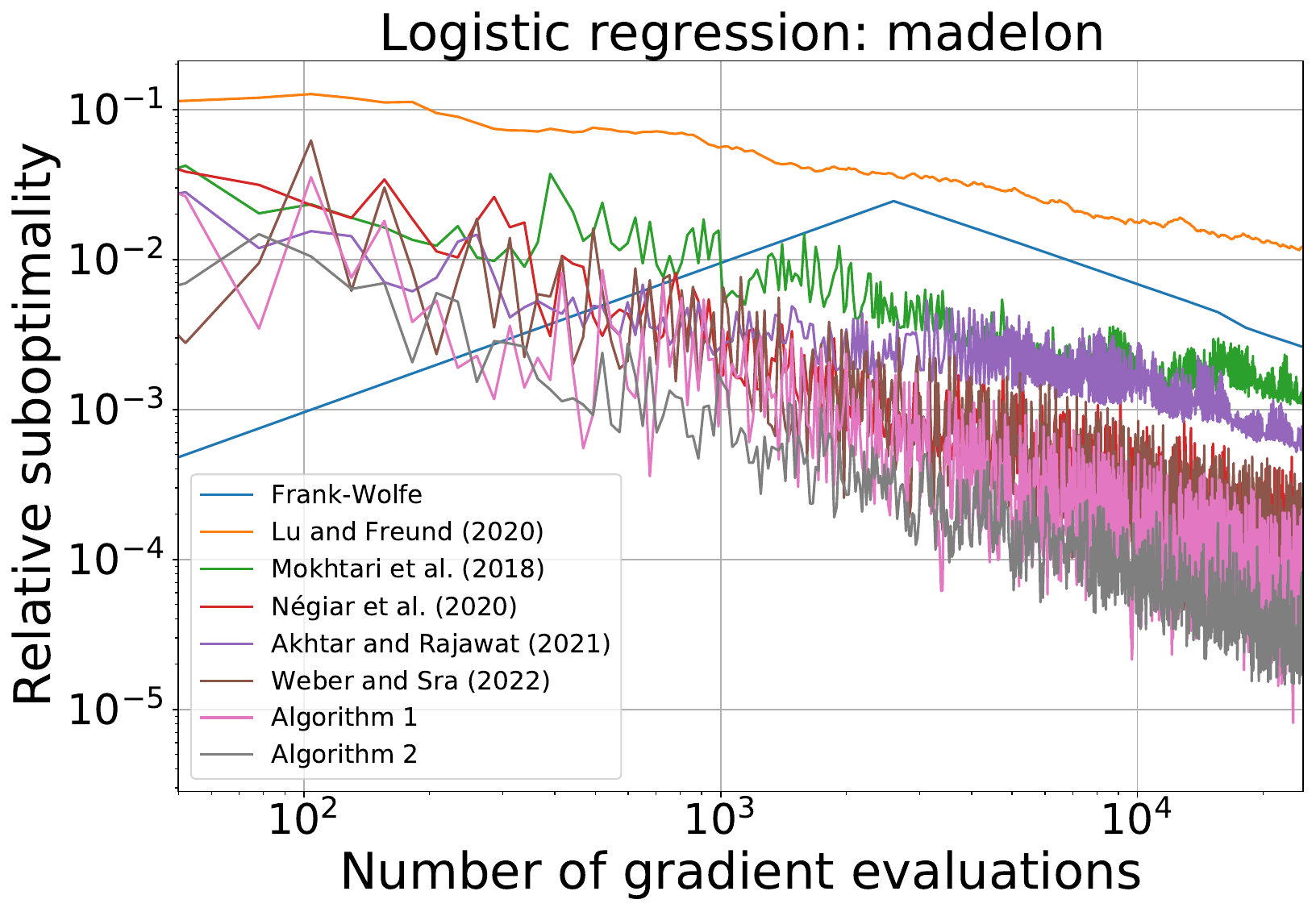}
\end{minipage}%
\begin{minipage}{0.24\textwidth}
  \centering
\includegraphics[width =  1\textwidth ]{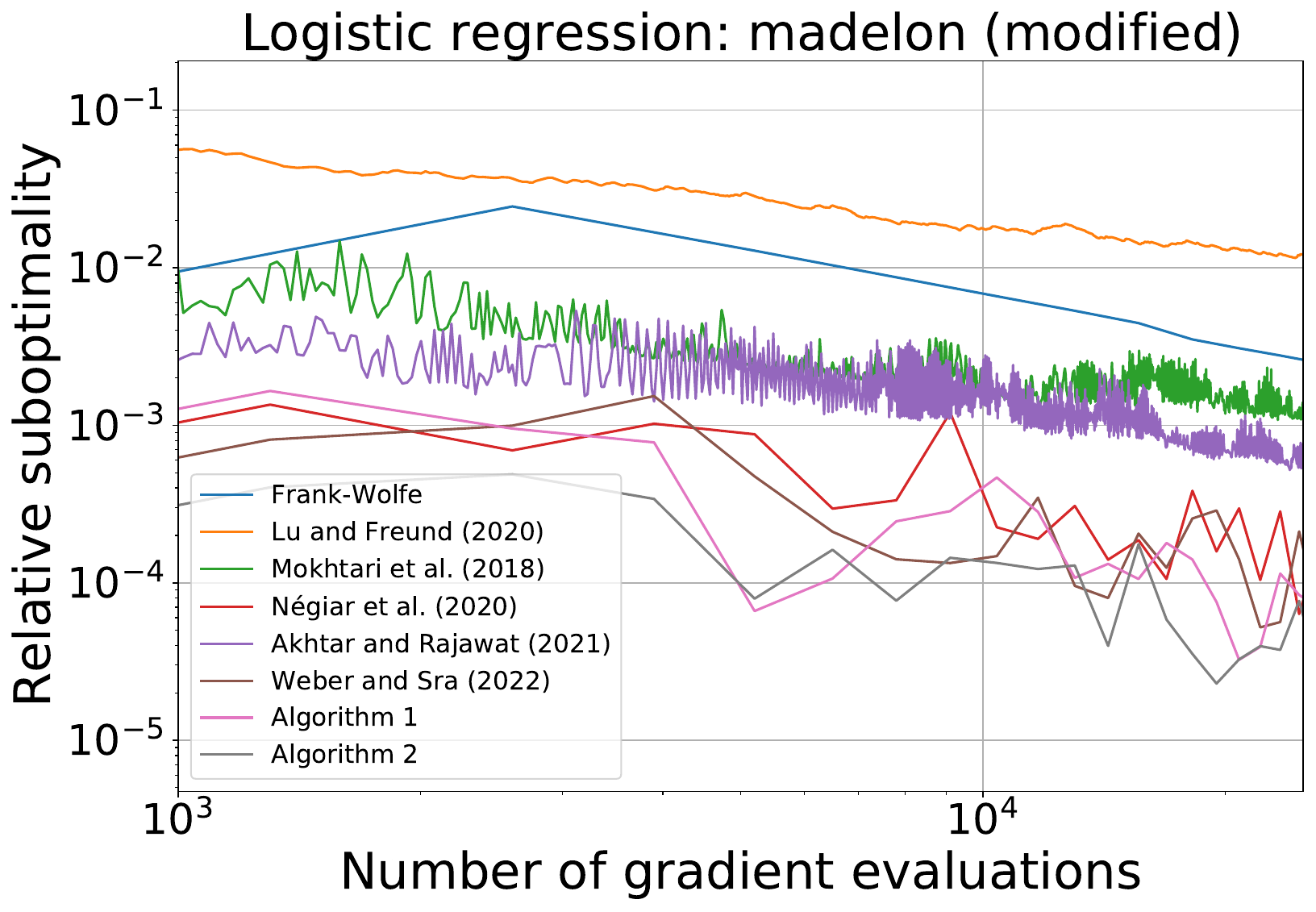}
\end{minipage}%
\begin{minipage}{0.24\textwidth}
  \centering
\includegraphics[width =  1\textwidth ]{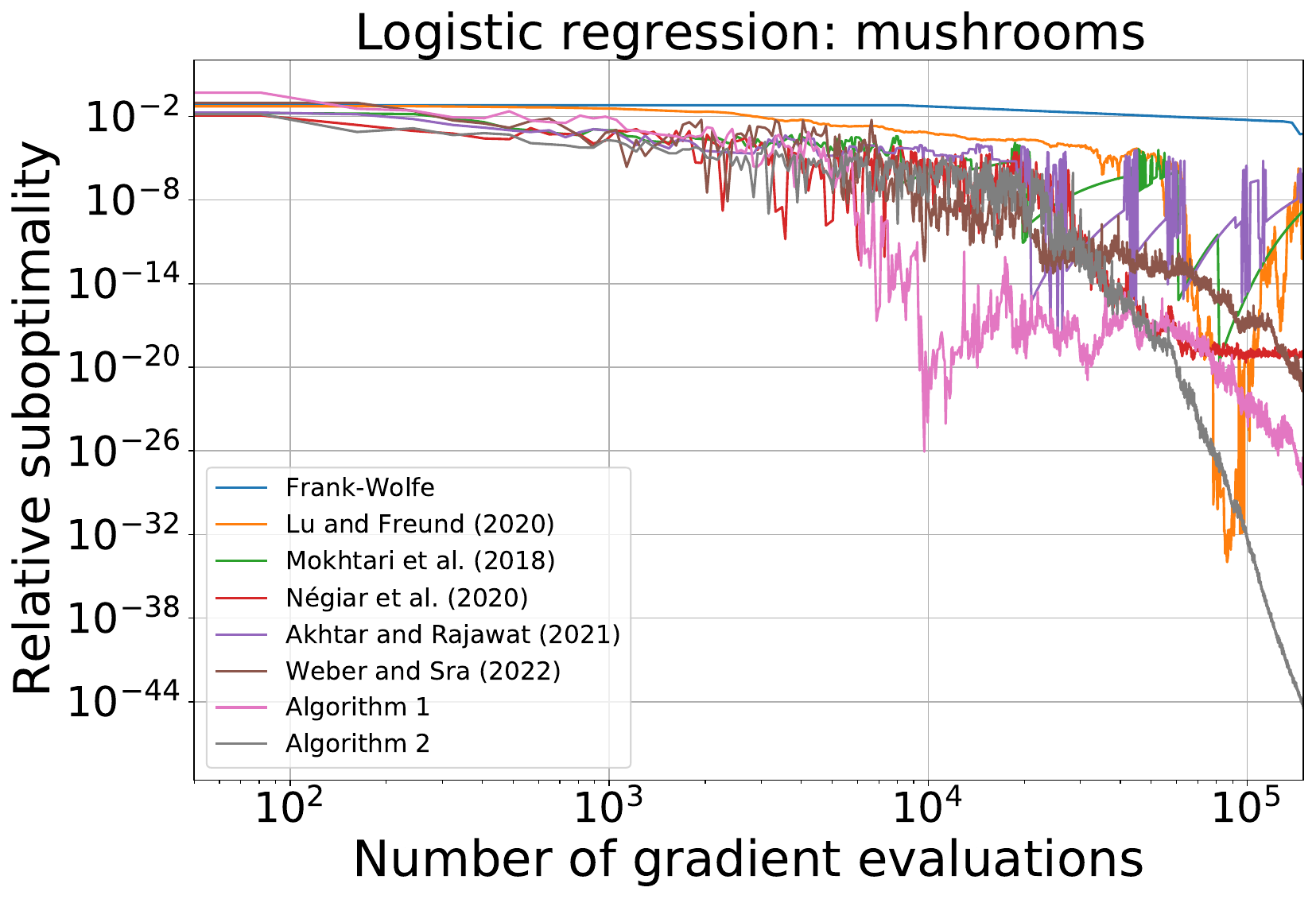}
\end{minipage}%
\begin{minipage}{0.24\textwidth}
  \centering
\includegraphics[width =  1\textwidth ]{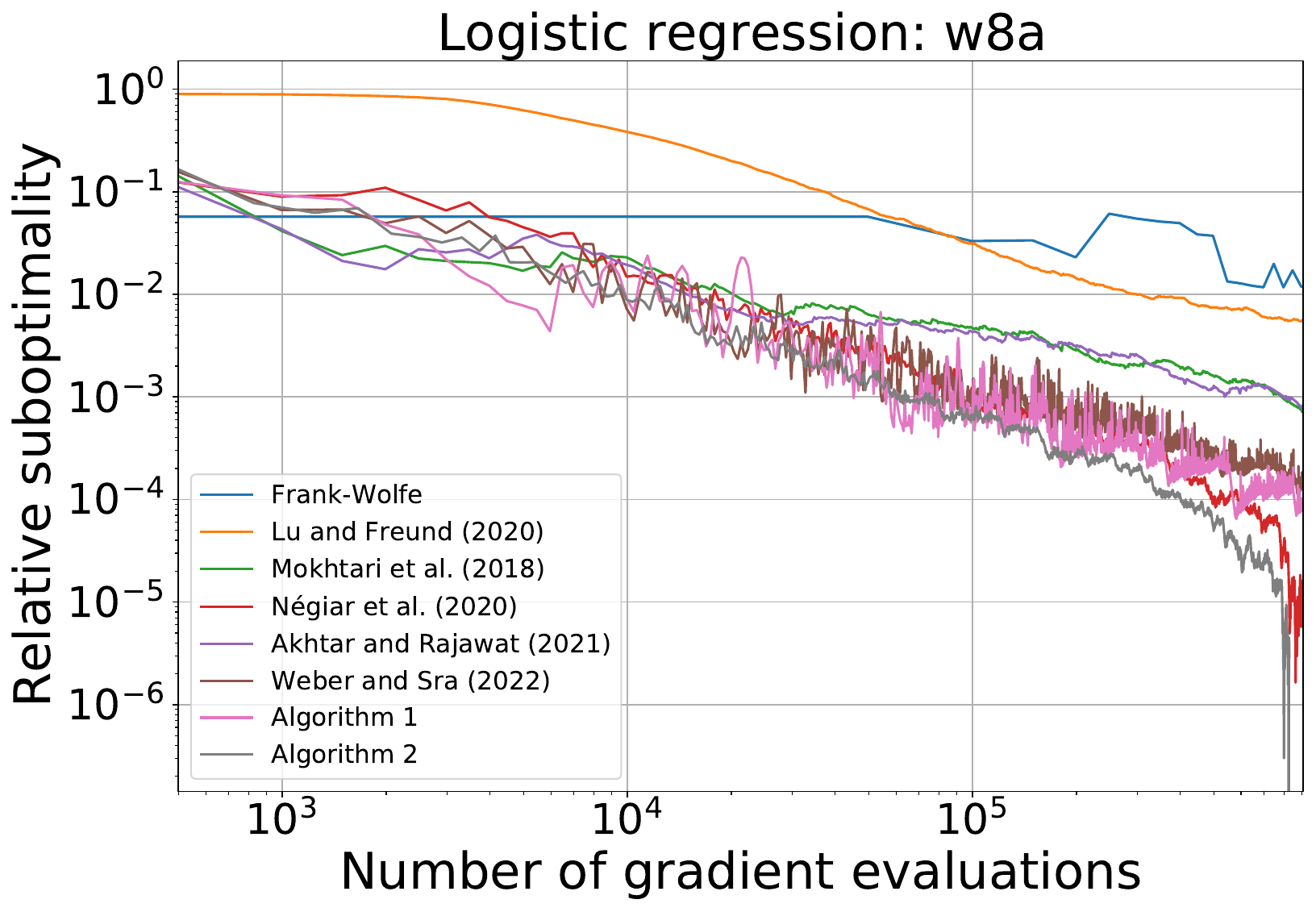}
\end{minipage}%
\\
\begin{minipage}{0.49\textwidth}
\centering
  (c) \texttt{madelon}
\end{minipage}%
\begin{minipage}{0.24\textwidth}
\centering
  (d) \texttt{mushrooms}
\end{minipage}%
\begin{minipage}{0.24\textwidth}
\centering
  (e) \texttt{w8a}
\end{minipage}%
\vspace{-0.3cm}
\caption{Comparison of state-of-the-art projection free methods with small batches for \eqref{eq;ls}. The comparison is made on the real datasets from LibSVM. The criterion is the number of full gradients computations. In the modified plots (the right plots in the first three lines), we left only every 100th point for \citep{negiar2020stochastic}, \citep{weber2022projection}, Algorithm \ref{alg:fw_sarah} and Algorithm \ref{alg:fw_zerosarah}.}
    \label{fig:min}
\end{figure*}
% \vspace{-0.4cm}
\begin{figure*}[h!]
\begin{minipage}{0.24\textwidth}
  \centering
\includegraphics[width =  \textwidth ]{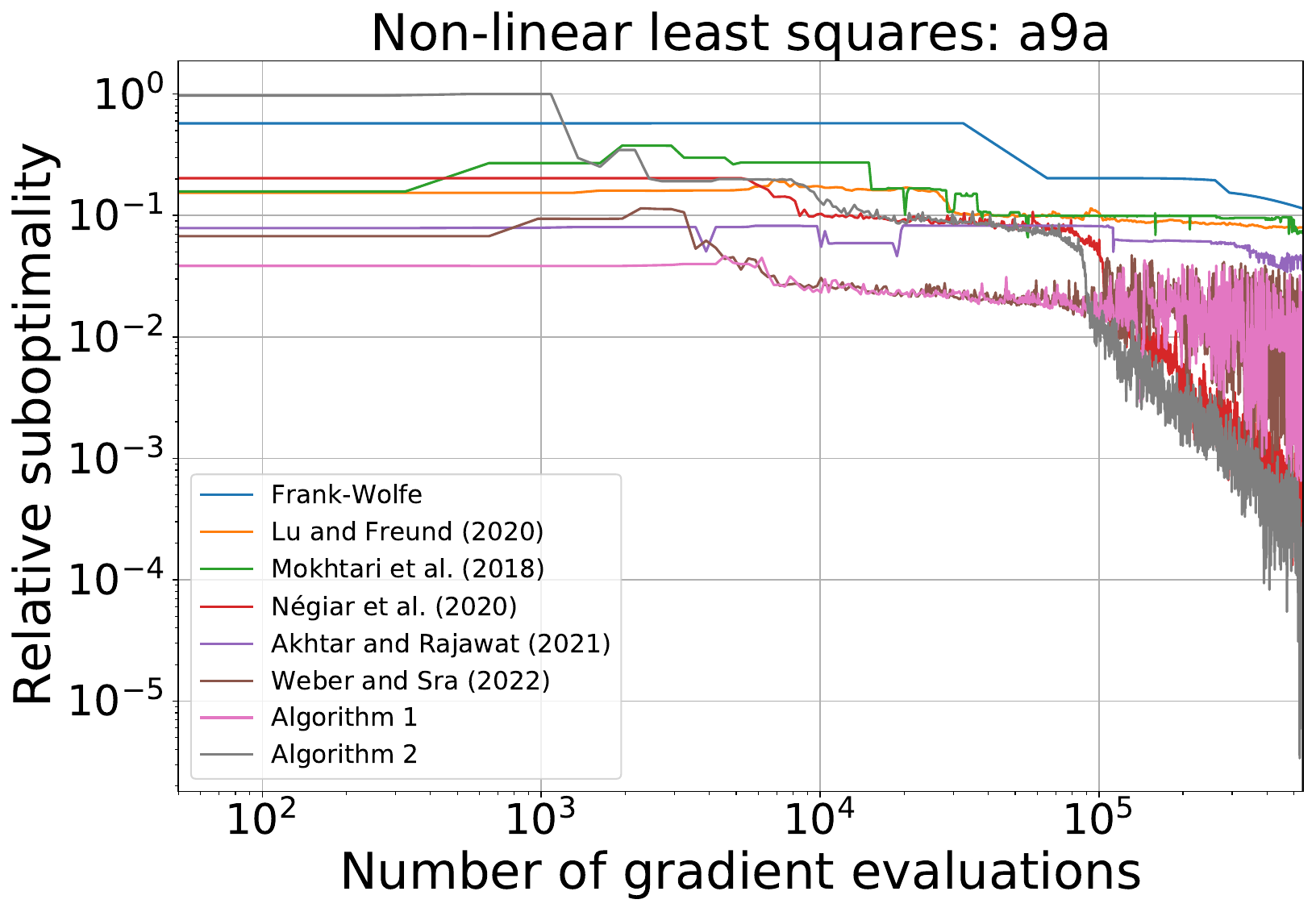}
\end{minipage}%
\begin{minipage}{0.24\textwidth}
  \centering
\includegraphics[width =  \textwidth ]{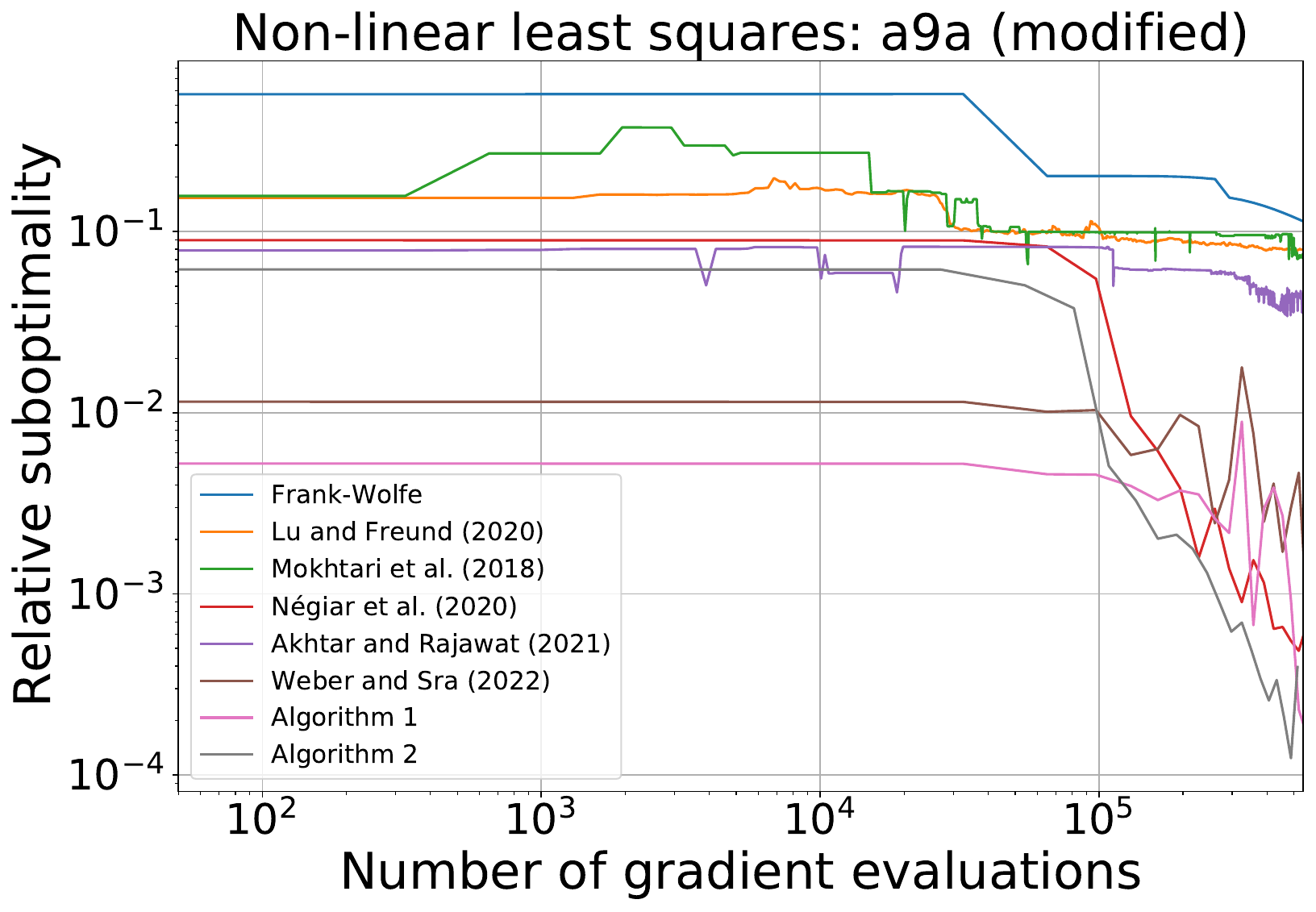}
\end{minipage}%
\begin{minipage}{0.24\textwidth}
  \centering
\includegraphics[width =  1\textwidth ]{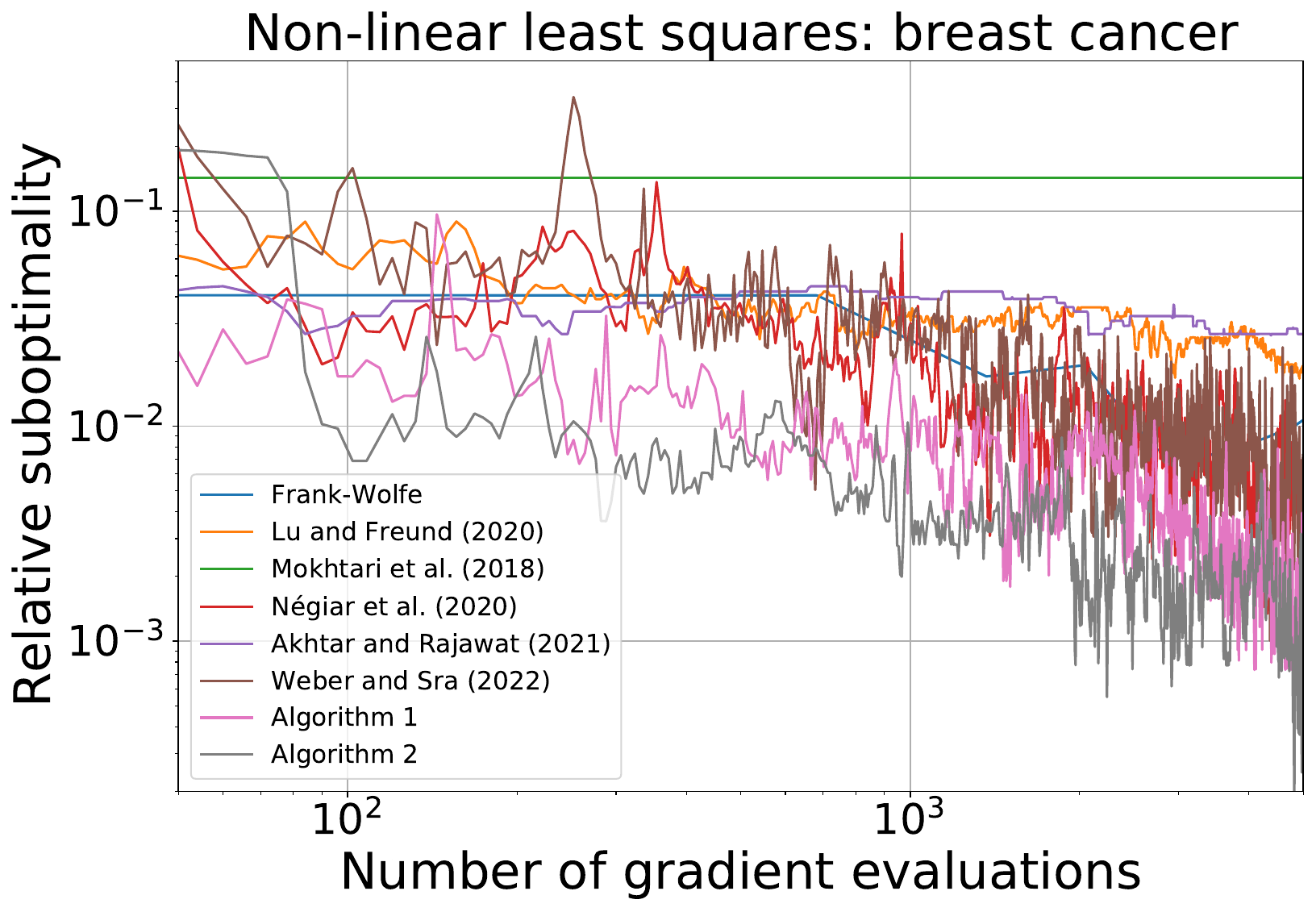}
\end{minipage}%
\begin{minipage}{0.24\textwidth}
  \centering
\includegraphics[width =  1\textwidth ]{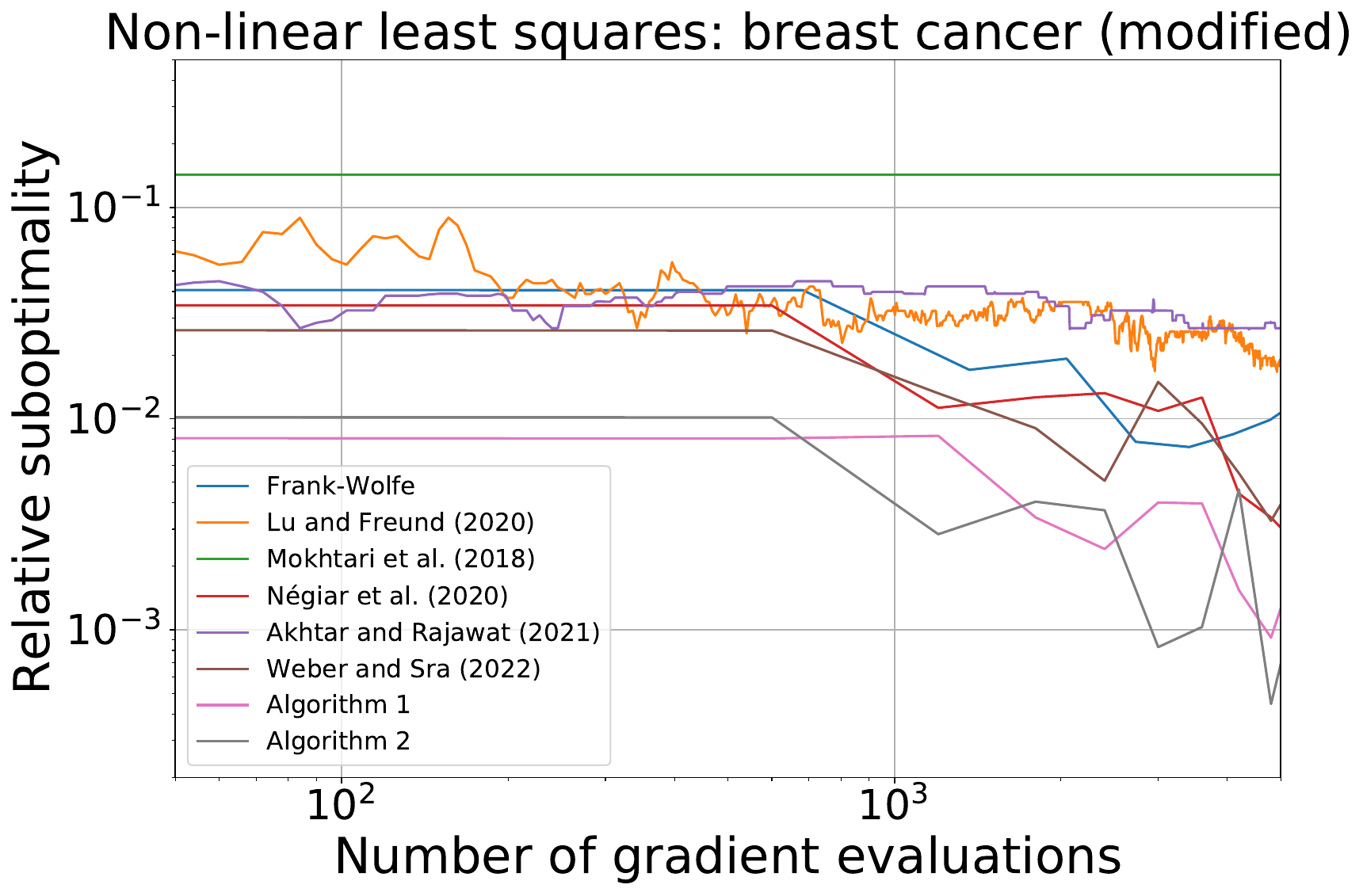}
\end{minipage}%
\\
\begin{minipage}{0.49\textwidth}
  \centering
(a) \texttt{a9a}
\end{minipage}%
\begin{minipage}{0.49\textwidth}
\centering
  (b) \texttt{breast cancer}
\end{minipage}%
\\
\begin{minipage}{0.24\textwidth}
\quad\quad\quad
\end{minipage}%
\begin{minipage}{0.24\textwidth}
  \centering
\includegraphics[width =  1\textwidth ]{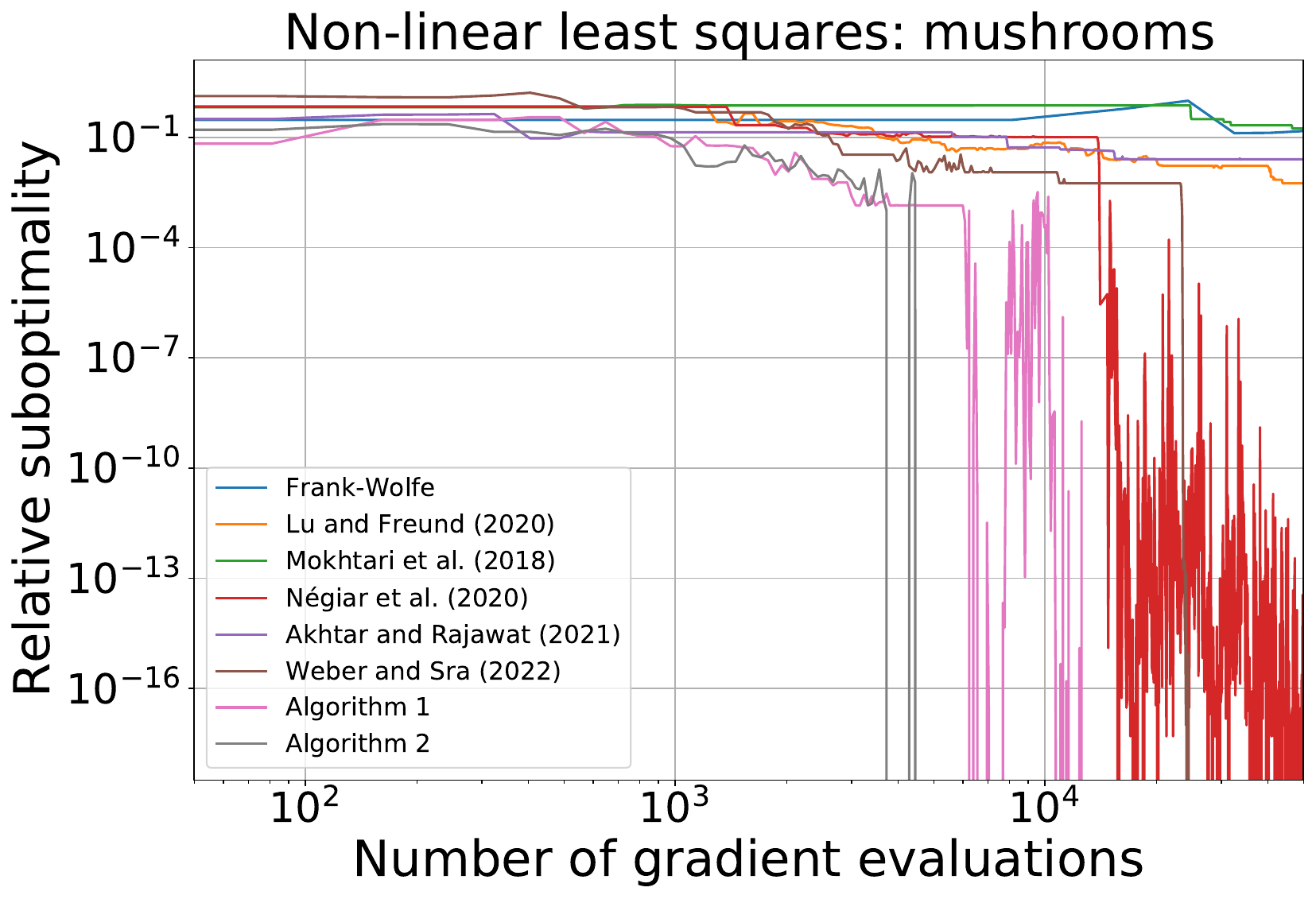}
\end{minipage}%
\begin{minipage}{0.24\textwidth}
  \centering
\includegraphics[width =  1\textwidth ]{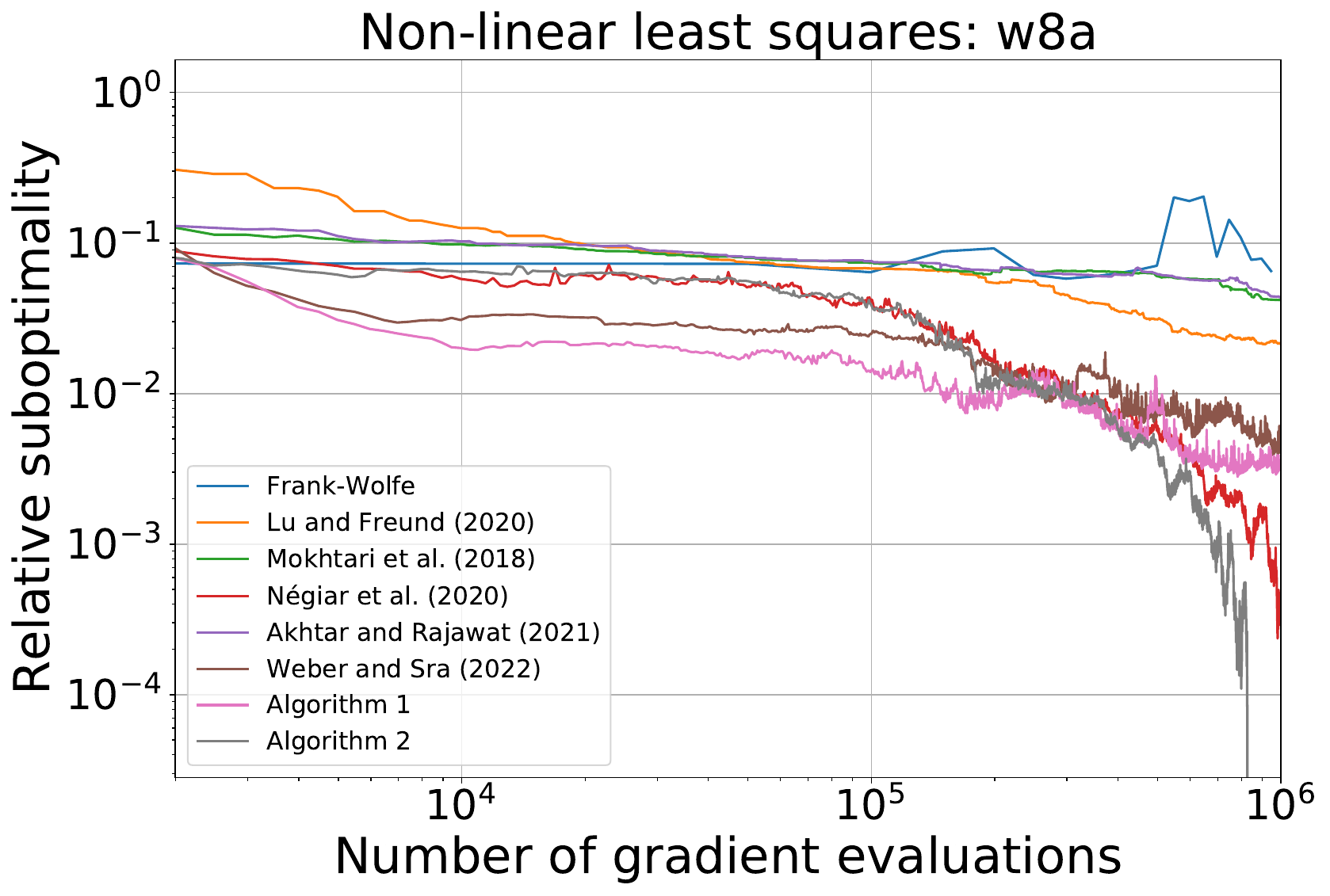}
\end{minipage}%
\\
\begin{minipage}{0.24\textwidth}
\quad\quad\quad
\end{minipage}%
\begin{minipage}{0.24\textwidth}
\centering
  (c) \texttt{mushrooms}
\end{minipage}%
\begin{minipage}{0.24\textwidth}
\centering
  (d) \texttt{w8a}
\end{minipage}%
\vspace{-0.3cm}
\caption{Comparison of state-of-the-art projection free methods with small batches for \eqref{eq;nls}. The comparison is made on the real datasets from LibSVM. The criterion is the number of full gradients computations. In the modified plots (the right plots in the first three lines), we left only every 100th point for \citep{negiar2020stochastic}, \citep{weber2022projection}, Algorithm \ref{alg:fw_sarah} and Algorithm \ref{alg:fw_zerosarah}.}
    \label{fig:min_nonconv}
\end{figure*}

We conduct our experiments on the constrained empirical risk for a linear model with weights $w$ and on training samples $\{x_i, y_i\}_{i=1}^n$. Then, the logistic regression problem is
\begin{equation}
    \label{eq;ls}
    \min_{w \in \C} f(w) = \frac{1}{n} \sum\limits_{i=1}^n \log(1 + \exp(-y_i w^T x_i )),
\end{equation}
where $y_i \in \{-1, 1\}$, and the non-linear least squares loss is
\begin{equation}
    \label{eq;nls}
    \min_{w \in \C} f(w) = \frac{1}{n} \sum\limits_{i=1}^n (y_i - 1 / (1 + \exp(w^T x_i)))^2,
\end{equation}

\setlength{\columnsep}{8pt}%
\begin{wrapfigure}[9]{r}{0.3\textwidth}
\renewcommand{\arraystretch}{1.2}
\vspace{-0.4cm}
    \centering
    % \small
%    \scriptsize
\hspace{-1cm}
\captionof{table}{Datasets from LibSVM in experiments.}
    \label{tab:data}   
    % \scriptsize
    \small
    % \resizebox{\linewidth}{!}{
  \begin{threeparttable}
   \centering
    \begin{tabular}{ccc}
    \hline
    \textbf{Dataset} & \textbf{$d$} & \textbf{$n$} \\
    \hline
    \texttt{a9a} & 123 & 22696 \\
    \texttt{breast cancer} & 10 & 683 \\
    \texttt{madelon} & 500 & 2000 \\
    \texttt{mushrooms} & 112 & 8124 \\
\texttt{w8a} & 300 & 49749
    \\\hline
    %%%%%%%%%%%%%%
    \end{tabular}   
    \end{threeparttable}
    % }
\end{wrapfigure} 

with $y_i \in \{0,1\}$. We consider two different loss functions to test our algorithms on both convex and nonconvex settings. We choose $\C$ as the $\ell_1$ norm ball with radius $R = 2 \cdot 10^3$ (see results with other $R$ in Section \ref{sec:add_exp}). The LMO for such a constraint set can be computed in a closed-form solution. We take 
LibSVM \citep{chang2011libsvm} datasets (see Table \ref{tab:data}). 

For comparison, we consider the methods from Table \ref{tab:comparison0}, which do not use large batches: \citep{mokhtari2020stochastic, negiar2020stochastic, lu2021generalized, 9483167, weber2022projection}. Our method is tuned according to the theory (see Sections \ref{sec:sarah} and \ref{sec:saga_sarah}), but we take the batchsize  $b = \left\lceil \tfrac{n}{100}\right\rceil$ (similarly as the baselines). For the baselines, we use the implementation of~\cite{negiar2020stochastic} and tune each method accordingly.\footnote{Note that the algorithms are not exactly implemented according to the theory of the corresponding works. In particular, instead of randomly selecting the batches, the authors sample them without replacement.} 

% \renewcommand{\arraystretch}{1.2}
% %\renewcommand{\tabcolsep}{6pt}
% %\setlength\extrarowheight{1pt}
% \begin{table}[!h]
% \vspace{-0.2cm}
%     \centering
%     % \small
% %    \scriptsize
% \captionof{table}{Datasets from LibSVM in experiments.}
% \vspace{-0.3cm}
%     \label{tab:data}   
%     % \scriptsize
%     \small
%     % \resizebox{\linewidth}{!}{
%   \begin{threeparttable}
%     \begin{tabular}{ccc}
%     \hline
%     \textbf{Dataset} & \textbf{$d$} & \textbf{$n$} \\
%     \hline
%     \texttt{madelon} & 500 & 2000 \\
%     \texttt{mushrooms} & 112 & 8124 \\
% \texttt{breast cancer} & 10 & 683 \\
% \texttt{a9a} & 123 & 22696 
%     \\\hline
%     %%%%%%%%%%%%%%
%     \end{tabular}   
%     \end{threeparttable}
%     % }
% \vspace{-0.2cm}
% \end{table}

In Figure \ref{fig:min} we plot the relative suboptimality which is defined a $\nicefrac{(f(x_t)-f_{\min})}{(f_{\max}-f_{\min})}$ where $f_{\max}$ is the
largest value observed along the optimization and $f_{\min}$ is obtained by running the best algorithm a bit longer than what is plotted. From our results, it is clear that our algorithms are superior or comparable to the baselines, despite the fact that some methods were specifically designed for linear models (e.g.~\citep{negiar2020stochastic, lu2021generalized})
which is not the case for Algorithm~\ref{alg:fw_sarah} and~\ref{alg:fw_zerosarah}.

%%%%%%%%%%%%%%%%%%%%%%%%%%%%%%%%%%%%%%%%
%%%%%%%%%%%%%%%%%%%%%%%%%%%%%%%%%%%%%%%%
%%%%%%%%%%%%%%%%%%%%%%%%%%%%%%%%%%%%%%%%
%%%%%%%%%%%%%%%%%%%%%%%%%%%%%%%%%%%%%%%%

\section{Conclusion and Future Works}

In this paper, we presented two new algorithms for stochastic finite-sum optimization. Our methods are based on the Frank-Wolfe and Sarah approaches. Both of our algorithms are free of target function parameters. In both convex and non-convex target cases, our algorithms have the best stochastic oracle complexity in the literature. Our methods do not need to resort to large batch computation. However, in the non-convex case, it is worth noticing that the methods with large batch sizes give a better oracle complexity estimate.
Moreover, Algorithm \ref{alg:fw_zerosarah} does not need to collect either large batches or full deterministic gradients at all. Our methods also perform well on different $\ell_1$ constrained logistic regression problems.

Ideas from~\citep{pmlr-v28-jaggi13}, and~\citep{lacoste2015global} can be noted as a starting point for the future research in order to get fast rates in the strongly convex case. The modifications presented in these papers make the Frank-Wolfe method more practical and faster. Combining these and our approaches can produce a strong synergy that results in new practical algorithms.

Recall that our estimates of LMO in the convex case and SFO in the non-convex setting achieve lower bounds, and thus can be considered as optimal and unimprovable. Meanwhile, the SFO results in the convex setup and the LMO in the non-convex case have an unclosed gap and potential for improvement in finding new algorithms or lower bounds, which give optimality of the current results.

% The question of obtaining lower bounds for stochastic projection-free methods with linear minimization oracle for constraint sets (in the deterministic case, a lower bound is known~\citep{pmlr-v28-jaggi13, lan2013complexity}) remains open. However, our results in Table \ref{tab:comparison0}, match the lower bound given for non-convex unconstrained minimization~\citep[Table 1]{pmlr-v139-li21a}). This indicates that our method is optimal (in terms of stochastic oracle call) for non-convex minimization.

Finally, as mentioned above, the analysis of the algorithms is based on constructing a recursive estimate not only on $f(x) - f(x^*)$, as for classical Frank-Wolfe, but also includes $\| g - \nabla f(x)\|^2$. But this kind of bound is applicable to a number of many modern methods. Here we can mention, for example, distributed optimization methods with compression \cite{NEURIPS2021_231141b3, pmlr-v139-gorbunov21a}, which are quite far from the setting of the current paper. A promising direction for future research is the generalization of the obtained results to a diverse array of optimization methods, with the aim of developing novel modifications of the Conditional Gradient method.

\section*{Acknowledgements}

The research of A. Beznosikov has been supported by the Analytical Center for the Government of the Russian Federation (Agreement No. 70-2021-00143 dd. 01.11.2021, IGK 000000D730321P5Q0002). 

\section*{Impact Statement}
This paper presents work whose goal is to advance the field of Machine Learning. There are many potential societal consequences of our work, none of which we feel must be specifically highlighted here.

\bibliography{lit}
\bibliographystyle{icml2024}

%%%%%%%%%%%%%%%%%%%%%%%%%%%%%%%%%%%%%%%%%%%%%%%%%%%%%%%%%%%%%%%%%%%%%%%%%%%%%%%
%%%%%%%%%%%%%%%%%%%%%%%%%%%%%%%%%%%%%%%%%%%%%%%%%%%%%%%%%%%%%%%%%%%%%%%%%%%%%%%
% APPENDIX
%%%%%%%%%%%%%%%%%%%%%%%%%%%%%%%%%%%%%%%%%%%%%%%%%%%%%%%%%%%%%%%%%%%%%%%%%%%%%%%
%%%%%%%%%%%%%%%%%%%%%%%%%%%%%%%%%%%%%%%%%%%%%%%%%%%%%%%%%%%%%%%%%%%%%%%%%%%%%%%
\newpage
\appendix
\onecolumn

\section{Technical Facts}

\begin{lemma} \label{lem:tech1}
For any $x_1, \ldots, x_n \in \R^d$ the following inequality holds:
\begin{equation*}
\left\| \sum\limits_{i=1}^n x_i \right\|^2 \leq n \sum\limits_{i=1}^n \| x_i\|^2.
\end{equation*}
\end{lemma}

\begin{lemma}[Lemma 1.2.3 from \citep{nesterov2003introductory}] \label{lem:tech2}
Suppose that $f$ is $L$-smooth. Then, for any $x,y \in \R^d$,
\begin{equation*}
f(x) \leq f(y) + \< \nabla f (y), x-y > + \frac{L}{2} \|x-y \|^2.
\end{equation*}
\end{lemma}

\begin{lemma}[Lemma 3 from \citep{stich2019unified}] \label{lem:tech3}
Let $\{r_k\}_{k \geq 0}$ is a non-negative sequence, which satisfies the relation
\begin{equation*}
r_{k+1} \leq (1 - \eta_k) r_k + c \eta_k^2.
\end{equation*}
Then there exists stepsizes $\eta_k \leq \tfrac{1}{d}$, such that:
\begin{equation*}
r_{K}  = \cO \left( dr_0 \exp\left(-\frac{K}{2d}\right) + \frac{c}{K}\right).
\end{equation*}
In particular, the step sizes $\{\eta\}_{k \geq 0}$ can be chosen as follows 
\begin{align*}
    &\text{if } ~~ K \leq d, && \hspace{-4cm}\eta_k = \frac{1}{d}, \\
    &\text{if } ~~ K > d ~~ \text{ and } ~~ k < k_0, && \hspace{-4cm} \eta_k = \frac{1}{d}, \\
    &\text{if } ~~ K > d ~~ \text{ and } ~~ k \geq k_0, && \hspace{-4cm} \eta_k = \frac{2}{(2d + k - k_0)},
\end{align*}
where $k_0 = \lceil \tfrac{K}{2}\rceil$.
\end{lemma}

%%%%%%%%%%%%%%%%%%%%%%%%%%%%%%%
%%%%%%%%%%%%%%%%%%%%%%%%%%%%%
%%%%%%%%%%%%%%%%%%%%%%%%%%%%%
%%%%%%%%%%%%%%%%%%%%%%%%%%%%%%%%%%%

\section{Missing Proofs}

\subsection{Proof of Theorem \ref{th:main1}} \label{sec:proof1}

\begin{theorem}[Theorem \ref{th:main1}] \label{th:1}
Let $\{x^k\}_{k\geq0}$ denote the iterates of Algorithm \ref{alg:fw_sarah} for solving problem \eqref{eq:main_problem}, which satisfies Assumptions \ref{as:lip}--\ref{as:set}. Let $x^*$ be the minimizer of $f$. Then for any $K$ one can choose $\{ \eta_k \}_{k \geq 0}$ as follows:
\begin{align*}
    &\text{if}~~  K \leq \frac{2}{p}, && \hspace{-4cm} \eta_k = \frac{p}{2}, \\
    &\text{if}~~  K > \frac{2}{p} ~~  \text{ and } ~~  k < k_0, && \hspace{-4cm}  \eta_k = \frac{p}{2}, \\
    &\text{if}~~  K > \frac{2}{p} ~~  \text{ and } ~~  k \geq k_0, && \hspace{-4cm} \eta_k = \frac{2}{(\nicefrac{4}{p} + k - \lceil \nicefrac{K}{2}\rceil)}.
\end{align*}
For this choice of $\eta_k$, we have the following convergence:
\begin{align*}
    \EE{f(x^{K}) - f(x^*)} =
    \mathcal{O}\left( \frac{1}{p} \left[f(x^{0}) - f(x^*) \right] \exp\left(-\frac{K p}{4} \right)   + \left[1 
    +  \frac{\tilde L}{L}\sqrt{\frac{1-p}{p b}}\right] \frac{LD^2}{K}\right).
\end{align*}
\end{theorem}

\textbf{Proof:} 
Let us start with Assumption \ref{as:lip} and Lemma \ref{lem:tech2}:
\begin{align*}
    f(x^{k+1}) \leq f(x^k) + \<\nabla f(x^k), x^{k+1} - x^k > + \frac{L}{2} \|x^{k+1} - x^k \|^2.
\end{align*}
Subtracting $f(x^*)$ from both sides, we get
\begin{align*}
    f(x^{k+1}) - f(x^*) \leq f(x^k) - f(x^*) + \<\nabla f(x^k), x^{k+1} - x^k > + \frac{L}{2} \|x^{k+1} - x^k \|^2.
\end{align*}
With the update of $x^{k+1}$ from line \ref{alg1:line3} of Algorithm \ref{alg:fw_sarah}, one can obtain
\begin{align*}
    f(x^{k+1}) - f(x^*) 
    \leq&
    f(x^k) - f(x^*) + \eta_k \<\nabla f(x^k), s^{k} - x^k > + \frac{L \eta_k^2}{2} \|s^{k} - x^k \|^2
    \\=&
    f(x^k) - f(x^*) + \eta_k \<g^k, s^{k} - x^k > + \eta_k \<\nabla f(x^k) - g^k, s^{k} - x^k > 
    \\&
    + \frac{L \eta_k^2}{2} \|s^{k} - x^k \|^2.
\end{align*}
The optimal choice of $s^k$ from line \ref{alg1:line2} gives that $\< g^k, s^k - x^k> \leq \< g^k, x^* - x^k>$. Then,
\begin{align*}
    f(x^{k+1}) - f(x^*) 
    \leq& 
    f(x^k) - f(x^*) + \eta_k \<g^k, x^* - x^k > + \eta_k \<\nabla f(x^k) - g^k, s^{k} - x^k > 
    \\&
    + \frac{L \eta_k^2}{2} \|s^{k} - x^k \|^2
    \\=& 
    f(x^k) - f(x^*) + \eta_k \<\nabla f(x^k), x^* - x^k > + \eta_k \<g^k - \nabla f(x^k), x^* - x^k > 
    \\&
    + \eta_k \<\nabla f(x^k) - g^k, s^{k} - x^k > + \frac{L \eta_k^2}{2} \|s^{k} - x^k \|^2
    \\=&
    f(x^k) - f(x^*) + \eta_k \<\nabla f(x^k), x^* - x^k > + \eta_k \<\nabla f(x^k) - g^k, s^{k} - x^* > 
    \\&
    + \frac{L \eta_k^2}{2} \|s^{k} - x^k \|^2.
\end{align*}
Applying the Cauchy-Schwartz inequality, we deduce $\<\tfrac{\sqrt{\alpha}}{\sqrt{L}}(\nabla f(x^k) - g^k), \tfrac{\sqrt{L}}{\sqrt{\alpha}}\eta_k (s^{k} - x^*) > \leq \frac{\alpha}{L}\| \nabla f(x^k) - g^k \|^2 + \frac{L \eta_k^2}{\alpha}\|s^{k} - x^* \|^2$ with some positive constant $\alpha$ (which we will define below). Thus, 
\begin{align*}
    f(x^{k+1}) - f(x^*) 
    \leq& 
    f(x^k) - f(x^*) + \eta_k \<\nabla f(x^k), x^* - x^k > + \frac{\alpha}{L}\| \nabla f(x^k) - g^k \|^2 
    \\&
    + \frac{L \eta_k^2}{\alpha}\|s^{k} - x^* \|^2
    + \frac{L \eta_k^2}{2} \|s^{k} - x^k \|^2.
\end{align*}
Using the convexity of the function $f$ (Assumption \ref{as:conv}): $\<\nabla f(x^k), x^* - x^k > \leq - (f(x^k) - f(x^*))$, we have
\begin{align*}
    f(x^{k+1}) - f(x^*) 
    \leq& 
    f(x^k) - f(x^*) - \eta_k (f(x^k) - f(x^*)) + \frac{\alpha}{L}\| \nabla f(x^k) - g^k \|^2 
    \nonumber\\&
    + \frac{L \eta_k^2}{\alpha}\|s^{k} - x^* \|^2
    + \frac{L \eta_k^2}{2} \|s^{k} - x^k \|^2
    \nonumber\\=&
    (1 - \eta_k)(f(x^k) - f(x^*)) + \frac{\alpha}{L}\| \nabla f(x^k) - g^k \|^2 
    + \frac{L \eta_k^2}{\alpha}\|s^{k} - x^* \|^2
    \\&
    + \frac{L \eta_k^2}{2} \|s^{k} - x^k \|^2.
\end{align*}
Taking the full mathematical expectation, one can obtain
\begin{align}
    \label{eq:temp1}
    \EE{f(x^{k+1}) - f(x^*)}
    \leq& 
    (1 - \eta_k) \EE{f(x^k) - f(x^*)} + \frac{\alpha}{L} \EE{\| \nabla f(x^k) - g^k \|^2}
    \nonumber\\&
    + \frac{L \eta_k^2}{\alpha} \EE{\|s^{k} - x^* \|^2}
    + \frac{L \eta_k^2}{2} \EE{\|s^{k} - x^k \|^2}.
\end{align}
For line \ref{alg1:line4} we use Lemma 3 from \citep{pmlr-v139-li21a}:
\begin{align*}
    \EE{\| \nabla f(x^{k+1}) - g^{k+1} \|^2} \leq (1 - p) \EE{\| \nabla f(x^{k}) - g^{k} \|^2} + \frac{1-p}{b} \EE{\| \nabla f_{i_k}(x^{k+1}) - \nabla f_{i_k}(x^{k})\|^2}.
\end{align*}
With Assumption \ref{as:lip}, we get 
\begin{align*}
    \EE{\| \nabla f(x^{k+1}) - g^{k+1} \|^2} 
    \leq& 
    (1 - p) \EE{\| \nabla f(x^{k}) - g^{k} \|^2} + \frac{1-p}{b}\EE{L^2_{i_k} \| x^{k+1} - x^{k}\|^2}
    \\=&
    (1 - p) \EE{\| \nabla f(x^{k}) - g^{k} \|^2} + \frac{1-p}{b} \EE{ \E_{i_k}[L^2_{i_k}] \cdot \| x^{k+1} - x^{k}\|^2}.
\end{align*}
In the last step, we use the independence of $i_k$ and $(x^{k+1} - x^k)$. Taking expectation on $i_k$, we obtain
\begin{align*}
    \EE{\| \nabla f(x^{k+1}) - g^{k+1} \|^2} \leq (1 - p) \EE{\| \nabla f(x^{k}) - g^{k} \|^2} + \frac{1-p}{b} \left( \frac{1}{n} \sum\limits_{i=1}^n L_i^2 \right) \EE{\| x^{k+1} - x^{k}\|^2}.
\end{align*}
The notation of $\tilde L$ provides
\begin{align}
    \label{eq:temp2}
    \EE{\| \nabla f(x^{k+1}) - g^{k+1} \|^2} \leq (1 - p) \EE{\| \nabla f(x^{k}) - g^{k} \|^2} + \frac{(1-p) \tilde L^2}{b}  \EE{\| x^{k+1} - x^{k}\|^2}.
\end{align}
Multiplying \eqref{eq:temp2} by the positive constant $M$ (which we will define below) and summing with \eqref{eq:temp1}, we have
\begin{align*}
    \E \big[f(x^{k+1}) - f(x^*) &+ M \cdot \| \nabla f(x^{k+1}) - g^{k+1} \|^2\big]
    \\\leq& 
    (1 - \eta_k) \EE{f(x^k) - f(x^*)} + \left( 1 - p + \frac{\alpha}{ML}\right)M\cdot \EE{\| \nabla f(x^k) - g^k \|^2} 
    \nonumber\\&
    + \frac{L \eta_k^2}{\alpha} \EE{\|s^{k} - x^* \|^2}
    + \frac{L \eta_k^2}{2} \EE{\|s^{k} - x^k \|^2} + \frac{M(1-p) \tilde L^2 \eta_k^2}{b} \EE{\| s^{k} - x^{k}\|^2}.
\end{align*}
The choice of $M = 2\alpha / (p L)$ gives
\begin{align*}
    \E \big[f(x^{k+1}) - f(x^*) &+ M \cdot \| \nabla f(x^{k+1}) - g^{k+1} \|^2\big]
    \\\leq& 
    (1 - \eta_k) \EE{f(x^k) - f(x^*)} + \left( 1 - \frac{p}{2} \right)M\cdot \EE{\| \nabla f(x^k) - g^k \|^2} 
    \nonumber\\&
    + \frac{L \eta_k^2}{\alpha}\EE{\|s^{k} - x^* \|^2}
    + \frac{L \eta_k^2}{2} \EE{\|s^{k} - x^k \|^2} + \frac{2\alpha (1-p) \tilde L^2 \eta_k^2}{p b L} \EE{\| s^{k} - x^{k}\|^2}
    \\\leq&
    \max\left\{ 1 - \eta_k, 1 - \frac{p}{2} \right\} \EE{ f(x^k) - f(x^*) + M\cdot \| \nabla f(x^k) - g^k \|^2 }
    \nonumber\\&
    + \frac{L \eta_k^2}{\alpha} \EE{\|s^{k} - x^* \|^2}
    + \frac{L \eta_k^2}{2} \EE{\|s^{k} - x^k \|^2} + \frac{2\alpha (1-p) \tilde L^2 \eta_k^2}{p b L} \EE{\| s^{k} - x^{k}\|^2}.
\end{align*}
With Assumption \ref{as:set} on the diameter $D$ of $\X$, we get
\begin{align*}
    \E \big[f(x^{k+1}) - f(x^*) &+ M \cdot \| \nabla f(x^{k+1}) - g^{k+1} \|^2\big]
    \\\leq& 
    \max\left\{ 1 - \eta_k, 1 - \frac{p}{2} \right\} \EE{ f(x^k) - f(x^*) + M\cdot \| \nabla f(x^k) - g^k \|^2 }
    \nonumber\\&
    + \frac{LD^2 \eta_k^2}{\alpha}
    + \frac{LD^2 \eta_k^2}{2} + \frac{2\alpha (1-p) \tilde L^2 D^2\eta_k^2}{p b L}
    \\=& 
    \max\left\{ 1 - \eta_k, 1 - \frac{p}{2} \right\} \EE{ f(x^k) - f(x^*) + M\cdot \| \nabla f(x^k) - g^k \|^2 }
    \nonumber\\&
    + LD^2 \eta_k^2 \left(\frac{1}{2} + \frac{1}{\alpha}
    +  \frac{2\alpha (1-p) \tilde L^2}{p b L^2}\right).
\end{align*}
If we choose $\eta_k \leq \frac{p}{2}$, $\alpha = \sqrt{\tfrac{p b L^2}{(1-p) \tilde L^2}}$, then we have
\begin{align*}
    \E \big[f(x^{k+1}) - f(x^*) &+ M \cdot \| \nabla f(x^{k+1}) - g^{k+1} \|^2\big]
    \\\leq& 
    \left( 1 - \eta_k \right) \EE{ f(x^k) - f(x^*) + M\cdot \| \nabla f(x^k) - g^k \|^2 }
    \nonumber\\&
    + LD^2 \eta_k^2 \left(\frac{1}{2} 
    +  \frac{3\tilde L}{L}\sqrt{\frac{1-p}{pb}}\right).
\end{align*}
It remains to use Lemma \ref{lem:tech3} with $c = LD^2\left(\tfrac{1}{2} 
    +  \tfrac{3\tilde L}{L} \sqrt{\tfrac{1-p}{pb}}\right)$, $d = \tfrac{2}{p}$ and obtain
\begin{align*}
    \E \big[f(x^{K}) & - f(x^*) + M \cdot \| \nabla f(x^{K}) - g^{K}\|^2 \big]
    \\=& 
    \mathcal{O}\left( \frac{1}{p} \left[f(x^{0}) - f(x^*) + M \cdot \| \nabla f(x^{0}) - g^{0} \|^2\right] \exp\left(-\frac{K p}{4} \right)   + \left[1 
    +  \frac{\tilde L}{L}\sqrt{\frac{1-p}{pb}}\right] \frac{LD^2}{K}\right) .
\end{align*}
Finally, we substitute $g^0 = \nabla f(x^0)$:
\begin{align*}
    \EE{f(x^{K}) - f(x^*)} =
    \mathcal{O}\left( \frac{1}{p} \left[f(x^{0}) - f(x^*) \right] \exp\left(-\frac{K p}{4} \right)   + \left[1 
    +  \frac{\tilde L}{L}\sqrt{\frac{1-p}{p b}}\right] \frac{LD^2}{K}\right).
\end{align*}
This completes the proof.
\EndProof

%%%%%%%%%%%%%%%%%%%%%%%%%%%%%%%%%%%%%%%%%%%%%%%%%%%%%%%%%%%%%%%%%%%%%%%%
%%%%%%%%%%%%%%%%%%%%%%%%%%%%%%%%%%%%%%%%%%%%%%%%%%%%%%%%%%%%%%%%%%%%%%%%
%%%%%%%%%%%%%%%%%%%%%%%%%%%%%%%%%%%%%%%%%%%%%%%%%%%%%%%%%%%%%%%%%%%%%%%%
%%%%%%%%%%%%%%%%%%%%%%%%%%%%%%%%%%%%%%%%%%%%%%%%%%%%%%%%%%%%%%%%%%%%%%%%

\subsection{Proof of Theorem \ref{th:main2}} \label{sec:proof2}

\begin{theorem}[Theorem \ref{th:main2}] \label{th:2}
Let $\{x^k\}_{k\geq0}$ denote the iterates of Algorithm \ref{alg:fw_sarah} for solving problem \eqref{eq:main_problem}, which satisfies Assumptions \ref{as:lip},\ref{as:set}. Let $x^*$ be the global (may be not unique) minimizer of $f$. Then, if we choose $\eta_k = \frac{1}{\sqrt{K}}$, we have the following convergence:
\begin{align*}
    \E \left[\min_{0\leq k \leq K-1}\gap(x^k)\right]
    = 
    \mathcal{O} \left(\frac{f(x^0) - f(x^*)}{\sqrt{K}} + \frac{LD^2}{\sqrt{K}} \left[1 
    +  \frac{\tilde L}{L} \sqrt{\frac{(1-p)}{pb}}\right]\right).
\end{align*}
\end{theorem}

\textbf{Proof:} 
Let us start with Assumption \ref{as:lip} and Lemma \ref{lem:tech2}:
\begin{align*}
    f(x^{k+1}) \leq f(x^k) + \<\nabla f(x^k), x^{k+1} - x^k > + \frac{L}{2} \|x^{k+1} - x^k \|^2.
\end{align*}
Subtracting $f(x^*)$ from both sides, we get
\begin{align*}
    f(x^{k+1}) - f(x^*) \leq f(x^k) - f(x^*) + \<\nabla f(x^k), x^{k+1} - x^k > + \frac{L}{2} \|x^{k+1} - x^k \|^2.
\end{align*}
With the update of $x^{k+1}$ from line \ref{alg1:line3} of Algorithm \ref{alg:fw_sarah}, one can obtain
\begin{align*}
    f(x^{k+1}) - f(x^*) 
    \leq& 
    f(x^k) - f(x^*) + \eta_k \<\nabla f(x^k), s^{k} - x^k > + \frac{L \eta_k^2}{2} \|s^{k} - x^k \|^2
    \\=&
    f(x^k) - f(x^*) + \eta_k \<g^k, s^{k} - x^k > + \eta_k \<\nabla f(x^k) - g^k, s^{k} - x^k > + \frac{L \eta_k^2}{2} \|s^{k} - x^k \|^2.
\end{align*}
The optimal choice of $s^k$ from line \ref{alg1:line2} gives that $\< g^k, s^k - x^k> \leq \< g^k, x - x^k>$ for any $x \in \X$. Then,
\begin{align*}
    f(x^{k+1}) - f(x^*) 
    \leq&
    f(x^k) - f(x^*) + \eta_k \<g^k, x - x^k > + \eta_k \<\nabla f(x^k) - g^k, s^{k} - x^k > + \frac{L \eta_k^2}{2} \|s^{k} - x^k \|^2
    \\=& 
    f(x^k) - f(x^*) + \eta_k \<\nabla f(x^k), x - x^k > + \eta_k \<g^k - \nabla f(x^k), x - x^k > 
    \\&
    + \eta_k \<\nabla f(x^k) - g^k, s^{k} - x^k > + \frac{L \eta_k^2}{2} \|s^{k} - x^k \|^2
    \\=&
    f(x^k) - f(x^*) + \eta_k \<\nabla f(x^k), x - x^k > + \eta_k \<\nabla f(x^k) - g^k, s^{k} - x > 
    \\&
    + \frac{L \eta_k^2}{2} \|s^{k} - x^k \|^2.
\end{align*}
Applying the Cauchy-Schwartz inequality, we deduce $\<\tfrac{\sqrt{\alpha}}{\sqrt{L}}(\nabla f(x^k) - g^k), \tfrac{\sqrt{L}}{\sqrt{\alpha}}\eta_k (s^{k} - x^*) > \leq \frac{\alpha}{L}\| \nabla f(x^k) - g^k \|^2 + \frac{L \eta_k^2}{\alpha}\|s^{k} - x^* \|^2$ with some positive constant $\alpha$ (which we will define below). Thus, 
\begin{align*}
    f(x^{k+1}) - f(x^*) 
    \leq& 
    f(x^k) - f(x^*) + \eta_k \<\nabla f(x^k), x - x^k > + \frac{\alpha}{L}\| \nabla f(x^k) - g^k \|^2 
    \nonumber\\&
    + \frac{L \eta_k^2}{\alpha}\|s^{k} - x \|^2
    + \frac{L \eta_k^2}{2} \|s^{k} - x^k \|^2.
\end{align*}
After small rearrangements one can obtain
\begin{align*}
    \eta_k \<\nabla f(x^k), x^k - x> 
    \leq& 
    f(x^k) - f(x^*) - \left(f(x^{k+1}) - f(x^*) \right)  + \frac{\alpha}{L}\| \nabla f(x^k) - g^k \|^2 
    \nonumber\\&
    + \frac{L \eta_k^2}{\alpha}\|s^{k} - x \|^2
    + \frac{L \eta_k^2}{2} \|s^{k} - x^k \|^2.
\end{align*}
Maximizing over all $x \in \X$ and taking the full mathematical expectation, we get
\begin{align}
    \label{eq:temp3}
    \eta_k \EE{ \max_{x \in \X}\<\nabla f(x^k), x^k - x >}
    \leq& 
    \EE{f(x^k) - f(x^*)} - \EE{f(x^{k+1}) - f(x^*)}  + \frac{\alpha}{L} \EE{\| \nabla f(x^k) - g^k \|^2} 
    \nonumber\\&
    + \frac{L \eta_k^2}{\alpha} \EE{\max_{x \in \X}\|s^{k} - x \|^2}
    + \frac{L \eta_k^2}{2} \EE{\|s^{k} - x^k \|^2}.
\end{align}
Multiplying \eqref{eq:temp2} by the positive constant $M$ (which we will define below) and summing with \eqref{eq:temp3}, we have
\begin{align*}
    \eta_k \EE{ \max_{x \in \X}\<\nabla f(x^k), x^k - x >} 
    \leq& 
    \EE{f(x^k) - f(x^*)} + \left( 1 - p + \frac{\alpha}{ML}\right) M \EE{\| \nabla f(x^k) - g^k \|^2} 
    \nonumber\\&
    - \E \big[f(x^{k+1}) - f(x^*) + M \| \nabla f(x^{k+1}) - g^{k+1}\|^2 \big]
    \nonumber\\&
    + \frac{L \eta_k^2}{\alpha} \EE{\max_{x \in \X} \|s^{k} - x \|^2}
    + \frac{L \eta_k^2}{2} \EE{\|s^{k} - x^k \|^2}
    \nonumber\\&
    + \frac{M (1-p) \tilde L^2 \eta_k^2}{b} \EE{\| s^{k} - x^{k}\|^2}.
\end{align*}
The choice of $M = \alpha / (p L)$ gives
\begin{align*}
    \eta_k \EE{ \max_{x \in \X}\<\nabla f(x^k), x^k - x >} 
    \leq& 
    \EE{f(x^k) - f(x^*) + M\| \nabla f(x^k) - g^k \|^2}
    \nonumber\\&
    - \E \big[f(x^{k+1}) - f(x^*) + M \| \nabla f(x^{k+1}) - g^{k+1}\|^2 \big]
    \nonumber\\&
    + \frac{L \eta_k^2}{\alpha} \EE{\max_{x \in \X} \|s^{k} - x \|^2}
    + \frac{L \eta_k^2}{2} \EE{\|s^{k} - x^k \|^2}
    \nonumber\\&
    + \frac{\alpha(1-p) \tilde L^2 \eta_k^2}{pbL} \EE{\| s^{k} - x^{k}\|^2}.
\end{align*}
With Assumption \ref{as:set} on the diameter $D$ of $\X$, we get
\begin{align*}
    \eta_k \EE{ \max_{x \in \X}\<\nabla f(x^k), x - x^k >} 
    \leq& 
    \EE{f(x^k) - f(x^*) + M\| \nabla f(x^k) - g^k \|^2}
    \nonumber\\&
    - \E \big[f(x^{k+1}) - f(x^*) + M \| \nabla f(x^{k+1}) - g^{k+1}\|^2 \big]
    \nonumber\\&
    + L D^2 \eta_k^2 \left( \frac{1}{2} + \frac{1}{\alpha} + \frac{\alpha (1-p) \tilde L^2 \eta_k^2}{p b L^2} \right).
\end{align*}
With the choice $\alpha = \tfrac{L}{\tilde L} \sqrt{\tfrac{pb}{1-p}}$, we have
\begin{align*}
    \eta_k \EE{\max_{x \in \X} \<\nabla f(x^k), x^k - x >}
    \leq& 
    \EE{f(x^k) - f(x^*) + M\| \nabla f(x^k) - g^k \|^2}
    \nonumber\\&
    -\EE{f(x^{k+1}) - f(x^*) + M\| \nabla f(x^{k+1}) - g^{k+1} \|^2}
    \nonumber\\&
    + LD^2 \eta_k^2 \left(\frac{1}{2} 
    +  \frac{2\tilde L}{L} \sqrt{\frac{(1-p)}{pb}}\right).
\end{align*}
Summing over all $k$ from $0$ to $K-1$, we have
\begin{align*}
    \sum\limits_{k=0}^{K-1}\eta_k \EE{\max_{x \in \X} \<\nabla f(x^k), x^k - x >}
    \leq& 
    f(x^0) - f(x^*) + M\| \nabla f(x^0) - g^0 \|^2
    \nonumber\\&
    -\EE{f(x^{K}) - f(x^*) + M \| \nabla f(x^{K}) - g^{K} \|^2}
    \nonumber\\&
    + LD^2 \left(\frac{1}{2} 
    +  \frac{2\tilde L}{L} \sqrt{\frac{(1-p)q}{pn}}\right) \sum\limits_{k=0}^{K-1} \eta_k^2
    \\\leq&
    f(x^0) - f(x^*) + \| \nabla f(x^0) - g^0 \|^2
    \nonumber\\&
    + LD^2 \left(\frac{1}{2} 
    +  \frac{2\tilde L}{L} \sqrt{\frac{(1-p)}{pb}}\right) \sum\limits_{k=0}^{K-1} \eta_k^2.
\end{align*}
If we take $\eta_k = \frac{1}{\sqrt{K}}$ and divide both sides by $\sqrt{K}$, then 
\begin{align*}
    \EE{\frac{1}{K}\sum\limits_{k=0}^{K-1} \max_{x \in \X} \<\nabla f(x^k), x^k - x >}
    \leq& 
    \frac{1}{\sqrt{K}} \cdot \left[f(x^0) - f(x^*) + M\| \nabla f(x^0) - g^0 \|^2\right]
    \nonumber\\&
    + \frac{LD^2}{\sqrt{K}} \left(\frac{1}{2} 
    +  \frac{2\tilde L}{L} \sqrt{\frac{(1-p)}{pb}}\right).
\end{align*}
Finally, we substitute $g^0 = \nabla f(x^0)$:
\begin{align*}
    \EE{\frac{1}{K}\sum\limits_{k=0}^{K-1} \max_{x \in \X} \<\nabla f(x^k), x^k - x >}
    \leq& 
    \frac{f(x^0) - f(x^*)}{\sqrt{K}} + \frac{LD^2}{\sqrt{K}} \left(\frac{1}{2} 
    +  \frac{2\tilde L}{L} \sqrt{\frac{(1-p)}{pb}}\right).
\end{align*}
The definition of \eqref{eq:gap} finishes the proof.
\EndProof

%%%%%%%%%%%%%%%%%%%%%%%%%%%%%%%%%%%%%%%%%%%%%%%%%%%%%%%%%%%%%%%%%%%%%%%%
%%%%%%%%%%%%%%%%%%%%%%%%%%%%%%%%%%%%%%%%%%%%%%%%%%%%%%%%%%%%%%%%%%%%%%%%
%%%%%%%%%%%%%%%%%%%%%%%%%%%%%%%%%%%%%%%%%%%%%%%%%%%%%%%%%%%%%%%%%%%%%%%%
%%%%%%%%%%%%%%%%%%%%%%%%%%%%%%%%%%%%%%%%%%%%%%%%%%%%%%%%%%%%%%%%%%%%%%%%

\subsection{Proof of Theorem \ref{th:main3}} \label{sec:proof3}

\begin{theorem}[Theorem \ref{th:main3}] \label{th:3}
Let $\{x^k\}_{k\geq0}$ denote the iterates of Algorithm \ref{alg:fw_zerosarah} for solving problem \eqref{eq:main_problem}, which satisfies Assumptions \ref{as:lip}--\ref{as:set}. Let $x^*$ be the minimizer of $f$. Then for any $K$ one can choose $\{ \eta_k \}_{k \geq 0}$ as follows:
\begin{align*}
    &\text{if}~~  K \leq \frac{4n}{b}, && \hspace{-4cm} \eta_k = \frac{b}{4n}, \\
    &\text{if}~~  K > \frac{4n}{b} ~~  \text{ and } ~~  k < k_0, && \hspace{-4cm} \eta_k = \frac{b}{4n}, \\
    &\text{if}~~  K > \frac{4n}{b} ~~  \text{ and } ~~  k \geq k_0, && \hspace{-4cm} \eta_k = \frac{2}{(\nicefrac{8n}{b} + k - \lceil \nicefrac{K}{2}\rceil)},
\end{align*}
and $\lambda = \tfrac{b}{2n}$. For this choice of $\eta_k$ and $\lambda$, we have the following convergence:
\begin{align*}
    \EE{f(x^{K}) - f(x^*)} 
    =
    \cO\Bigg( \frac{n}{b} \left[f(x^{0}) - f(x^*) \right] \exp\left(-\frac{bK}{8n} \right)   + \left[1 
    +  \frac{\tilde L \sqrt{n}}{L b}\right] \frac{LD^2}{K}\Bigg).
\end{align*}
\end{theorem}

\textbf{Proof:} The first steps of the proof are the same with Theorem \ref{th:main1} (Theorem \ref{th:1} ), therefore we can start from \eqref{eq:temp1}. For line \ref{alg2:line4} we use Lemma 2 from \citep{li2021zerosarah}:
\begin{align*}
    \EE{\| \nabla f(x^{k+1}) - g^{k+1} \|^2} 
    \leq& 
    (1 - \lambda) \EE{\| \nabla f(x^{k}) - g^{k} \|^2} + \frac{2}{b} \EE{\| \nabla f_{i_k}(x^{k+1}) - \nabla f_{i_k}(x^{k})\|^2} 
    \\&
    + \frac{2 \lambda^2}{b} \cdot \frac{1}{n} \sum\limits_{j=1}^n \EE{\|\nabla f_j (x^k) - y_j^k \|^2}.
\end{align*}
With Assumption \ref{as:lip}, we get 
\begin{align*}
    \EE{\| \nabla f(x^{k+1}) - g^{k+1} \|^2} 
    \leq& 
    (1 - \lambda) \EE{\| \nabla f(x^{k}) - g^{k} \|^2} + \frac{2}{b}\EE{L^2_{i_k} \| x^{k+1} - x^{k}\|^2}
    \\&
    + \frac{2 \lambda^2}{b} \cdot \frac{1}{n} \sum\limits_{j=1}^n \EE{\|\nabla f_j (x^k) - y_j^k \|^2}
    \\=&
    (1 - \lambda) \EE{\| \nabla f(x^{k}) - g^{k} \|^2} + \frac{2}{b} \EE{ \E_{i_k}[L^2_{i_k}] \cdot \| x^{k+1} - x^{k}\|^2}
    \\&
    + \frac{2 \lambda^2}{b} \cdot \frac{1}{n} \sum\limits_{j=1}^n \EE{\|\nabla f_j (x^k) - y_j^k \|^2}.
\end{align*}
In the last step, we use the independence of $i_k$ and $(x^{k+1} - x^k)$. Taking expectation on $i_k$, we obtain
\begin{align*}
    \EE{\| \nabla f(x^{k+1}) - g^{k+1} \|^2} 
    \leq&
    (1 - \lambda) \EE{\| \nabla f(x^{k}) - g^{k} \|^2} + \frac{2}{b} \left( \frac{1}{n} \sum\limits_{i=1}^n L_i^2 \right) \EE{\| x^{k+1} - x^{k}\|^2}
    \\&
    + \frac{2 \lambda^2}{b} \cdot \frac{1}{n} \sum\limits_{j=1}^n \EE{\|\nabla f_j (x^k) - y_j^k \|^2}.
\end{align*}
The notation of $\tilde L$ provides
\begin{align}
    \label{eq:temp4}
    \EE{\| \nabla f(x^{k+1}) - g^{k+1} \|^2} 
    \leq&
    (1 - \lambda) \EE{\| \nabla f(x^{k}) - g^{k} \|^2} + \frac{2 \tilde L^2}{b} \EE{\| x^{k+1} - x^{k}\|^2} 
    \nonumber\\&
    + \frac{2 \lambda^2}{b} \cdot \frac{1}{n} \sum\limits_{j=1}^n \EE{\|\nabla f_j (x^k) - y_j^k \|^2}.
\end{align}
Additionally, we need Lemma 3 from \citep{li2021zerosarah} with $\beta_k = \tfrac{b}{2n}$:
\begin{align*}
    \frac{1}{n} \sum\limits_{j=1}^n \EE{\|\nabla f_j (x^{k+1}) - y_j^{k+1} \|^2} 
    \leq&
    \left(1 - \frac{b}{2n} \right)\cdot \frac{1}{n} \sum\limits_{j=1}^n \EE{\|\nabla f_j (x^{k}) - y_j^{k} \|^2} 
    \\&
    + \frac{2n}{b} \cdot \frac{1}{n} \sum\limits_{i=1}^n \EE{\| \nabla f_{j}(x^{k+1}) - \nabla f_{j}(x^{k})\|^2}.
\end{align*}
With Assumption \ref{as:lip} and the notation of $\tilde L$, we get 
\begin{align}
    \label{eq:temp5}
    \frac{1}{n} \sum\limits_{j=1}^n \EE{\|\nabla f_j (x^{k+1}) - y_j^{k+1} \|^2} 
    \leq& 
    \left(1 - \frac{b}{2n} \right)\cdot \frac{1}{n} \sum\limits_{j=1}^n \EE{\|\nabla f_j (x^{k}) - y_j^{k} \|^2}
    \nonumber\\&
    + \frac{2n \tilde L^2}{b} \EE{\| x^{k+1} - x^{k}\|^2}.
\end{align}
Multiplying \eqref{eq:temp4} by the positive constant $M_1$, \eqref{eq:temp5} by the positive constant $M_2$ ($M_1$, $M_2$ will be defined below) and summing with \eqref{eq:temp1}, we have
\begin{align*}
    \E \Bigg[f(x^{k+1}) - f(x^*) &+ M_1 \cdot \| \nabla f(x^{k+1}) - g^{k+1} \|^2 + M_2 \cdot \frac{1}{n} \sum\limits_{j=1}^n \|\nabla f_j (x^{k+1}) - y_j^{k+1} \|^2 \Bigg]
    \\\leq& 
    (1 - \eta_k)\EE{f(x^k) - f(x^*)} 
    \\&
    +\left(1 - \lambda + \frac{\alpha}{ M_1 L}\right) M_1 \cdot \EE{\| \nabla f(x^{k}) - g^{k} \|^2}
    \\&
    + \left(1 - \frac{b}{2n} + \frac{2 M_1 \lambda^2}{M_2 b}\right) M_2\cdot \frac{1}{n} \sum\limits_{j=1}^n \EE{\|\nabla f_j (x^{k}) - y_j^{k} \|^2}
    \\&
    + \frac{2( M_1 + n M_2) \tilde L^2}{b} \EE{\| x^{k+1} - x^{k}\|^2}
    + \frac{L \eta_k^2}{\alpha} \EE{\|s^{k} - x^* \|^2}
    \\&
    + \frac{L \eta_k^2}{2} \EE{\|s^{k} - x^k \|^2}.
\end{align*}
With $M_1 = \tfrac{2\alpha}{\lambda L}$ and $M_2 = \tfrac{8M_1 \lambda^2 n}{b^2}$, we obtain
\begin{align*}
    \E \Bigg[f(x^{k+1}) & - f(x^*) + M_1 \cdot \| \nabla f(x^{k+1}) - g^{k+1} \|^2 + M_2 \cdot \frac{1}{n} \sum\limits_{j=1}^n \|\nabla f_j (x^{k+1}) - y_j^{k+1} \|^2 \Bigg]
    \\\leq& 
    (1 - \eta_k) \EE{f(x^k) - f(x^*)} 
    \\&
    +\left(1 - \frac{\lambda}{ 2}\right) M_1 \cdot \EE{\| \nabla f(x^{k}) - g^{k} \|^2}
    + \left(1 - \frac{b}{4n}\right) M_2\cdot \frac{1}{n} \sum\limits_{j=1}^n \EE{\|\nabla f_j (x^{k}) - y_j^{k} \|^2}
    \\&
    + \frac{4\alpha \tilde L^2}{b L \lambda} \left( 1 + \frac{8 \lambda^2 n^2}{b^2} \right)  \EE{\| x^{k+1} - x^{k}\|^2}
    + \frac{L \eta_k^2}{\alpha} \EE{\|s^{k} - x^* \|^2}
    + \frac{L \eta_k^2}{2} \EE{\|s^{k} - x^k \|^2}
    \\\leq& 
    (1 - \eta_k) \EE{f(x^k) - f(x^*)}
    \\&
    +\left(1 - \frac{\lambda}{ 2}\right) M_1 \cdot \EE{\| \nabla f(x^{k}) - g^{k} \|^2} 
    + \left(1 - \frac{b}{4n}\right) M_2\cdot \frac{1}{n} \sum\limits_{j=1}^n \EE{\|\nabla f_j (x^{k}) - y_j^{k} \|^2}
    \\&
    + \frac{4\alpha \tilde L^2 \eta_k^2}{b L \lambda} \left( 1 + \frac{8 \lambda^2 n^2}{b^2} \right)  \EE{\| s^{k} - x^{k}\|^2}
    + \frac{L \eta_k^2}{\alpha}\EE{\|s^{k} - x^* \|^2}
    + \frac{L \eta_k^2}{2} \EE{\|s^{k} - x^k \|^2}.
\end{align*}
With Assumption \ref{as:set} on the diameter $D$ of $\X$, we get
\begin{align*}
    \E \Bigg[f(x^{k+1}) &- f(x^*) + M_1 \cdot \| \nabla f(x^{k+1}) - g^{k+1} \|^2 + M_2 \cdot \frac{1}{n} \sum\limits_{j=1}^n \|\nabla f_j (x^{k+1}) - y_j^{k+1} \|^2\Bigg]
    \\\leq& 
    (1 - \eta_k) \EE{f(x^k) - f(x^*)} 
    \\&
    +\left(1 - \frac{\lambda}{ 2}\right) M_1 \cdot \EE{\| \nabla f(x^{k}) - g^{k} \|^2}
    + \left(1 - \frac{b}{4n}\right) M_2\cdot \frac{1}{n} \sum\limits_{j=1}^n \EE{\|\nabla f_j (x^{k}) - y_j^{k} \|^2}
    \\&
    + LD^2 \eta_k^2 \left( \frac{1}{2} + \frac{1}{\alpha} +  \frac{4\alpha \tilde L^2}{b L^2 \lambda} \left[ 1 + \frac{8 \lambda^2 n^2}{b^2} \right] \right).
\end{align*}
The choices of $\lambda = \tfrac{b}{2n}$ and $\alpha = \tfrac{b L}{5 \tilde L \sqrt{n}}$ provides
\begin{align*}
    \E \Bigg[f(x^{k+1}) &- f(x^*) + M_1 \cdot \| \nabla f(x^{k+1}) - g^{k+1} \|^2 + M_2 \cdot \frac{1}{n} \sum\limits_{j=1}^n \|\nabla f_j (x^{k+1}) - y_j^{k+1} \|^2 \Bigg]
    \\\leq& 
    (1 - \eta_k) \EE{f(x^k) - f(x^*)} 
    \\&
    +\left(1 - \frac{b}{4n}\right) M_1 \cdot \EE{\| \nabla f(x^{k}) - g^{k} \|^2}
    + \left(1 - \frac{1}{4n}\right) M_2\cdot \frac{1}{n} \sum\limits_{j=1}^n \EE{\|\nabla f_j (x^{k}) - y_j^{k} \|^2}
    \\&
    + LD^2 \eta_k^2 \left( \frac{1}{2} + \frac{10 \tilde L \sqrt{n}}{L}\right)
    \\\leq& 
    \max\left\{ 1 - \eta_k,  1 - \frac{b}{4n}\right\} \mathbb{E}\Bigg[f(x^k) - f(x^*) + M_1 \cdot \| \nabla f(x^{k}) - g^{k} \|^2 
    \\
    &+ M_2\cdot \frac{1}{n} \sum\limits_{j=1}^n \|\nabla f_j (x^{k}) - y_j^{k} \|^2 \Bigg]
    + LD^2 \eta_k^2 \left( \frac{1}{2} + \frac{10 \tilde L \sqrt{n}}{Lb}\right).
\end{align*}
If we choose $\eta_k \leq \frac{b}{4n}$, we have
\begin{align*}
    \E \Bigg[f(x^{k+1}) &- f(x^*) + M_1 \cdot \| \nabla f(x^{k+1}) - g^{k+1} \|^2 + M_2 \cdot \frac{1}{n} \sum\limits_{j=1}^n \|\nabla f_j (x^{k+1}) - y_j^{k+1} \|^2\Bigg]
    \\\leq& 
    \left( 1 - \eta_k\right) \EE{f(x^k) - f(x^*) + M_1 \cdot \| \nabla f(x^{k}) - g^{k} \|^2 
    + M_2\cdot \frac{1}{n} \sum\limits_{j=1}^n \|\nabla f_j (x^{k}) - y_j^{k} \|^2 }
    \\&
    + LD^2 \eta_k^2 \left( \frac{1}{2} + \frac{10 \tilde L \sqrt{n}}{Lb}\right).
\end{align*}
It remains to use Lemma \ref{lem:tech3} with $c = LD^2\left(\tfrac{1}{2} 
    +  \tfrac{10\tilde L \sqrt{n}}{L b}\right)$, $d = \tfrac{4n}{b}$ and obtain
\begin{align*}
    \EE{f(x^{K}) - f(x^*)} 
    =&
    \cO\Bigg( \frac{n}{b} \bigg[f(x^{0}) - f(x^*) + M_1 \cdot \| \nabla f(x^{0}) - g^{0} \|^2 
    \\&\quad \quad 
    + M_2 \cdot \frac{1}{n} \sum\limits_{j=1}^n \|\nabla f_j (x^{0}) - y_j^{0} \|^2 \bigg] \exp\left(-\frac{bK}{8n} \right)   + \left[1 
    +  \frac{\tilde L \sqrt{n}}{L b}\right] \frac{LD^2}{K}\Bigg).
\end{align*}
Finally, we substitute $g^0 = \nabla f(x^0)$, $y_j^0 = \nabla f_j (x^0)$ and get
\begin{align*}
    \EE{f(x^{K}) - f(x^*)} 
    =&
    \cO\Bigg( \frac{n}{b} \left[f(x^{0}) - f(x^*) \right] \exp\left(-\frac{bK}{8n} \right)   + \left[1 
    +  \frac{\tilde L \sqrt{n}}{L b}\right] \frac{LD^2}{K}\Bigg).
\end{align*}
This completes the proof.
\EndProof

%%%%%%%%%%%%%%%%%%%%%%%%%%%%%%%%%%%%%%%%%%%%%%%%%%%%%%%%%%%%%%%%%%%%%%%%
%%%%%%%%%%%%%%%%%%%%%%%%%%%%%%%%%%%%%%%%%%%%%%%%%%%%%%%%%%%%%%%%%%%%%%%%
%%%%%%%%%%%%%%%%%%%%%%%%%%%%%%%%%%%%%%%%%%%%%%%%%%%%%%%%%%%%%%%%%%%%%%%%
%%%%%%%%%%%%%%%%%%%%%%%%%%%%%%%%%%%%%%%%%%%%%%%%%%%%%%%%%%%%%%%%%%%%%%%%

\subsection{Proof of Theorem \ref{th:main4}} \label{sec:proof4}

\begin{theorem}[Theorem \ref{th:main4}] \label{th:4}
Let $\{x^k\}_{k\geq0}$ denote the iterates of Algorithm \ref{alg:fw_zerosarah} for solving problem \eqref{eq:main_problem}, which satisfies Assumptions \ref{as:lip},\ref{as:set}. Let $x^*$ be the global (may be not unique) minimizer of $f$ on $\X$. Then, if we choose $\eta_k = \frac{1}{\sqrt{K}}$ and $\lambda = \tfrac{b}{2n}$, we have the following convergence:
\begin{align*}
    \EE{\min_{0\leq k \leq K-1}\gap(x^k)}
    = \mathcal{O}\left(
    \frac{f(x^0) - f(x^*)}{\sqrt{K}} + \frac{LD^2}{\sqrt{K}} \left[1 
    +  \frac{\tilde L \sqrt{n}}{Lb}\right]\right).
\end{align*}
\end{theorem}

\textbf{Proof:} Since lines \ref{alg2:line2} and \ref{alg2:line3} of Algorithms \ref{alg:fw_sarah} and \ref{alg:fw_zerosarah} are the same, we start the proof from \eqref{eq:temp3}. Multiplying \eqref{eq:temp4} by the positive constant $M_1$, \eqref{eq:temp5} by the positive constant $M_2$ ($M_1$, $M_2$ will be defined below) and summing with \eqref{eq:temp3}, we have
\begin{align*}
    \eta_k \E\big[\max_{x \in \X}&\<\nabla f(x^k), x^k - x >\big]
    \\\leq&
    (1 - \eta_k) \EE{f(x^k) - f(x^*)}
    \\&
    +\left(1 - \lambda + \frac{\alpha}{ M_1 L}\right) M_1 \cdot \EE{\| \nabla f(x^{k}) - g^{k} \|^2} 
    \\&
    + \left(1 - \frac{b}{2n} + \frac{2 M_1 \lambda^2}{M_2 b}\right) M_2\cdot \frac{1}{n} \sum\limits_{j=1}^n \EE{\|\nabla f_j (x^{k}) - y_j^{k} \|^2}
    \\& 
    -\EE{f(x^{k+1}) - f(x^*) + M_1 \cdot \| \nabla f(x^{k+1}) - g^{k+1} \|^2 + M_2 \cdot \frac{1}{n} \sum\limits_{j=1}^n \|\nabla f_j (x^{k+1}) - y_j^{k+1} \|^2 }
    \\&
    +  \frac{2( M_1 + n M_2) \tilde L^2 \eta_k}{b}  \EE{\| s^{k} - x^{k}\|^2}
    + \frac{L \eta_k^2}{\alpha}\EE{\max_{x \in \X} \|s^{k} - x \|^2}
    + \frac{L \eta_k^2}{2} \EE{\|s^{k} - x^k \|^2}.
\end{align*}
With $M_1 = \tfrac{\alpha}{\lambda L}$ and $M_2 = \tfrac{4M_1 \lambda^2 n}{b^2}$, we obtain
\begin{align*}
    \eta_k \E\big[\max_{x \in \X}&\<\nabla f(x^k), x^k - x >\big]
    \\\leq&
    \EE{f(x^k) - f(x^*)
    +M_1 \cdot \| \nabla f(x^{k}) - g^{k} \|^2
    +M_2\cdot \frac{1}{n} \sum\limits_{j=1}^n \|\nabla f_j (x^{k}) - y_j^{k} \|^2}
    \\& 
    -\EE{f(x^{k+1}) - f(x^*) + M_1 \cdot \| \nabla f(x^{k+1}) - g^{k+1} \|^2 + M_2 \cdot \frac{1}{n} \sum\limits_{j=1}^n \|\nabla f_j (x^{k+1}) - y_j^{k+1} \|^2 }
    \\&
    + \frac{2\alpha \tilde L^2 \eta_k^2}{b L \lambda} \left( 1 + \frac{4 \lambda^2 n^2}{b^2} \right) \EE{\| s^{k} - x^{k}\|^2}
    + \frac{L \eta_k^2}{\alpha}\EE{\|s^{k} - x^* \|^2}
    + \frac{L \eta_k^2}{2} \EE{\|s^{k} - x^k \|^2}.
\end{align*}
Assumption \ref{as:set} on the diameter $D$ of $\X$ gives
\begin{align*}
   \eta_k \E\big[\max_{x \in \X}&\<\nabla f(x^k), x^k - x >\big]
    \\\leq&
    \EE{f(x^k) - f(x^*)
    +M_1 \cdot \| \nabla f(x^{k}) - g^{k} \|^2
    +M_2\cdot \frac{1}{n} \sum\limits_{j=1}^n \|\nabla f_j (x^{k}) - y_j^{k} \|^2}
    \\& 
    -\EE{f(x^{k+1}) - f(x^*) + M_1 \cdot \| \nabla f(x^{k+1}) - g^{k+1} \|^2 + M_2 \cdot \frac{1}{n} \sum\limits_{j=1}^n \|\nabla f_j (x^{k+1}) - y_j^{k+1} \|^2 }
    \nonumber\\&
    + LD^2 \eta_k^2 \left(\frac{1}{2} + \frac{1}{\alpha}
    +  \frac{2\alpha \tilde L^2}{b L^2 \lambda} \left[ 1 + \frac{4 \lambda^2 n^2}{b^2} \right]\right).
\end{align*}
The choices of $\lambda = \tfrac{b}{2n}$ and $\alpha = \tfrac{b L}{3 \tilde L \sqrt{n}}$ provides
\begin{align*}
    \eta_k \E\big[\max_{x \in \X}&\<\nabla f(x^k), x^k - x >\big]
    \\\leq&
    \EE{f(x^k) - f(x^*)
    +M_1 \cdot \| \nabla f(x^{k}) - g^{k} \|^2
    +M_2\cdot \frac{1}{n} \sum\limits_{j=1}^n \|\nabla f_j (x^{k}) - y_j^{k} \|^2}
    \\& 
    -\EE{f(x^{k+1}) - f(x^*) + M_1 \cdot \| \nabla f(x^{k+1}) - g^{k+1} \|^2 + M_2 \cdot \frac{1}{n} \sum\limits_{j=1}^n \|\nabla f_j (x^{k+1}) - y_j^{k+1} \|^2 }
    \nonumber\\&
    + LD^2 \eta_k^2 \left( \frac{1}{2} + \frac{6 \tilde L \sqrt{n}}{Lb}\right).
\end{align*}
Summing over all $k$ from $0$ to $K-1$, taking $\eta_k = \frac{1}{\sqrt{K}}$ and dividing both sides by 
\begin{align*}
    \E \Bigg[\frac{1}{K} & \sum\limits_{k=0}^{K-1} \max_{x \in \X} \<\nabla f(x^k), x^k - x >\Bigg]
    \\\leq& 
    \frac{1}{\sqrt{K}} \cdot \left[f(x^0) - f(x^*)
    +M_1 \cdot \| \nabla f(x^{0}) - g^{0} \|^2
    +M_2\cdot \frac{1}{n} \sum\limits_{j=1}^n \|\nabla f_j (x^{0}) - y_j^{0} \|^2 \right]
    \\&
    + \frac{LD^2}{\sqrt{K}} \left(\frac{1}{2} 
    +  \frac{6 \tilde L \sqrt{n}}{Lb} \right).
\end{align*}
Finally, we substitute $g^0 = \nabla f(x^0)$, $y_j^0 = \nabla f_j (x^0)$ and get
\begin{align*}
    \EE{\frac{1}{K}\sum\limits_{k=0}^{K-1} \max_{x \in \X} \<\nabla f(x^k), x^k - x >}
    \leq 
    \frac{f(x^0) - f(x^*)}{\sqrt{K}} + \frac{LD^2}{\sqrt{K}} \left(\frac{1}{2} 
    +  \frac{6 \tilde L \sqrt{n}}{Lb} \right).
\end{align*}
The definition of \eqref{eq:gap} finishes the proof.
\EndProof

\section{Additional Comments}

\subsection{On the Convergence Criterion in \citep{pmlr-v80-qu18a, pmlr-v119-gao20b}} \label{sec:crit}
In Table \ref{tab:comparison0}, we indicates that the papers by \cite{pmlr-v80-qu18a, pmlr-v119-gao20b} considers $\| \nabla f(x)\|^2$ as a convergence criterion. But the authors actually use a more complex criterion $\| G(x, \nabla f(x), \gamma) \|^2$, where 
$
G(x, \nabla f(x), \gamma) = \frac{1}{\gamma} (x - \psi (x, \nabla f(x), \gamma) ) \quad \text{with} \quad  \psi (x, \nabla f(x), \gamma) = \arg\min_{y \in C} \left( \langle \nabla f(x), y\rangle  + \frac{1}{2\gamma} \| x - y \|^2 \right).
$
Let us simplify this criterion. If $C$ is large enough, one can assume that we work on unconstrained setting and $\arg\min_{y \in C} \left(\langle \nabla f(x), y \rangle + \frac{1}{2\gamma} \| x - y \|^2 \right) \approx \arg\min_{y \in R^d} \left(\langle \nabla f(x), y \rangle + \frac{1}{2\gamma} \| x - y \|^2 \right)$, i.e. $\psi (x, \nabla f(x), \gamma) = x - \gamma \nabla f(x)$. Therefore, we get that $ G(x, \nabla f(x), \gamma) \approx \nabla f(x)$ and $\| G(x, \nabla f(x), \gamma) \|^2 \approx \|\nabla f(x) \|^2$.
As we noted in Table \ref{tab:comparison0}, to avoid discrepancies with the lower bounds from \citep{pmlr-v139-li21a}, we slightly modified the result of \cite{pmlr-v80-qu18a, pmlr-v119-gao20b}. In more details, following \cite{pmlr-v139-li21a}, we assume that want to achieve $\| G(x, \nabla f(x), \gamma) \|^2 \sim \varepsilon^2$. In the original papers \citep{pmlr-v80-qu18a, pmlr-v119-gao20b}, the authors uses $\| G(x, \nabla f(x), \gamma) \|^2 \sim \varepsilon$.

\subsection{Incorrect Proof of Theorem 4 from \citep{reddi2016stochastic}} \label{sec:error}

As we noted in Section \ref{sec:related}, the paper provides another algorithm (Algorithm 4). This method is a modification of the SAGA technique. The proof of convergence of this algorithm, in our opinion, contains a mistake. The authors introduce an additional technical sequence:
\begin{equation}
    \label{eq:c_t}
   c_t = (1 - \rho) c_{t+1} + \frac{L D \gamma \sqrt{n}}{\sqrt{b}} \quad \text{with} \quad c_T = 0,
\end{equation}
and claim that the following estimate is valid:
\begin{equation*}
   \sum\limits_{t=1}^T c_t \leq \frac{L D \gamma \sqrt{n}}{\rho \sqrt{b}}.
\end{equation*}
But if we consider the simplest case with $\rho = 1$, we have that 
\begin{equation*}
   c_t = \frac{L D \gamma \sqrt{n}}{\sqrt{b}} \quad \text{ and } \quad \sum\limits_{t=1}^T c_t \leq \frac{L D \gamma \sqrt{n}}{\sqrt{b}} \cdot (T-1),
\end{equation*}
which is larger than the authors' estimate: $\frac{L D \gamma \sqrt{n}}{\sqrt{b}}$.

Let us try to correct this error. Running the recursion \eqref{eq:c_t}, we get for all $t = 0, \ldots, (T-1)$
\begin{equation*}
    c_t = \frac{L D \gamma \sqrt{n}}{\sqrt{b}} \cdot \sum\limits_{i = T}^{t+1} (1-\rho)^{T - i} \leq \frac{L D \gamma \sqrt{n}}{\rho \sqrt{b}}
\end{equation*}
And then,
\begin{equation*}
   \sum\limits_{t=1}^T c_t \leq \frac{L D \gamma \sqrt{n}}{\rho \sqrt{b}} \cdot (T-1).
\end{equation*}
With (13) from \citep{reddi2016stochastic}: $\rho \geq \tfrac{b}{2n}$, one can obtain
\begin{equation*}
    \sum\limits_{t=1}^T c_t \leq \frac{L D \gamma n^{3/2}}{b^{3/2}} \cdot (T-1).
\end{equation*}
The result is the following estimate
\begin{align*}
    \eta \sum\limits_{k=0}^{K-1}\EE{\max_{x \in \X} \<\nabla f(x^k), x^k - x >}
    \leq
    f(x^0) - f(x^*) 
    + LD^2 \eta^2 K \left(1 
    +  \frac{n^{3/2}}{b^{3/2}}\right).
\end{align*}
If we take $\eta = \frac{1}{\sqrt{K}}$ and divide both sides by $\sqrt{K}$, then 
\begin{align*}
    \frac{1}{K} \sum\limits_{k=0}^{K-1}\EE{\max_{x \in \X} \<\nabla f(x^k), x^k - x >}
    \leq
    \frac{f(x^0) - f(x^*)}{\sqrt{K}} 
    + \frac{LD^2}{\sqrt{K}} \left(1 
    +  \frac{n^{3/2}}{b^{3/2}}\right).
\end{align*}
The definition of \eqref{eq:gap} gives
\begin{align*}
    \EE{\min_{0\leq k \leq K-1}\gap(x^k)}
    \leq
    \frac{f(x^0) - f(x^*)}{\sqrt{K}} 
    + \frac{LD^2}{\sqrt{K}} \left(1 
    +  \frac{n^{3/2}}{b^{3/2}}\right).
\end{align*}
Therefore, the following number of the stochastic oracle calls is needed to achieve the accuracy $\varepsilon$:
\begin{align*} 
    \mathcal{O}\left( b\left[\frac{f(x^0) - f(x^*)}{\varepsilon}\right]^2 + \left[\frac{LD^2}{\varepsilon}\right]^2 \left[b +  \frac{n^3}{b^2}\right]\right) .
\end{align*}
With the optimal choice of $b = n$, we get
\begin{align*} 
    \mathcal{O}\left( n\left[\frac{f(x^0) - f(x^*)}{\varepsilon}\right]^2 + n\left[\frac{LD^2}{\varepsilon}\right]^2\right) .
\end{align*}

\subsection{Big Batches and LMO Complexities} \label{sec:batches}

Here we study the use of large batches to obtain better estimates on LMO by using of lower bounds. For this purposes, we need to introduce a class of algorithms for which the lower bounds will be valid. Since we work with projection-free methods, the following definition takes this into account.

\begin{definition}
    We have local memory $\mathcal{M}$ with initialization $\mathcal{M} = \{0\}$.
    In addition, we also have an auxiliary buffer $\mathcal{H}$, which is also initially equal to $\{0\}$. These memory $\mathcal{M}$ and buffer $\mathcal{H}$ are updated as follows. 
    
    $\bullet$ One can sample uniformly and independently batch $S$ of any size $b$ (if $b=n$ -- it means that we call the full gradient) from $\{f_{i}\}$ at some point $x \in \mathcal{M}$, compute stochastic gradient $g = \sum_{j \in S} \nabla f_j (x)$ and adds it to the buffer $\mathcal{H}$ as linear combination of existing vectors in the buffer: 
    $$ \mathcal{H} = \text{span} \{\mathcal{H}, g\}. $$ 
    We can repeat this operation with different batches with any sizes and different points $x \in \mathcal{M}$.
    
    $\bullet$ Using information in the buffer, we can update our memory $\mathcal{M}$ by adding to $\mathcal{M}$ a finite number of points $x'$, satisfying 
    $$ x' \in \text{span} \{x, \text{LMO} (h, \mathcal{X}) \}, $$ 
    where we can take any $x \in \mathcal{M}$, $h \in \mathcal{H}$ and LMO is the linear minimization oracle.

    $\bullet$ The final global output is calculated as $x \in \mathcal{M}$.
\end{definition} 
 
The next step in constructing lower bounds is to create a "bad" problem on which all methods perform poorly. In our case, the problem consists of two parts: a function $f$ with its division by $f_i$, and an optimization set $\mathcal{C}$. The functions can be taken from works on lower bounds for the unconstrained case, e.g., in the convex \cite{han2024lower} and non-convex \cite{fang2018spider, pmlr-v139-li21a} setups. As an optimization set we choose $\ell_1$-ball with $0$ center and size $R=1$. For this ball: 
$$ 
\text{LMO} (h, \mathcal{X}) = - \text{sign} (h_i) e_i \quad \text{with} \quad i = \arg\max_j |g_j|. 
$$ 
If the solution of argmax is not unique, we choose the smallest one.

The essence of the lower bounds is classical \cite{nesterov2013introductory} -- how the final output close to the real solution is measured in the number of non-zero coordinates in the output. On non-zero coordinates we can (in the best case) get a number corresponding to the real solution, and  on zero coordinates we cannot. How it works for unconstrained optimization methods: each gradient call "open" a new non-zero coordinate, then as many gradients we compute --  that is how many non-zero coordinates we have.

In our case, we have also LMO in the update rule. If $\mathcal{H} = \text{span}\{ e_1, \ldots, e_k \}$, then we can guarantee that LMO for any vector from $\mathcal{H}$ lies in $\text{span}\{ e_1, \ldots, e_k \}$, since LMO for $\ell_1$-ball is a maximum absolute value of coordinates and if there are two maximums we choose the smaller one. It means that LMO for our set does give any progress in terms of non-zero coordinates, but we can make this progress in $\mathcal{H}$ and using LMO we transfer a new non-zero coordinate to $\mathcal{M}$. Therefore, we need to understand how batching affects $\mathcal{H}$, and then we immediately understand how it affects LMO.

\begin{proposition} \label{prop:1}
After $K$ LMO computations, only the first $K$ coordinates of the global output can be non-zero while the rest of the $d-K$ coordinates are strictly equal to zero.
\end{proposition}
\textbf{Proof:}
We compute full gradients in the best case, these full gradients add one more non-zero coordinate in $\mathcal{H}$ and then using LMO in $\mathcal{M}$.
\EndProof

\begin{proposition} \label{prop:2}
If between LMO computations and $\mathcal{M}$ updates we collect batch size of $1$ in only one point, then after $K$ LMO computations, in expectation only the first $\frac{K}{n}$ coordinates of the global output can be non-zero while the rest of the $d-\frac{K}{n}$ coordinates are strictly equal to zero.
\end{proposition}
\textbf{Proof:}
Now we cannot compute full gradients, we compute only batch of size 1. The key problem with such a batch is that we set $f_i$ in such a way that different parts of the whole $f$ are stored on different $f_i$, and only one of all $f_i$s can increase the number of non-zero coordinates \cite{fang2018spider,pmlr-v139-li21a,han2024lower}. But by virtue of the fact that we choose this function randomly, then we get a new non-zero coordinate with probability $1/n$. It turns out that we can call LMO for nothing and not get new non-zero coordinates in the output of $\mathcal{M}$.
\EndProof

Propositions \ref{prop:1} and \ref{prop:2} show that with different batching we can get different number of non-zero coordinates. These propositions are not a rigorous justification for the use of large (e.g., $\sqrt{n}$) or even full batches, but they provide a semi-formalized intuition for why this might be the case.

\newpage

\section{Additional Experiments} \label{sec:add_exp}

Here we give experiments on logistic regression but unlike the main part we consider other sizes of the $\ell_1$-ball $R = 200, 20, 2$. These experiments also verify the superiority of our algorithms. 

\begin{figure*}[h!]
\begin{minipage}{0.24\textwidth}
  \centering
\includegraphics[width =  \textwidth ]{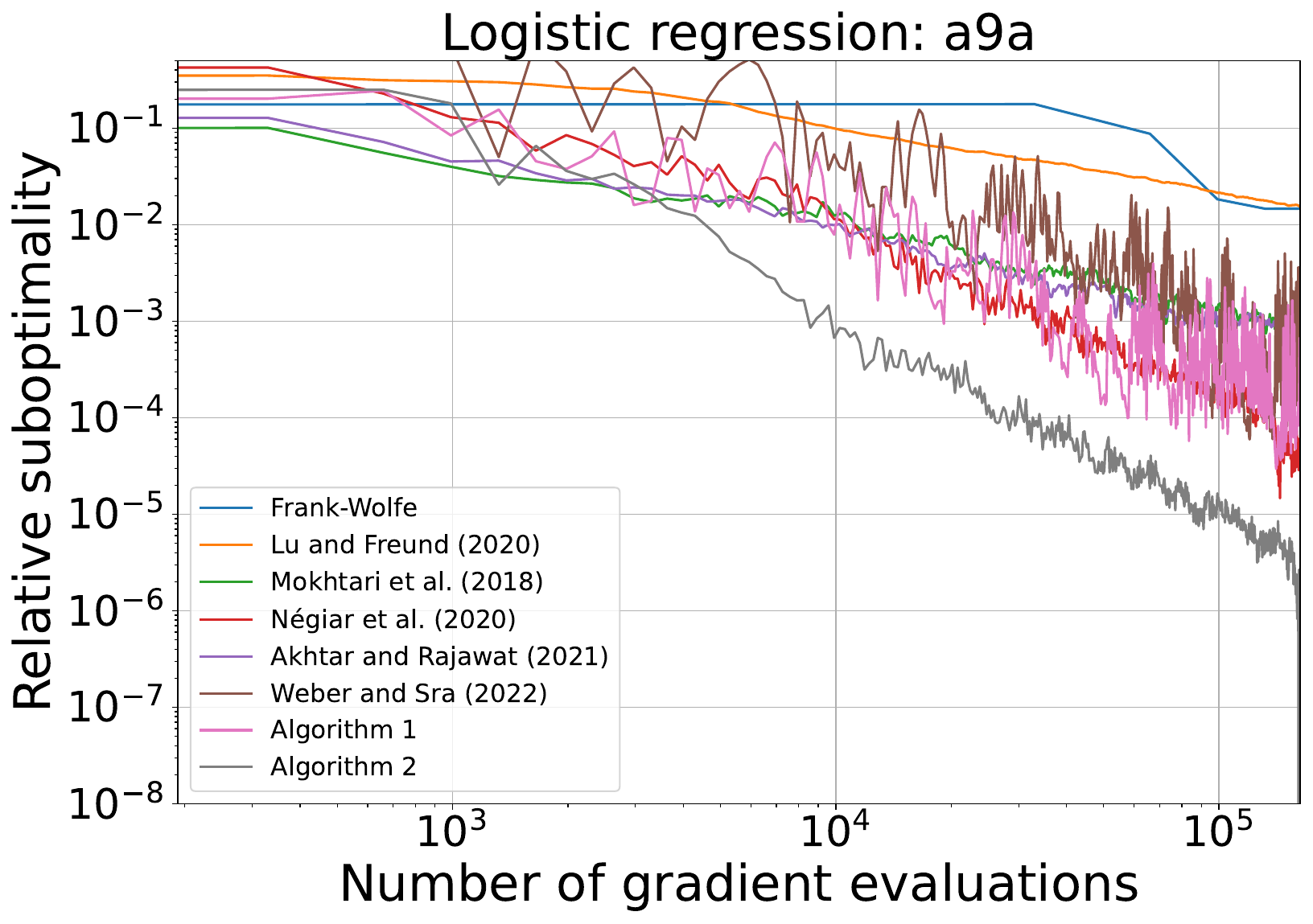}
\end{minipage}%
\begin{minipage}{0.24\textwidth}
  \centering
\includegraphics[width =  1\textwidth ]{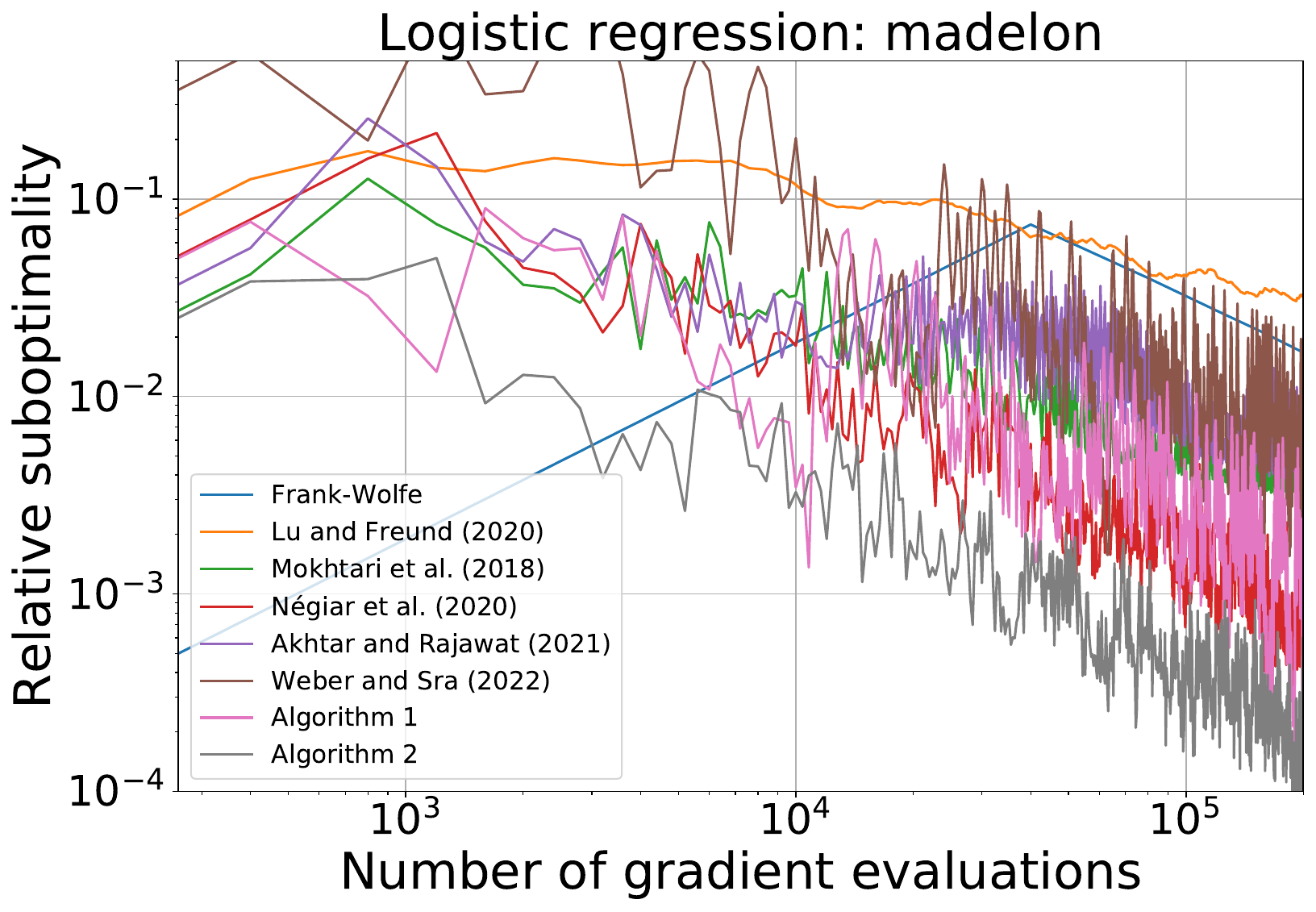}
\end{minipage}%
\begin{minipage}{0.24\textwidth}
  \centering
\includegraphics[width =  1\textwidth ]{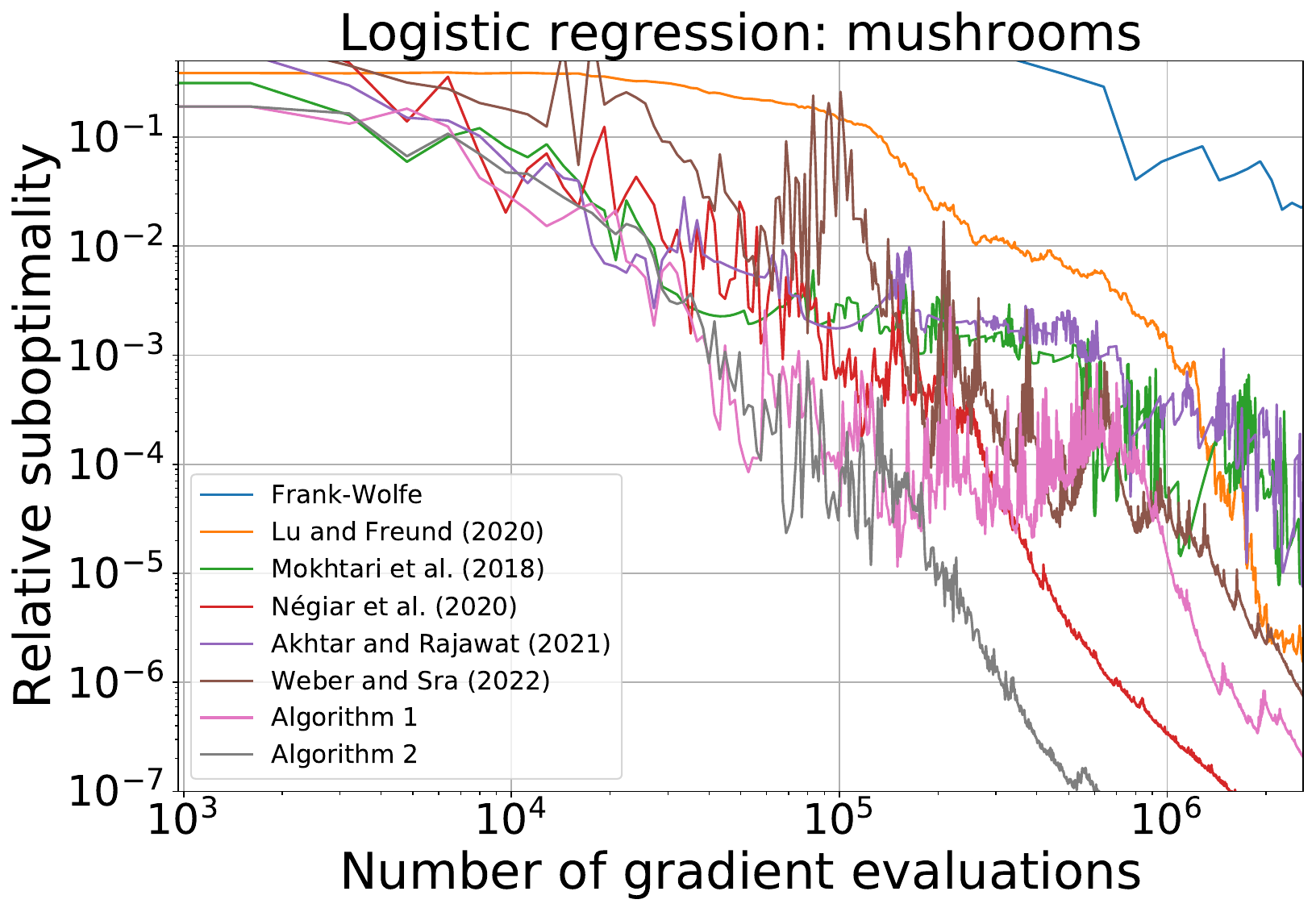}
\end{minipage}%
\begin{minipage}{0.24\textwidth}
  \centering
\includegraphics[width =  1\textwidth ]{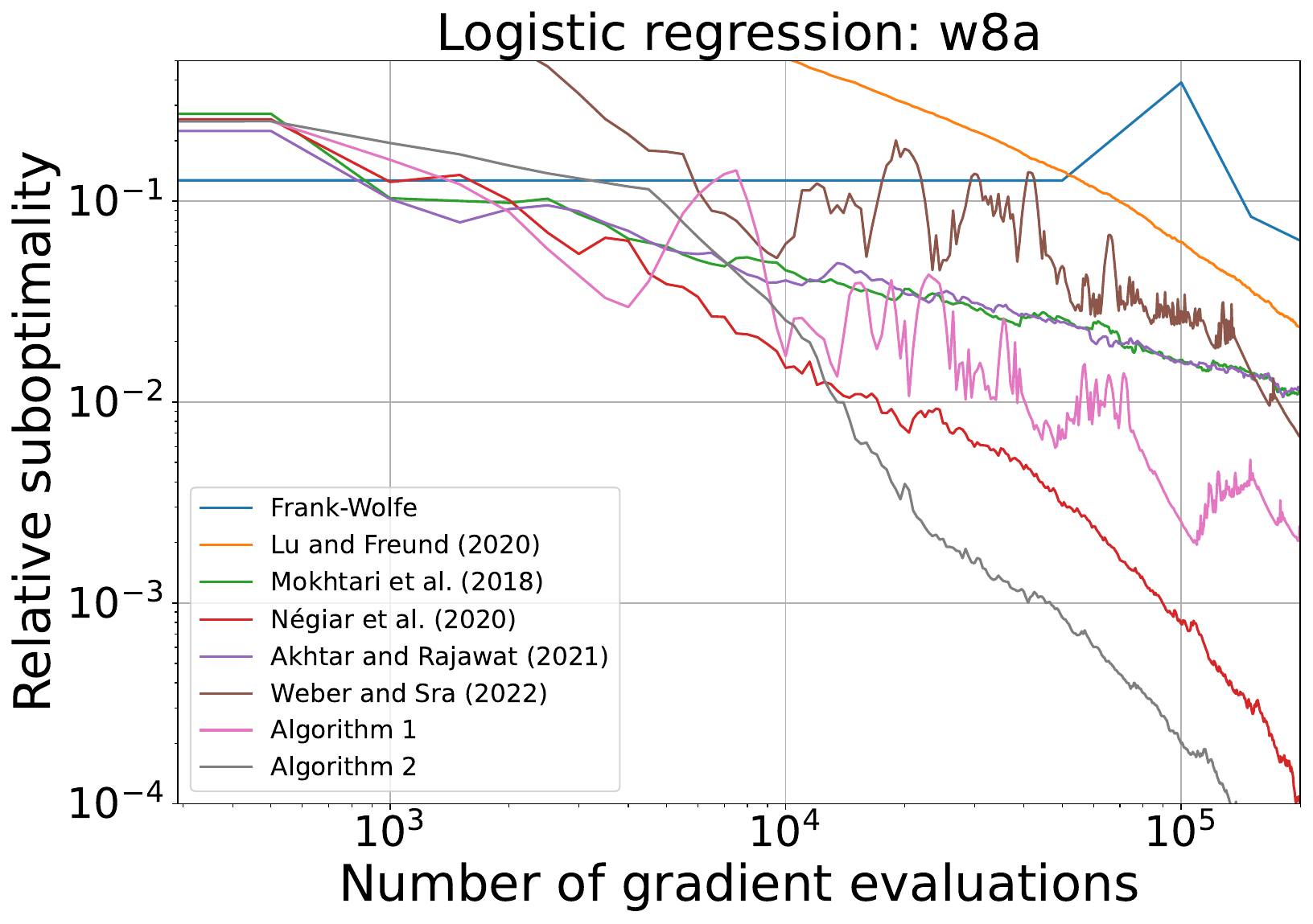}
\end{minipage}%
\\
\begin{minipage}{0.24\textwidth}
  \centering
(a) \texttt{a9a}
\end{minipage}%
\begin{minipage}{0.24\textwidth}
\centering
(b) \texttt{breast cancer}
\end{minipage}%
\begin{minipage}{0.24\textwidth}
\centering
  (c) \texttt{mushrooms}
\end{minipage}%
\begin{minipage}{0.24\textwidth}
\centering
  (d) \texttt{w8a}
\end{minipage}%
\vspace{-0.3cm}
\caption{Comparison of state-of-the-art projection free methods with small batches for \eqref{eq;ls} with $R = 200$. The comparison is made on the real datasets from LibSVM. The criterion is the number of full gradients computations.}
\end{figure*}

\begin{figure*}[h!]
\vspace{-0.3cm}
\begin{minipage}{0.24\textwidth}
  \centering
\includegraphics[width =  \textwidth ]{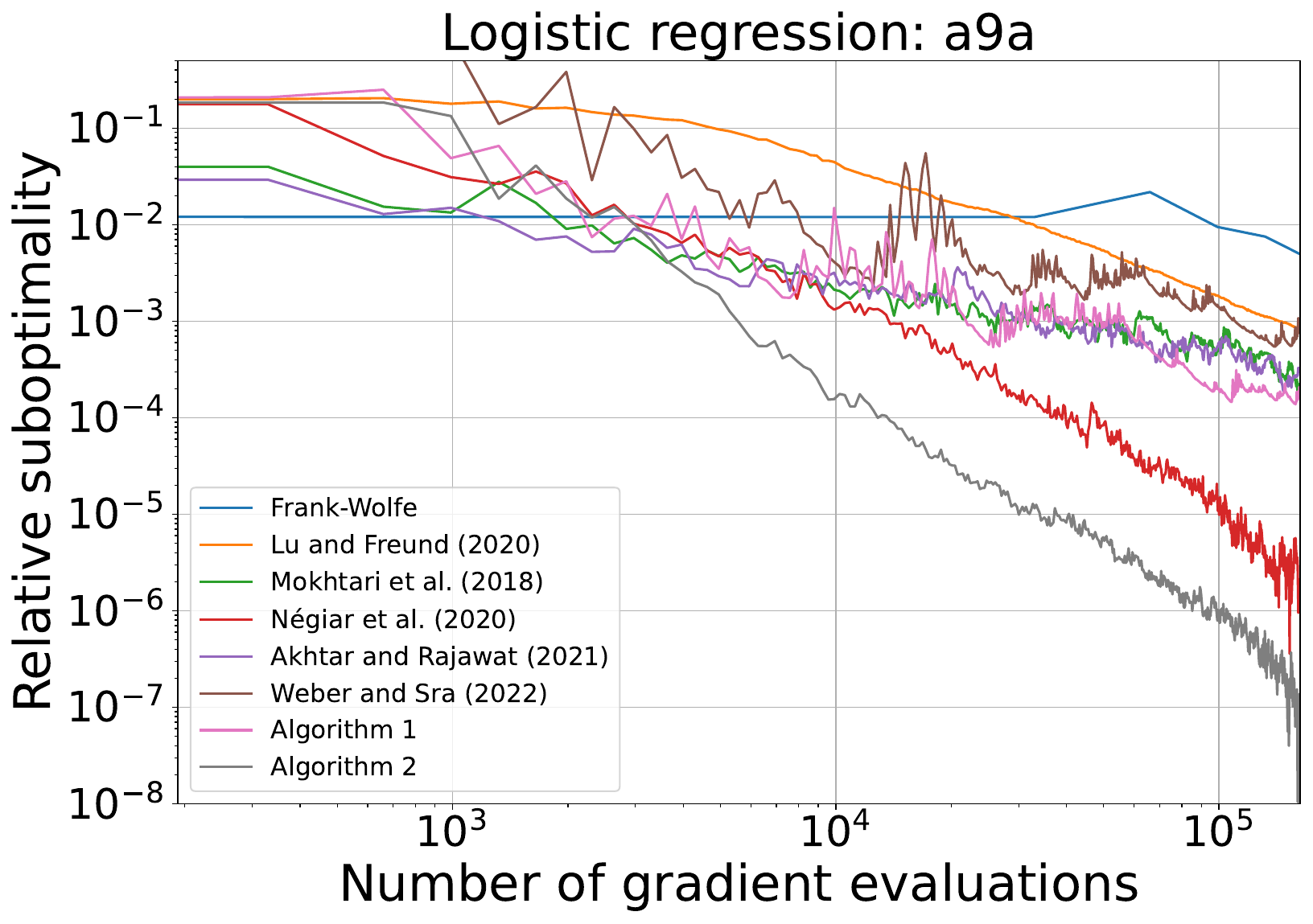}
\end{minipage}%
\begin{minipage}{0.24\textwidth}
  \centering
\includegraphics[width =  1\textwidth ]{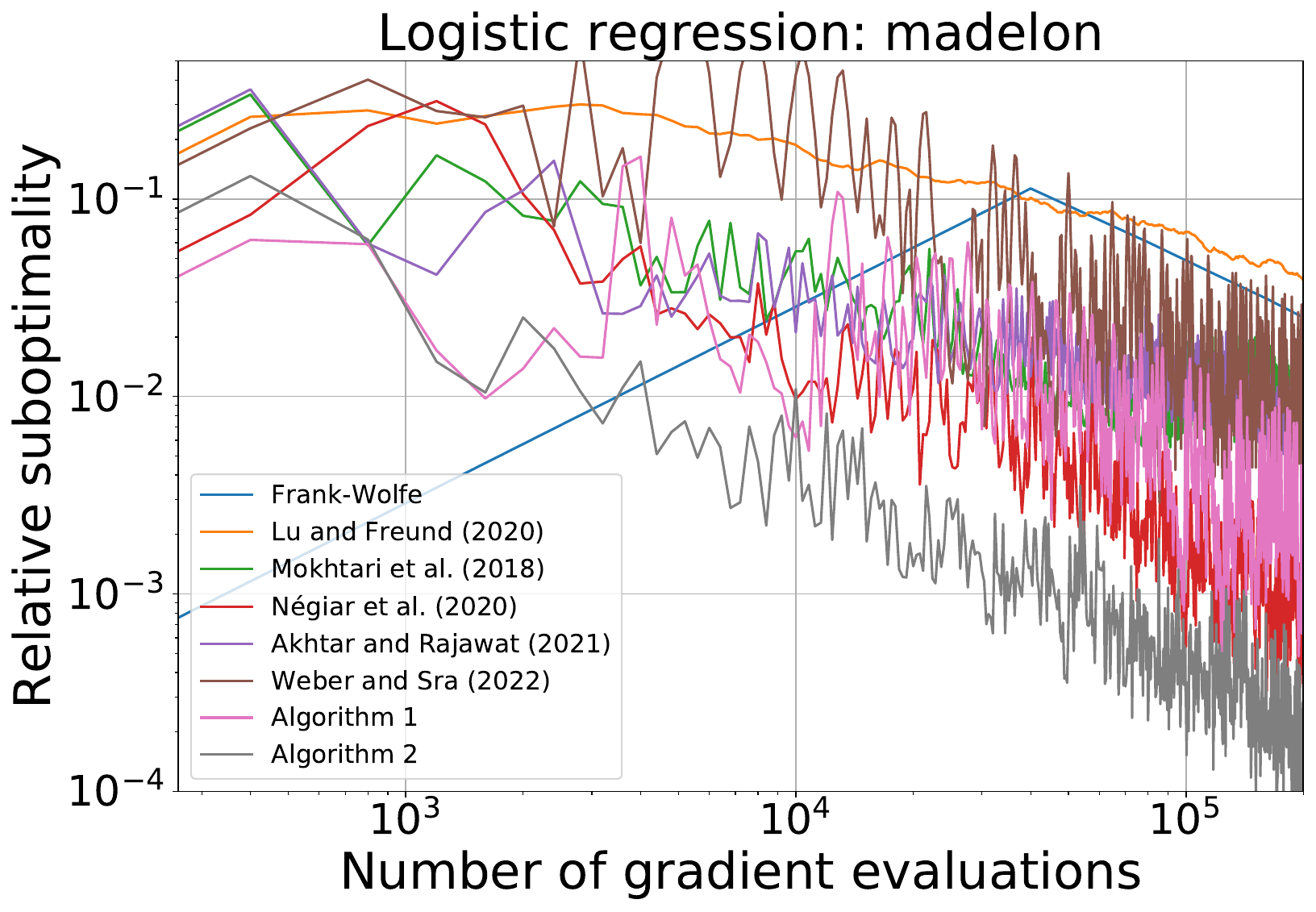}
\end{minipage}%
\begin{minipage}{0.24\textwidth}
  \centering
\includegraphics[width =  1\textwidth ]{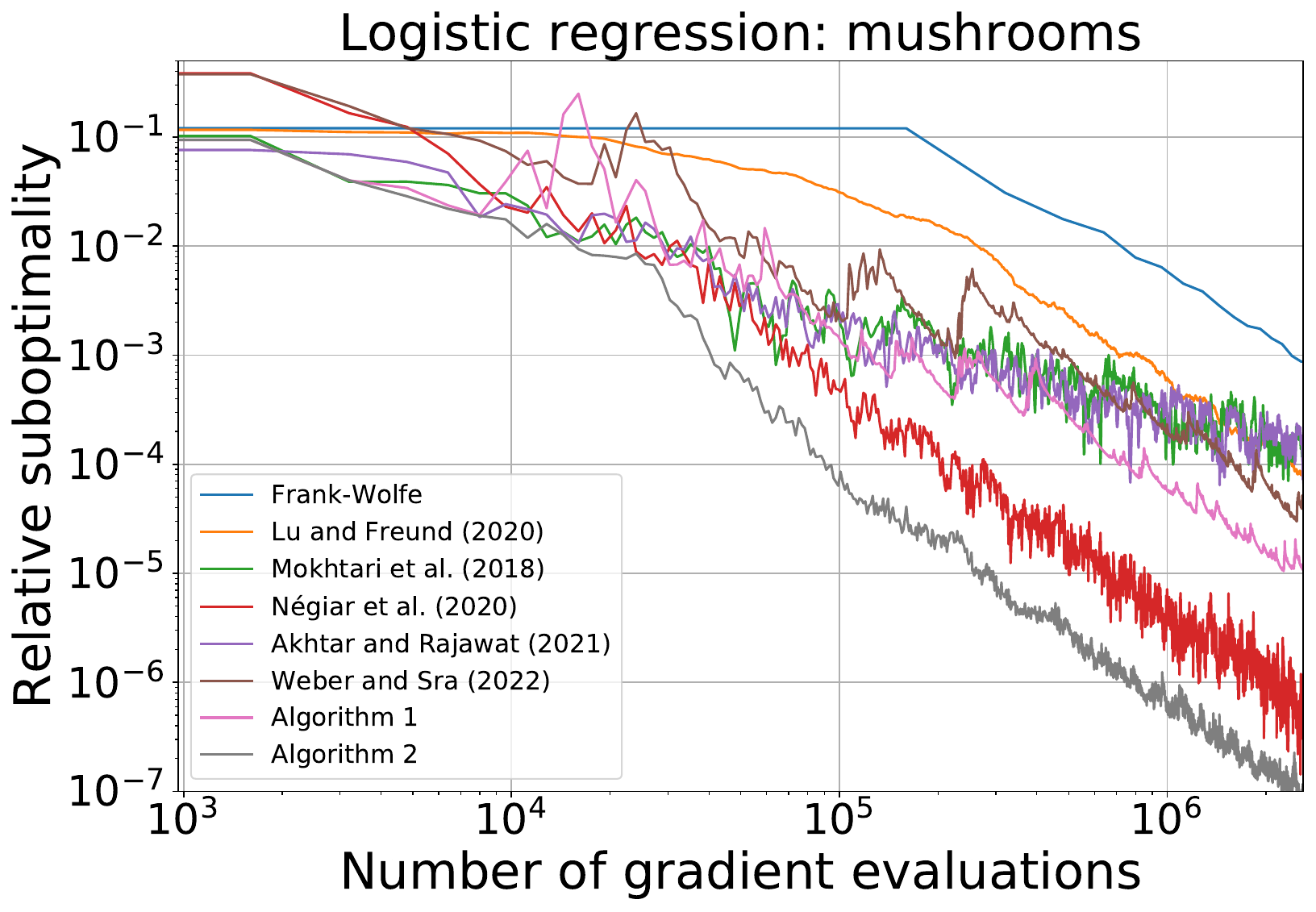}
\end{minipage}%
\begin{minipage}{0.24\textwidth}
  \centering
\includegraphics[width =  1\textwidth ]{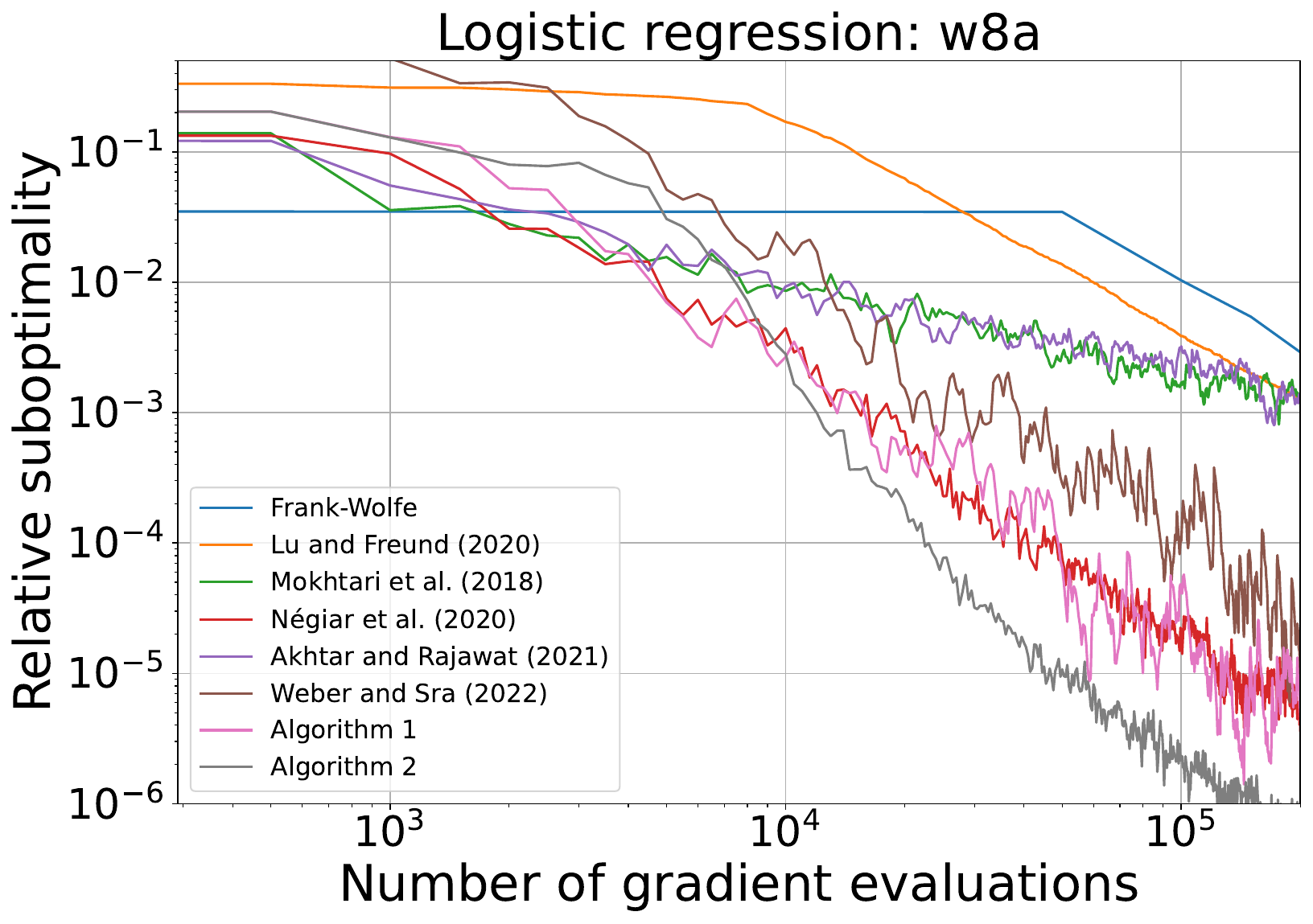}
\end{minipage}%
\\
\begin{minipage}{0.24\textwidth}
  \centering
(a) \texttt{a9a}
\end{minipage}%
\begin{minipage}{0.24\textwidth}
\centering
(b) \texttt{breast cancer}
\end{minipage}%
\begin{minipage}{0.24\textwidth}
\centering
  (c) \texttt{mushrooms}
\end{minipage}%
\begin{minipage}{0.24\textwidth}
\centering
  (d) \texttt{w8a}
\end{minipage}%
\vspace{-0.3cm}
\caption{Comparison of state-of-the-art projection free methods with small batches for \eqref{eq;ls} with $R = 20$. The comparison is made on the real datasets from LibSVM. The criterion is the number of full gradients computations.}
\end{figure*}

\begin{figure*}[h!]
\vspace{-0.3cm}
\begin{minipage}{0.24\textwidth}
  \centering
\includegraphics[width =  \textwidth ]{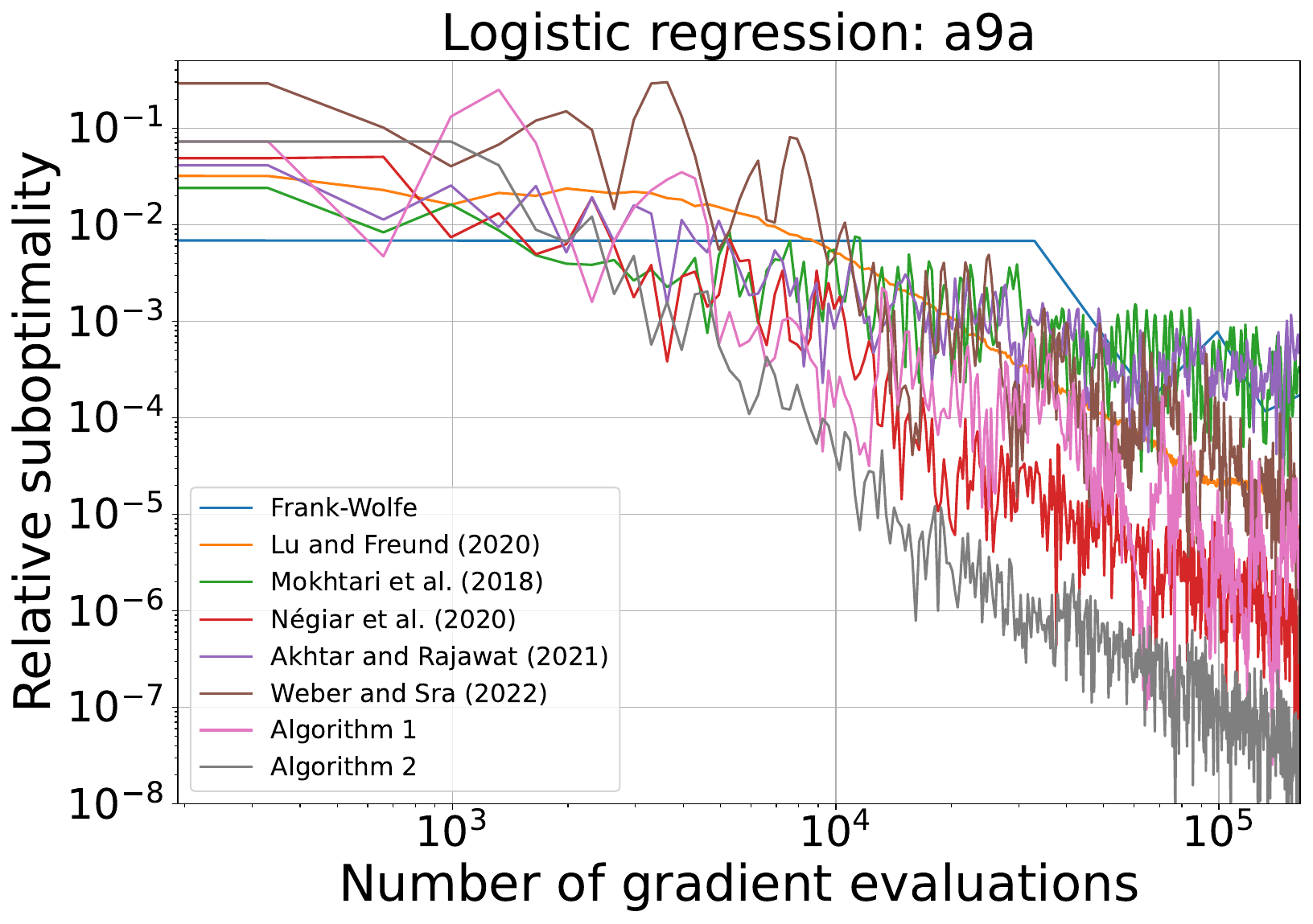}
\end{minipage}%
\begin{minipage}{0.24\textwidth}
  \centering
\includegraphics[width =  1\textwidth ]{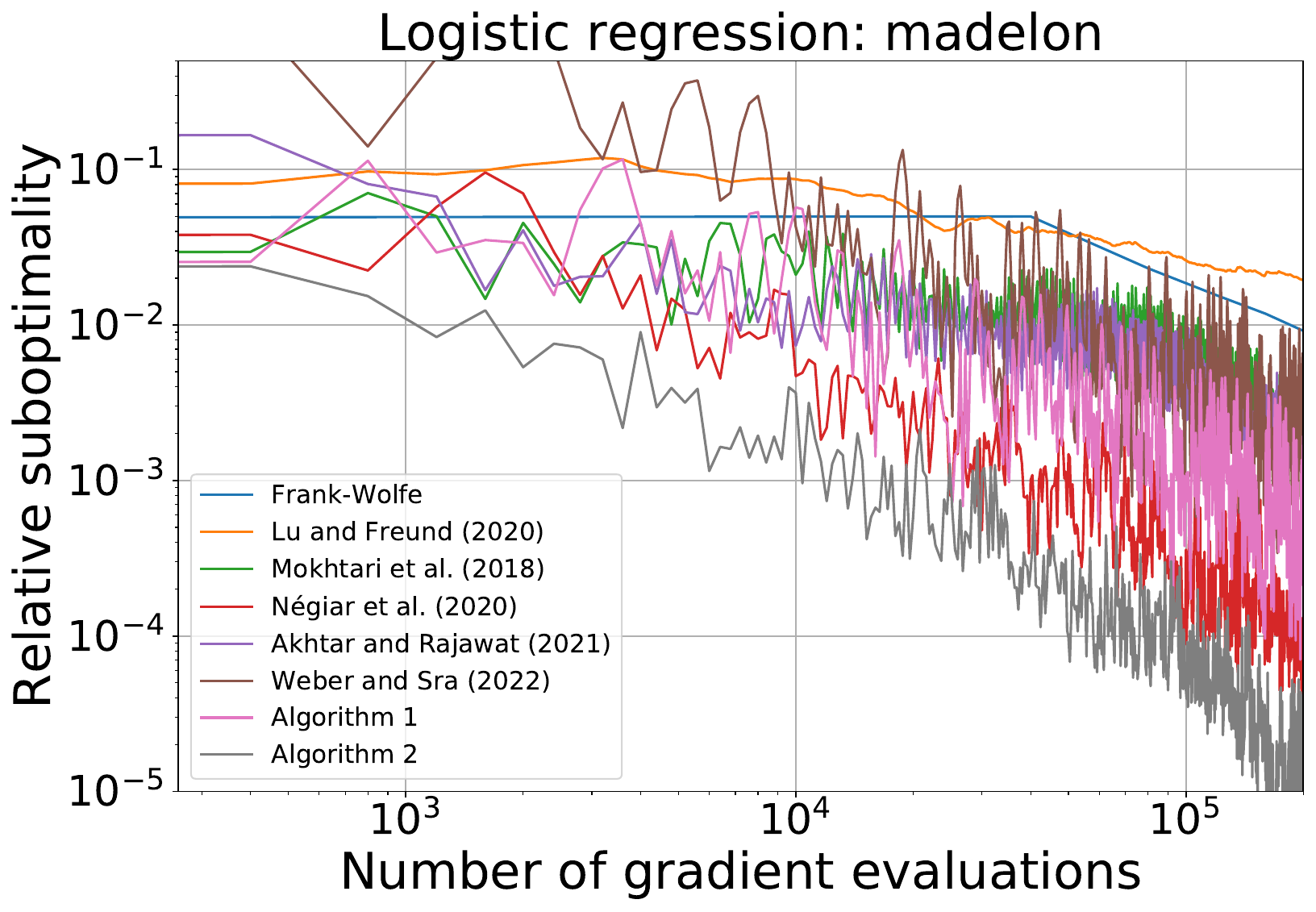}
\end{minipage}%
\begin{minipage}{0.24\textwidth}
  \centering
\includegraphics[width =  1\textwidth ]{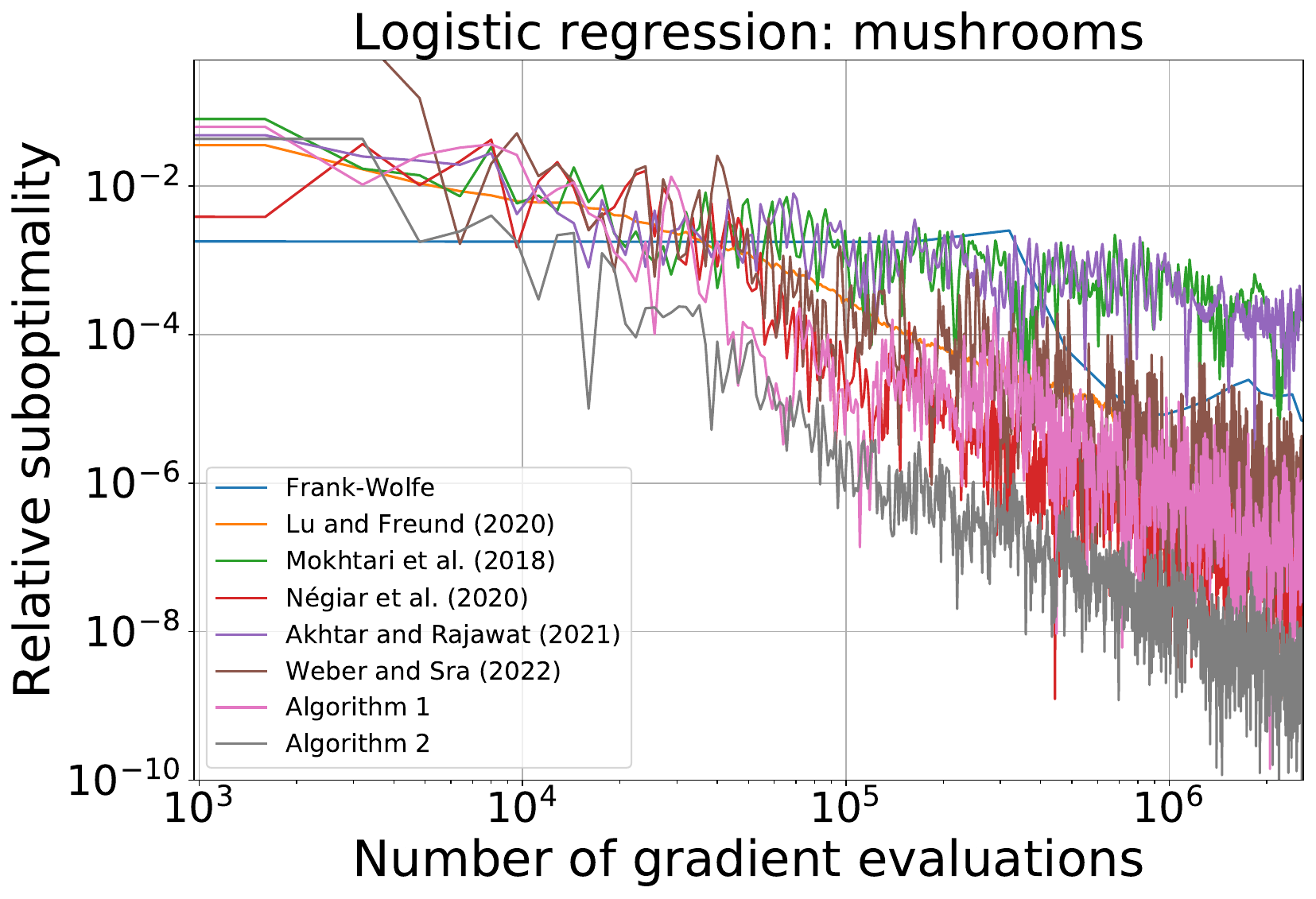}
\end{minipage}%
\begin{minipage}{0.24\textwidth}
  \centering
\includegraphics[width =  1\textwidth ]{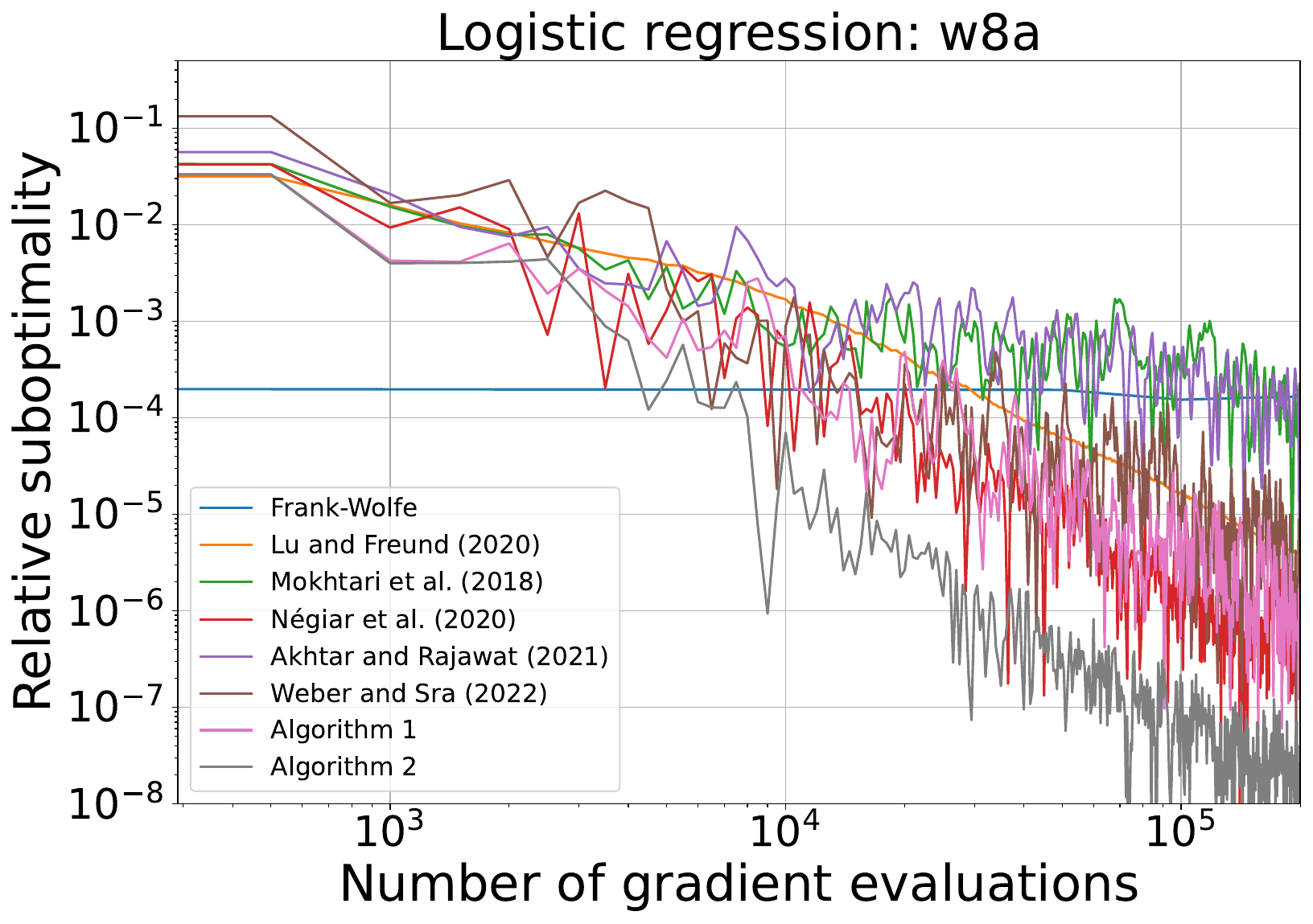}
\end{minipage}%
\\
\begin{minipage}{0.24\textwidth}
  \centering
(a) \texttt{a9a}
\end{minipage}%
\begin{minipage}{0.24\textwidth}
\centering
(b) \texttt{breast cancer}
\end{minipage}%
\begin{minipage}{0.24\textwidth}
\centering
  (c) \texttt{mushrooms}
\end{minipage}%
\begin{minipage}{0.24\textwidth}
\centering
  (d) \texttt{w8a}
\end{minipage}%
\vspace{-0.3cm}
\caption{Comparison of state-of-the-art projection free methods with small batches for \eqref{eq;ls} with $R = 2$. The comparison is made on the real datasets from LibSVM. The criterion is the number of full gradients computations.}
\end{figure*}

\end{document}